\documentclass[12pt,a4paper,onecolumn]{article}
\usepackage{textcomp}
\usepackage{tikz}
\usepackage{amssymb}
\usepackage{caption}
\usepackage{color}
\usepackage{verbatim}
\usepackage{fancybox}
\usepackage{titlesec} 
\renewcommand\thesection{\arabic{section}} 
\usepackage{esvect}
\usepackage{cite}
\usepackage{setspace}
\usepackage{mathtools}
\usepackage{amsmath}
\usepackage{mathrsfs}
\usetikzlibrary{positioning,fit,calc}
\usepackage[a4paper,left=1in,top=1in,right=1in,bottom=1in,nohead]{geometry} 
\usepackage[hidelinks]{hyperref}
\usepackage{array}
\usepackage{bm}
\usepackage{multirow}
\usepackage{hhline}

\usepackage{amsthm}
\usepackage[symbol]{footmisc} 
\usepackage{subcaption} 
\usetikzlibrary{decorations.pathmorphing}
\usetikzlibrary{arrows}
\numberwithin{equation}{section}
\usepackage{rotating}
\usepackage{ragged2e}
\DeclareFontFamily{U}{mathb}{}
\DeclareFontShape{U}{mathb}{m}{n}{<-5.5> mathb5 <5.5-6.5> mathb6 
<6.5-7.5> mathb7 <7.5-8.5> mathb8 <8.5-9.5> mathb9 <9.5-11> mathb10 
<11-> mathb12}{}
\DeclareSymbolFont{mathb}{U}{mathb}{m}{n}
\DeclareFontSubstitution{U}{mathb}{m}{n}
\DeclareMathSymbol{\blackdiamond}{\mathbin}{mathb}{"0C}
\tikzset{E/.style = {rectangle, fill=none, draw=#1}, E/.default=black}
\usetikzlibrary{positioning, shapes.geometric}

\def\correspondingauthor{\footnote{Corresponding author. Email: williewong088@gmail.com.}}

\tikzset{block/.style={draw,thick,text width=2cm,minimum height=1cm,align=center},
         line/.style={-latex}}
\newcolumntype{P}[1]{>{\centering\arraybackslash}m{#1}} 

\titleformat{\section}[block]{\large\scshape\bfseries}{\thesection.}{1em}{} 
\titleformat{\subsection}[block]{\bfseries}{\thesubsection.}{1em}{} 

\newtheorem{thm}{Theorem}[section]
\newtheorem{ppn}[thm]{Proposition}
\newtheorem{cor}[thm]{Corollary}
\newtheorem{lem}[thm]{Lemma}

\theoremstyle{definition}
\newtheorem{defn}[thm]{Definition}

\newtheorem{rmk}[thm]{Remark}
\newtheorem{eg}[thm]{Example}

\newtheorem{prob}[thm]{Problem}

\newtheorem{customcon}{Conjecture}
\renewcommand{\thecustomcon}{\Alph{customcon}} 

\begin{document}
\pagenumbering{arabic}
\begin{center}
    \textbf{\Large Optimal orientations of vertex-multiplications\\of Cartesian products of graphs}
\vspace{0.1 in} 
    \\{\large H.W. Willie Wong\correspondingauthor{}, E.G. Tay}
\vspace{0.1 in} 
\\National Institute of Education\\Nanyang Technological University\\Singapore
\end{center}
\begin{abstract}
\noindent Koh and Tay proved a fundamental classification of $G$ vertex-multiplications into three classes $\mathscr{C}_0, \mathscr{C}_1$ and $\mathscr{C}_2$. In this paper, we prove that vertex-multiplications of Cartesian products of graphs $G\times H$ lie in $\mathscr{C}_0$ ($\mathscr{C}_0\cup \mathscr{C}_1$ resp.) if $G^{(2)}\in \mathscr{C}_0$ ($\mathscr{C}_1$ resp.), $d(G)\ge 2$ and $d(G\times H)\ge 4$, providing further support for a conjecture by Koh and Tay. We also focus on Cartesian products involving trees, paths and cycles and show that most of them lie in $\mathscr{C}_0$.
\end{abstract}
\section{Introduction}
Let $G$ be a graph with vertex set $V(G)$ and edge set $E(G)$. In this paper, we consider only graphs with no loops or parallel edges. For any vertices $v,x\in V(G)$, the \textit{distance} from $v$ to $x$, $d_G(v,x)$, is defined as the length of a shortest path from $v$ to $x$. For $v\in V(G)$, its \textit{eccentricity} $e_G(v)$ is defined as $e_G(v)=\max\{d_G(v,x)\mid x\in V(G)\}$. A vertex $x$ is called an \textit{eccentric vertex} of $v$ if $d_G(v,x)=e_G(v)$. The \textit{diameter} of $G$, denoted by $d(G)$, is defined as $d(G)=\max\{e_G(v)\mid v\in V(G)\}$ while the \textit{radius} of $G$, denoted by $r(G)$, is defined as $r(G)=\min\{e_G(v)\mid v\in V(G)\}$. The above notions are defined similarly for a digraph $D$; and we refer the reader to \cite{BJJ GG} for any undefined terminology. For a digraph $D$, a vertex $x$ is said to be \textit{reachable} from another vertex $v$ if $d_D(v,x)<\infty$. The \textit{outset} and \textit{inset} of a vertex $v\in V(D)$ are defined to be $O_D(v)=\{x\in V(D)\mid v\rightarrow x\}$ and $I_D(v)=\{y\in V(D)\mid y\rightarrow v\}$ respectively. If there is no ambiguity, we shall omit the subscript for the above notation.
\indent\par An $\textit{orientation}$ $D$ of a graph $G$ is a digraph obtained from $G$ by assigning a direction to every edge $e\in E(G)$. An orientation $D$ of $G$ is said to be \textit{strong} if every two vertices in $V(D)$ are mutually reachable. An edge $e\in E(G)$ is a \textit{bridge} if $G-e$ is disconnected. Robbins' well-known One-way Street Theorem  \cite{RHE} states that a connected graph $G$ has a strong orientation if and only if $G$ is bridgeless.
\indent\par Given a connected and bridgeless graph $G$, let $\mathscr{D}(G)$ be the family of strong orientations of $G$. The $\textit{orientation number}$ of $G$ is defined as 
\begin{align*}
\bar{d}(G)=\min\{d(D)\mid D\in \mathscr{D}(G)\}.
\end{align*}
Any orientation $D$ in $\mathscr{D}(G)$ with $d(D)=\bar{d}(G)$ is called an \textit{optimal orientation} of $G$. The general problem of finding the orientation number of a connected and bridgeless graph is very difficult. Moreover, Chv{\'a}tal and Thomassen \cite{CV TC} proved that it is NP-hard to determine whether a graph admits an orientation of diameter 2. Hence, it is natural to focus on special classes of graphs. The orientation number was evaluated for various classes of graphs, such as the complete graphs \cite{BF TR,MSB,PJ 1} and complete bipartite graphs \cite{GG 1,SL}. 
\indent\par In 2000, Koh and Tay \cite{KKM TEG 8} introduced a new family of graphs, $G$ vertex-multiplications, and extended the results on the orientation number of complete $n$-partite graphs. Let $G$ be a given connected graph with vertex set $V(G)=\{v_1,v_2,\ldots, v_n\}$. For any sequence of $n$ positive integers $(s_i)$, a $G$ \textit{vertex-multiplication}, denoted by $G(s_1, s_2,\ldots, s_n)$, is the graph with vertex set $V^*=\bigcup_{i=1}^n{V_i}$ and edge set $E^*$, where $V_i$'s are pairwise disjoint sets with $|V_i|=s_i$, for $i=1,2,\ldots,n$; and for any $u,v\in V^*$, $uv\in E^*$ if and only if $u\in V_i$ and $v\in V_j$ for some $i,j\in \{1,2,\ldots, n\}$ with $i\neq j$ such that $v_i v_j\in E(G)$. For instance, if $G\cong K_n$, then the graph $G(s_1, s_2,\ldots, s_n)$ is a complete $n$-partite graph with partite sizes $s_1, s_2,\ldots, s_n$. Also, we say $G$ is a \textit{parent graph} of a graph $H$ if $H\cong G(s_1, s_2,\ldots, s_n)$ for some sequence $(s_i)$ of positive integers.
\indent\par For $i=1,2,\ldots, n$, we denote the $x$-th vertex in $V_i$ by $(x,v_i)$, i.e. $V_i=\{(x,v_i)\mid x=1,2,\ldots,s_i\}$. Hence, two vertices $(x,v_i)$ and $(y,v_j)$ in $V^*$ are adjacent in $G(s_1,s_2,\ldots, s_n)$ if and only if $i\neq j$ and $v_i v_j\in E(G)$. We will loosely use the two denotations of a vertex, for example, if $v_i=j$, then $v_i=v_j$ and $s_i=s_{j}$. For convenience, we write $G^{(s)}$ in place of $G(s,s,\ldots,s)$ for any positive integer $s$, and it is understood that the number of $s$'s is equal to the order $n$ of $G$. Thus, $G^{(1)}$ is simply the graph $G$ itself.
\indent\par $G$ vertex-multiplications are a natural generalisation of complete multipartite graphs. Optimal orientations minimising the diameter can also be used to solve a variant of the Gossip Problem on a graph $G$. The Gossip Problem attributed to Boyd by Hajnal et al. \cite{HA MEC SE} is stated as follows:
\begin{justify}
``There are $n$ ladies, and each one of them knows an item of scandal which is not known to any of the others. They communicate by telephone, and whenever two ladies make a call, they pass on to each other, as much scandal as they know at that time.  How many calls are needed before all ladies know all the scandal?"
\end{justify}
\indent\par The Problem has been the source of many papers that have studied the spread of information by telephone calls, conference calls, letters and computer networks. One can imagine a network of people modelled by a $G$ vertex-multiplication where the underlying graph is $G$ and persons within a partite set are not allowed to communicate directly with each other, for perhaps secrecy or disease containment reasons.
\indent\par The following theorem by Koh and Tay \cite{KKM TEG 8} provides a fundamental classification on $G$ vertex-multiplications.
\begin{thm} (Koh and Tay \cite{KKM TEG 8}) \label{thmA8.1.1} Let $G$ be a connected graph of order $n\ge 3$. If $s_i\ge 2$ for $i=1,2,\ldots, n$, then $d(G)\le \bar{d}(G(s_1,s_2,\ldots,s_n))\le d(G)+2$.
\end{thm}
In view of Theorem \ref{thmA8.1.1}, all graphs of the form $G(s_1,s_2,\ldots, s_n)$, with $s_i\ge 2$ for all $i=1,2,\ldots, n$, can be classified into three classes $\mathscr{C}_j$, where 
\begin{align*}
\mathscr{C}_j=\{G(s_1,s_2,\ldots, s_n)\mid \bar{d}(G(s_1,s_2,\ldots, s_n))=d(G)+j\},
\end{align*}
for $j=0,1,2$. Henceforth, we assume $s_i\ge 2$ for $i=1,2,\ldots, n$. The following lemma was found useful in proving Theorem \ref{thmA8.1.1}.

\begin{lem}(Koh and Tay \cite{KKM TEG 8}) \label{lemA8.1.2} Let $s_i,t_i$ be integers such that $s_i\le t_i$ for $i=1,2,\ldots, n$. If the graph $G(s_1,s_2,\ldots, s_n)$ admits an orientation $F$ in which every vertex $v$ lies on a cycle of length not exceeding $m$, then $\bar{d}(G(t_1, t_2,\ldots, t_n))\le \max\{m, d(F)\}$.
\end{lem}
Koh and Tay \cite{KKM TEG 8} made the following conjecture and proved it for some families of graphs, including cycles.
\begin{customcon}(Koh and Tay \cite{KKM TEG 8}) \label{conA8.1.A} If $G$ is a graph such that $d(G)\ge 3$ and $s_i\ge2$ for $i=1,2,\ldots,n$, then $G(s_1,s_2,\ldots,s_n)\not\in \mathscr{C}_2$.
\end{customcon}
These results and conjecture were generalised to digraphs by Gutin et al. \cite{GG KKM TEG AY}. Ng and Koh \cite{NKL KKM} examined cycle vertex-multiplications and Koh and Tay \cite{KKM TEG 11} investigated tree vertex-multiplications. Since trees with diameter at most 2 are parent graphs of complete bipartite graphs and are completely solved, they considered trees of diameter at least 3 and proved Conjecture \ref{conA8.1.A} for trees.
\begin{thm}(Koh and Tay \cite{KKM TEG 11})\label{thmA8.1.3}
\\If $T$ is a tree of order $n$ and $3\le d(T)\le 5$, then $T(s_1,s_2,\ldots ,s_n)\in \mathscr{C}_0\cup \mathscr{C}_1$.
\end{thm}
\begin{thm}(Koh and Tay \cite{KKM TEG 11})\label{thmA8.1.4}
\\If $T$ is a tree of order $n$ and $d(T)\ge 6$, then $T(s_1,s_2, \ldots ,s_n)\in \mathscr{C}_0$.
\end{thm}
In \cite{WHW TEG 3A}, Wong and Tay proved a characterisation for vertex-multiplications of trees with diameter $5$ in $\mathscr{C}_0$ and $\mathscr{C}_1$. In \cite{WHW TEG 6A}, they almost completely characterised vertex-multiplications of trees with diameter $4$.
\indent\par In this paper, we examine vertex-multiplications of Cartesian products of graphs and provide further support for Conjecture \ref{conA8.1.A}. Our main approach is a nimble application of Lemma \ref{lemA8.1.2} via some elementary orientations (see Definition \ref{defnA8.2.1}) leveraging on the neat structure enjoyed by Cartesian products of graphs. We also focus on Cartesian products involving trees, paths and cycles and show that most of them lie in $\mathscr{C}_0$. The Cartesian product of two graphs $G$ and $H$ is denoted by $G \times H$, where $V(G\times H)=\{\langle u,x\rangle\mid u\in V(G), x\in V(H)\}$ and $E(G)=\{\langle u,x\rangle \langle v,y\rangle\mid u=v$ and $xy\in E(H)$, or $uv\in E(G)$ and $x=y\}$. Since Cartesian products of disconnected graphs are disconnected, we concern ourselves with only connected graphs. We shall denote a path (cycle, complete graph resp.) of order $n$ as $P_n$ ($C_n$, $K_n$ resp.) while $T_d$ represents a tree of diameter $d$. Since the orientation number of complete bipartite graphs $K(p,q)$ has been characterised by {\v S}olt{\'e}s \cite{SL} and Gutin \cite{GG 1} and $P_2\times P_2\cong K(2,2)$, we shall exclude $P_2\times P_2$ from our discussion. In Section 2, we consider Cartesian products of graphs in the general setting.
\begin{thm} \label{thmA8.1.5} Let $G$ and $H$ be connected graphs with order at least two. If $d(G)\ge 2$ and $G^{(2)} \in \mathscr{C}_0$ ($\mathscr{C}_1$ resp.), then $(G\times H)^{(2)} \in \mathscr{C}_0$ ($\mathscr{C}_0 \cup \mathscr{C}_1$ resp.).
\end{thm}
\begin{cor}\label{corA8.1.6} Let $G$ and $H$ be connected graphs with order at least two. If $d(G\times H)\ge 4$, $d(G)\ge 2$ and $G^{(2)} \in \mathscr{C}_0$ ($\mathscr{C}_1$ resp.), then $(G\times H)(s_1,s_2,\ldots,s_n) \in \mathscr{C}_0$ ($\mathscr{C}_0 \cup \mathscr{C}_1$ resp.).
\end{cor}
In Section 3, we prove that the vertex-multiplications of Cartesian products of two trees are mostly in $\mathscr{C}_0$.
\begin{thm}\label{thmA8.1.7}
If $\lambda\ge 2$ and $\mu \ge 3$, then $(T_\lambda\times T_\mu)(s_1,s_2,\ldots,s_n)\in\mathscr{C}_0$.
\end{thm}
For trees with diameter $2$, the same conclusion holds if both trees are paths, i.e., $P_3\times P_3$.
\begin{thm}\label{thmA8.1.8}
~\\(a) $(P_3\times P_2)^{(2)}\in\mathscr{C}_1$.
\\(b)
\begin{equation}
  (P_\lambda\times P_\mu)(s_1,s_2,\ldots,s_n)\in\left\{
  \begin{array}{@{}ll@{}}
    \mathscr{C}_0, & \text{if } \lambda\ge 4, \mu=2, \text{ or } \lambda\ge \mu\ge3,  \nonumber\\
    \mathscr{C}_0\cup\mathscr{C}_1, & \text{if } (\lambda,\mu)=(3,2). \nonumber\\
  \end{array}\right.
\end{equation}
\end{thm}
We also prove an analogue on the hypercube graph $Q_\lambda=\overbrace{K_2\times K_2\times \cdots\times K_2}^{\lambda}$, $\lambda\in\mathbb{Z}^+$.
\begin{ppn}\label{ppnA8.1.9}
~\\(a) $Q_3^{(2)}\in\mathscr{C}_0$ and $Q_3(s_1,s_2,\ldots,s_n)\in\mathscr{C}_0\cup\mathscr{C}_1$.
\\(b) For $\lambda\ge 4$, $Q_\lambda(s_1,s_2,\ldots,s_n)\in\mathscr{C}_0$.
\end{ppn}
In Sections 4 and 5, we examine the Cartesian products of a tree and a cycle, and two cycles respectively. 
\begin{thm} \label{thmA8.1.10}
If $\lambda\ge 2$ and $\mu\ge 4$ or $\lambda=\mu=3$, then $(T_\lambda\times C_\mu)(s_1,s_2,\ldots,s_n)\in\mathscr{C}_0$.
\end{thm}

\begin{thm}\label{thmA8.1.11}
~\\(a) For $\lambda\ge 4, \mu\ge 3$, $(C_\lambda\times C_\mu)^{(2)}\in\mathscr{C}_0$.
\\(b)
\begin{equation}
  (C_\lambda\times C_\mu)(s_1,s_2,\ldots,s_n)\in\left\{
  \begin{array}{@{}ll@{}}
    \mathscr{C}_0, & \text{if } \lambda\ge\mu\ge 4,  \nonumber\\
    \mathscr{C}_0\cup\mathscr{C}_1, & \text{if } (\lambda,\mu)=(3,3) \text{ or } (4,3). \nonumber\\
  \end{array}\right.
\end{equation}
\end{thm}
\section{General results}
In defining an orientation, we use the following notation to write succinctly. For any orientation $D$, $\tilde{D}$ denotes the orientation satisfying $u\rightarrow v$ in $\tilde{D} \iff v\rightarrow u$ in $D$. 
\begin{defn}\label{defnA8.2.1}
Suppose $uv$ and $wx$ are edges of a graph $G$ and $D$ is an orientation of $G^{(2)}$. We denote 
\\(a) $u\rightrightarrows v$ if $\{(1,u),(2,u)\}\rightarrow \{(1,v),(2,v)\}$ in $D$ (see Figure \ref{figA8.1.1}(a)).
\\(b) $u\rightsquigarrow v$ if $(1,u)\rightarrow (1,v)\rightarrow (2,u)\rightarrow (2,v)\rightarrow (1,u)$ in $D$ (see Figure \ref{figA8.1.1}(b)).
\\(c) $u\overset{1}\twoheadrightarrow v$ if $(2,v)\rightarrow \{(1,u), (2,u)\}\rightarrow (1,v)$ and $w\overset{2}\twoheadrightarrow x$ if $(1,x)\rightarrow \{(1,w), (2,w)\}\rightarrow (2,x)$ in $D$ (see Figure \ref{figA8.1.1}(c)).
\end{defn}
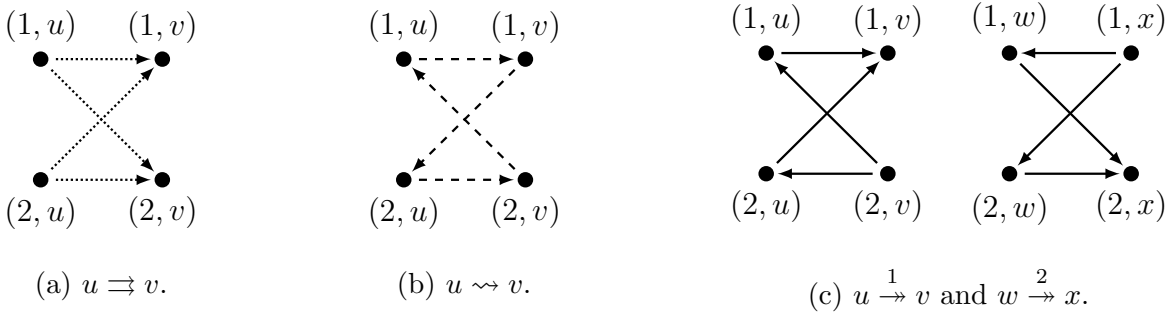
\begin{figure}[h]
\begin{subfigure}{.2\textwidth}
\begin{center}
\tikzstyle{every node}=[circle, draw, fill=black!100,
                       inner sep=0pt, minimum width=5pt]
\begin{tikzpicture}[thick,scale=0.8]%
\draw(-4,0)node[label={[yshift=0.2cm]270:{$(2,u)$}}]{};
\draw(-2,0)node[label={[yshift=0.2cm]270:{$(2,v)$}}]{};
\draw(-2,2)node[label={[yshift=-0.2cm]90:{$(1,v)$}}]{};
\draw(-4,2)node[label={[yshift=-0.2cm]90:{$(1,u)$}}]{};

\draw[densely dotted, ->, line width=0.3mm, >=latex, shorten <= 0.2cm, shorten >= 0.15cm](-4,0)--(-2,0);
\draw[densely dotted, ->, line width=0.3mm, >=latex, shorten <= 0.2cm, shorten >= 0.15cm](-4,2)--(-2,0);
\draw[densely dotted, ->, line width=0.3mm, >=latex, shorten <= 0.2cm, shorten >= 0.15cm](-4,2)--(-2,2);
\draw[densely dotted, ->, line width=0.3mm, >=latex, shorten <= 0.2cm, shorten >= 0.15cm](-4,0)--(-2,2);
\end{tikzpicture}
{\caption{$u\rightrightarrows v$.}}
\end{center}
\end{subfigure}\hfill%
\begin{subfigure}{.2\textwidth}
\begin{center}
\tikzstyle{every node}=[circle, draw, fill=black!100,
                       inner sep=0pt, minimum width=5pt]
\begin{tikzpicture}[thick,scale=0.8]%
\draw(-4,0)node[label={[yshift=0.2cm]270:{$(2,u)$}}]{};
\draw(-2,0)node[label={[yshift=0.2cm]270:{$(2,v)$}}]{};
\draw(-2,2)node[label={[yshift=-0.2cm]90:{$(1,v)$}}]{};
\draw(-4,2)node[label={[yshift=-0.2cm]90:{$(1,u)$}}]{};

\draw[dashed, ->, line width=0.3mm, >=latex, shorten <= 0.2cm, shorten >= 0.15cm](-4,0)--(-2,0);
\draw[dashed, ->, line width=0.3mm, >=latex, shorten <= 0.2cm, shorten >= 0.15cm](-2,0)--(-4,2);
\draw[dashed, ->, line width=0.3mm, >=latex, shorten <= 0.2cm, shorten >= 0.15cm](-4,2)--(-2,2);
\draw[dashed, ->, line width=0.3mm, >=latex, shorten <= 0.2cm, shorten >= 0.15cm](-2,2)--(-4,0);
\end{tikzpicture}
{\caption{$u\rightsquigarrow v$.}}
\end{center}
\end{subfigure}\hfill
\begin{subfigure}{.4\textwidth}
\begin{center}
\tikzstyle{every node}=[circle, draw, fill=black!100,
                       inner sep=0pt, minimum width=5pt]
\begin{tikzpicture}[thick,scale=0.8]%
\draw(-4,0)node[label={[yshift=0.2cm] 270:{$(2,u)$}}]{};
\draw(-2,0)node[label={[yshift=0.2cm] 270:{$(2,v)$}}]{};
\draw(-4,2)node[label={[yshift=-0.2cm] 90:{$(1,u)$}}]{};
\draw(-2,2)node[label={[yshift=-0.2cm] 90:{$(1,v)$}}]{};

\draw(0,0)node[label={[yshift=0.2cm] 270:{$(2,w)$}}]{};
\draw(2,0)node[label={[yshift=0.2cm] 270:{$(2,x)$}}]{};
\draw(0,2)node[label={[yshift=-0.2cm] 90:{$(1,w)$}}]{};
\draw(2,2)node[label={[yshift=-0.2cm] 90:{$(1,x)$}}]{};

\draw[->, line width=0.3mm, >=latex, shorten <= 0.2cm, shorten >= 0.15cm](-4,0)--(-2,2);
\draw[->, line width=0.3mm, >=latex, shorten <= 0.2cm, shorten >= 0.15cm](-4,2)--(-2,2);
\draw[->, line width=0.3mm, >=latex, shorten <= 0.2cm, shorten >= 0.15cm](-2,0)--(-4,2);
\draw[->, line width=0.3mm, >=latex, shorten <= 0.2cm, shorten >= 0.15cm](-2,0)--(-4,0);

\draw[->, line width=0.3mm, >=latex, shorten <= 0.2cm, shorten >= 0.15cm](0,0)--(2,0);
\draw[->, line width=0.3mm, >=latex, shorten <= 0.2cm, shorten >= 0.15cm](0,2)--(2,0);
\draw[->, line width=0.3mm, >=latex, shorten <= 0.2cm, shorten >= 0.15cm](2,2)--(0,2);
\draw[->, line width=0.3mm, >=latex, shorten <= 0.2cm, shorten >= 0.15cm](2,2)--(0,0);
\end{tikzpicture}
{\captionof{figure}{$u\overset{1}\twoheadrightarrow v$ and $w\overset{2}\twoheadrightarrow x$.}}
\end{center}
\end{subfigure}
{\caption{Notation for orientations.}\label{figA8.1.1}}
\end{figure}

\noindent\textit{Proof of Theorem \ref{thmA8.1.5}}: Since $G^{(2)} \in \mathscr{C}_0$, there exists an orientation $D$ of $G^{(2)}$ such that $d(D)=d(G)$. Define an orientation $D^*$ of $(G\times H)^{(2)}$ as follows.
\\For any $u,v\in V(G)$, any $x\in V(H)$ and any $p,q=1,2$,
\begin{align}
(p,\langle u,x\rangle)\rightarrow (q,\langle v,x\rangle)\iff (p,u)\rightarrow (q,v) \text{ in }D,\label{eqA8.2.1}
\end{align}
i.e. each copy of $G^{(2)}$ is oriented similarly to $D$.
\\For any $u\in V(G)$ and any $x,y\in V(H)$,
\begin{align}
\langle u,x\rangle\rightsquigarrow \langle u,y\rangle \iff xy\in E(H). \label{eqA8.2.2}
\end{align}
\noindent\par We remark that the definition in (\ref{eqA8.2.2}) is arbitrary since $xy\in E(H)$ is equivalent to $yx\in E(H)$. However, this shall not affect the following argument. For a well-defined orientation, one may linearly order the vertices in $V(H)$ before applying (\ref{eqA8.2.2}) with the condition $x$ precedes $y$.
\noindent\par We claim that $d_{D^*}((p,\langle u,x\rangle),(q,\langle v,y\rangle))\le d(G\times H)=d(G)+d(H)$ for $p,q=1,2$. If $x=y$, then by (\ref{eqA8.2.1}), $d_{D^*}((p,\langle u,x\rangle),(q,\langle v,x\rangle))\le d(D)=d(G)< d(G)+d(H)$.
\noindent\par Suppose $x\neq y$. In view of (\ref{eqA8.2.2}), there exists a $(p,\langle u,x\rangle)-(r,\langle u,y\rangle)$ path of length $d_{D^*}((p,\langle u,x\rangle),(r,\langle u,y\rangle))\le d_H(x,y)\le d(H)$ for some $r=1,2$. If $(q,\langle v,y\rangle)=(r,\langle u,y\rangle)$ (i.e. $q=r$, $u=v$), then we are done. If $(q,\langle v,y\rangle)= (3-r,\langle u,y\rangle)$, then $d_{D^*}((p,\langle u,x\rangle),(q,\langle v,y\rangle))\le d_{D^*}((p,\langle u,x\rangle),(r,\langle u,y\rangle))+2\le d(H)+2\le d(H)+d(G)$ by (\ref{eqA8.2.2}). Finally, if $u\neq v$, then
$d_{D^*}((p,\langle u,x\rangle),(q,\langle v,y\rangle))\le d_{D^*}((p,\langle u,x\rangle),(r,\langle u,y\rangle))+d_{D^*}((r,\langle u,y\rangle),(q,\langle v,y\rangle))\le d(H)+d(D)=d(H)+d(G)$.
\indent\par The proof is similar if $G^{(2)}\in\mathscr{C}_1$.
\qed
\\
\\\textit{Proof of Corollary \ref{corA8.1.6}}: Since every vertex lies in a directed $C_4$ in the orientation $D^*$ because of (\ref{eqA8.2.2}), it follows from Lemma \ref{lemA8.1.2} that $(G\times H)(s_1,s_2,\ldots, s_n) \in \mathscr{C}_0$. The proof is similar if $G^{(2)}\in\mathscr{C}_1$.
\qed

\begin{cor}\label{corA8.2.2} For all $i\in \mathbb{Z}^+$, let $G_i$ be a connected graph with order at least two. If $(G_1\times G_2)^{(2)}\in \mathscr{C}_0$ ($\mathscr{C}_1$ resp.), then
~\\(a) for $j\ge 3$, $(\prod\limits_{i=1}^{j}G_i)^{(2)}\in\mathscr{C}_0$ ($\mathscr{C}_0\cup \mathscr{C}_1$ resp.), and 
\\(b) for $k\ge 4$, $(\prod\limits_{i=1}^{k}G_i)(s_1,s_2,\ldots,s_n)\in\mathscr{C}_0$ ($\mathscr{C}_0\cup \mathscr{C}_1$ resp.).
\end{cor}
\noindent\textit{Proof}: (a) Since $d(G_1\times G_2)\ge 2$, the result follows from Theorem \ref{thmA8.1.5}.
\\(b) Since $d(\prod\limits_{i=1}^{k}G_i)\ge 4$, $d(G_1\times G_2)\ge 2$, the result follows from Corollary \ref{corA8.1.6}.
\qed

\begin{cor}\label{corA8.2.3} Let $G$ be a connected graph with order at least two.
~\\(a) If $3\le d\le 5$, then $(T_d\times G)(s_1,s_2,\ldots,s_n)\in\mathscr{C}_0\cup\mathscr{C}_1$.
\\(b) If $d\ge 6$, then $(T_d\times G)(s_1,s_2,\ldots,s_n)\in\mathscr{C}_0$.
\end{cor}
\noindent\textit{Proof}: Since $d(T_d\times G)\ge 4$ and by Corollary \ref{corA8.1.6}, (a) and (b) follow from Theorems \ref{thmA8.1.3} and \ref{thmA8.1.4} respectively.
\qed
\section{Cartesian product of trees $T_\lambda\times T_\mu$}
In this section, we shall show that Corollary \ref{corA8.2.3}(a) can be further improved in the case of $T_\lambda\times T_\mu$. Before that, we introduce a notation for trees $T_d$ with $d\le 5$. Whenever we speak of a tree with even diameter $d$, we denote by $\mathtt{c}$, the unique central vertex of $T_d$, i.e. $e_{T_d}(\mathtt{c})=r(T_d)$, and the neighbours of $\mathtt{c}$ by $[i]$, i.e. $N_{T_d}(\mathtt{c})=\{[i]\mid i=1,2,\ldots, deg_{T_d}(\mathtt{c}) \}$. For each $i =1,2,\ldots, deg_{T_d}(\mathtt{c})$, we further denote the neighbours of $[i]$, excluding $\mathtt{c}$, by $[\alpha, i]$, i.e. $N_{T_d}([i])-\{\mathtt{c}\}=\{[\alpha, i]\mid \alpha=1,2,\ldots, deg_{T_d}([i])-1\}$.  
\indent\par If $d$ is odd, we let $\mathtt{c}_1$ and $\mathtt{c}_2$ be the two central vertices of $T_d$, i.e. $e_{T_d}(\mathtt{c}_k)=r(T_d)$ for $k=1,2$. For $k=1,2$, denote the neighbours of $\mathtt{c}_k$, excluding $\mathtt{c}_{3-k}$, by $[i]_k$. i.e. $N_{T_d}(\mathtt{c}_k)-\{\mathtt{c}_{3-k}\}=\{[i]_k\mid i=1,2,\ldots, deg_{T_d}(\mathtt{c}_k)-1 \}$. For each $i =1,2,\ldots, deg_{T_d}(\mathtt{c}_k)-1$, we denote the neighbours of $[i]_k$, excluding $\mathtt{c}_k$, by $[\alpha, i]_k$. i.e. $N_{T_d}([i]_k)-\{\mathtt{c}_k\}=\{[\alpha, i]_k\mid \alpha=1,2,\ldots, deg_{T_d}([i]_k)-1\}$. Figures \ref{figA8.3.2} and \ref{figA8.3.3} illustrate the use of this notation.
\begin{center}
\tikzstyle{every node}=[circle, draw, fill=black!100,
                       inner sep=0pt, minimum width=6pt]
\begin{tikzpicture}[thick,scale=0.7]%
\draw(0,0)node[label={[yshift=-0.2cm]270:{$\mathtt{c}$}}](c){};

\draw(-6,2)node[label={[yshift=0.2cm] 270:{$[1,1]$}}](11){};
\draw(-6,0)node[label={[yshift=0.2cm] 270:{$[2,1]$}}](12){};
\draw(-6,-2)node[label={[yshift=0.2cm] 270:{$[1,2]$}}](21){};

\draw(-3,1)node[label={[yshift=-0.2cm] 270:{$[1]$}}](1){};
\draw(-3,-1)node[label={[yshift=-0.2cm] 270:{$[2]$}}](2){};
\draw(3,1)node[label={[yshift=-0.2cm] 270:{$[3]$}}](3){};
\draw(3,-1)node[label={[yshift=-0.2cm] 270:{$[4]$}}](4){};

\draw(6,-2)node[label={[yshift=0.2cm] 270:{$[1,4]$}}](41){};

\draw(c)--(1);
\draw(c)--(2);
\draw(c)--(3);
\draw(c)--(4);

\draw(1)--(11);
\draw(1)--(12);
\draw(2)--(21);
\draw(4)--(41);
\end{tikzpicture}
{\captionof{figure}{Labelling vertices in a $T_4$}\label{figA8.3.2}}

\begin{tikzpicture}[thick,scale=0.7]%

\draw(-6,2)node[label={[yshift=0.2cm] 270:{$[1,1]_{1}$}}](11_c1){};
\draw(-6,0)node[label={[yshift=0.2cm] 270:{$[2,1]_{1}$}}](12_c1){};
\draw(-6,-2)node[label={[yshift=0.2cm] 270:{$[1,2]_{1}$}}](21_c1){};

\draw(-3,1)node[label={[yshift=0cm] 270:{$[1]_{1}$}}](1_c1){};
\draw(-3,-1)node[label={[yshift=0cm] 270:{$[2]_{1}$}}](2_c1){};
\draw(0,0)node[label={[yshift=-0.2cm]270:{$\mathtt{c}_1$}}](c1){};

\draw(3,0)node[label={[yshift=-0.2cm] 270:{$\mathtt{c}_2$}}](c2){};
\draw(6,1)node[label={[yshift=0cm] 270:{$[1]_{2}$}}](1_c2){};

\draw(6,-1)node[label={[yshift=0cm] 270:{$[2]_{2}$}}](2_c2){};
\draw(9,1)node[label={[yshift=0.1cm] 270:{$[1,2]_{2}$}}](12_c2){};
\draw(9,-1)node[label={[yshift=0.2cm] 270:{$[2,2]_{2}$}}](22_c2){};
\draw(9,-3)node[label={[yshift=0.2cm] 270:{$[3,2]_{2}$}}](32_c2){};

\draw(c1)--(c2);

\draw(c1)--(1_c1);
\draw(c1)--(2_c1);
\draw(1_c1)--(11_c1);
\draw(1_c1)--(12_c1);
\draw(2_c1)--(21_c1);

\draw(c2)--(1_c2);
\draw(c2)--(2_c2);
\draw(2_c2)--(12_c2);
\draw(2_c2)--(22_c2);
\draw(2_c2)--(32_c2);
\end{tikzpicture}
{\captionof{figure}{Labelling vertices in a $T_5$}\label{figA8.3.3}}
\end{center}

With this, we prove Theorem \ref{thmA8.1.7}.
\\
\\\textit{Proof of Theorem \ref{thmA8.1.7}}:
Let $G=T_\lambda\times T_\mu$. By Corollary \ref{corA8.2.3}(b), it suffices to consider $\lambda,\mu\le 5$. Define an orientation $D_{(\lambda,\mu)}$ for $G^{(2)}$ as follows.
\\
\\Case 1. $\lambda$ is even and $\mu$ is odd, i.e. $\lambda=2,4$ and $\mu=3,5$.
\begin{align}
&\langle\mathtt{c},\mathtt{c}_2\rangle 
\rightrightarrows \langle[y],\mathtt{c}_2\rangle
\rightrightarrows \langle[y],\mathtt{c}_1\rangle
\rightrightarrows \langle\mathtt{c},\mathtt{c}_1\rangle
\rightrightarrows \langle\mathtt{c},\mathtt{c}_2\rangle\label{eqA8.3.1}
\end{align} 
for all $[y]\in N_{T_\lambda}(\mathtt{c})$. Excluding the edges defined above, for each $[i]_1\in N_{T_\mu}(\mathtt{c}_1)-\{\mathtt{c}_2\}$, each $\alpha=1,2,\ldots,$ $deg_{T_\mu}([i]_1)-1$, each $[j]_2\in N_{T_\mu}(\mathtt{c}_2)-\{\mathtt{c}_1\}$, and each $\beta=1,2,\ldots,deg_{T_\mu}([j]_2)-1$,
\begin{align}
\langle x,[\alpha,i]_1\rangle \rightsquigarrow \langle x,[i]_1\rangle \rightsquigarrow \langle x,\mathtt{c}_1\rangle,
\langle x,[\beta,j]_2\rangle \rightsquigarrow \langle x,[j]_2\rangle \rightsquigarrow \langle x,\mathtt{c}_2\rangle\rightsquigarrow \langle x,\mathtt{c}_1\rangle, \label{eqA8.3.2}
\end{align}
for all $x\in V(T_\lambda)$, and

\begin{equation}
\left. \begin{array}{@{}ll@{}}
&\langle [\gamma,k],[\alpha,i]_1\rangle\rightsquigarrow \langle [k],[\alpha,i]_1\rangle \rightsquigarrow \langle \mathtt{c},[\alpha,i]_1\rangle,\\
&\langle [\gamma,k],[i]_1\rangle\rightsquigarrow \langle [k],[i]_1\rangle \rightsquigarrow \langle \mathtt{c},[i]_1\rangle,\\
&\langle [\gamma,k],[\beta,j]_2\rangle\rightsquigarrow\langle [k],[\beta,j]_2\rangle \rightsquigarrow \langle \mathtt{c},[\beta,j]_2\rangle,\\
&\langle [\gamma,k],[j]_2\rangle\rightsquigarrow\langle [k],[j]_2\rangle \rightsquigarrow \langle \mathtt{c},[j]_2\rangle,\\
&\langle [\gamma,k],\mathtt{c}_t\rangle\rightsquigarrow \langle [k],\mathtt{c}_t\rangle \text{ for } t=1,2,
  \end{array}\right\}
\label{eqA8.3.3}
\end{equation}
for all $[k]\in N_{T_\lambda}(\mathtt{c})$ and $\gamma=1,2,\ldots, deg_{T_\lambda}([k])-1$. (See Figures \ref{figA8.3.4}-\ref{figA8.3.7}.)
\indent\par We claim that $d(D_{(\lambda,\mu)})=\lambda+\mu=d(G)$. Let $u,v\in V(G)$ and $P:=w_0 w_1 \ldots w_l$ be a shortest $u-v$ path in $G$ with $u=w_0$ and $v=w_l$. If $d_G(u,v)\le d(G)-2$ and $P$ satisfies 
\begin{align}
w_i\rightsquigarrow w_{i+1}\text{ or }w_{i+1}\rightsquigarrow w_i\text{ for all }i=0,1,\ldots, l-1,\label{eqA8.3.4}
\end{align}
then $d_{D_{(\lambda,\mu)}}((p,u),(q,v))\le d_G(u,v)+2\le d(G)$ for $p,q=1,2$. Particularly, this holds for $u=\langle [\gamma_1,k],y_1\rangle$, $v=\langle[ \gamma_2,k],y_2\rangle$ with $\gamma_1\neq \gamma_2$ in $T_4\times T_\mu$. So, by symmetry of (\ref{eqA8.3.1})-(\ref{eqA8.3.3}), we may assume without loss of generality that $\mathtt{c}$ has two eccentric vertices in $T_\lambda$, i.e. $T_\lambda=P_3$ if $\lambda=2$, and $T_\lambda=P_5$ if $\lambda=4$. Furthermore, by symmetry of (\ref{eqA8.3.2}), we may assume $\mathtt{c}_i$ has two eccentric vertices for $i=1,2$, in $T_\mu$.
\indent\par For the pairs of $u,v$ that do not satisfy (\ref{eqA8.3.4}), we claim that there exists a path $P$ with length at most $d(G)$ that satisfies 
\begin{equation}
\left.\begin{array}{@{}ll@{}}
&w_i\rightrightarrows w_{i+1}\text{ for some }i=0,1,\ldots, l-1\text{ and }\\
&w_{j+1}\rightrightarrows w_j\text{ for none of }j=0,1,\ldots, l-1.
  \end{array}\right\}
\label{eqA8.3.5}
\end{equation}
Then, we can conclude $d_{D_{(\lambda,\mu)}}((p,u),(q,v))\le d(G)$.
\\
\\Subcase 1.1. $\lambda=2$ and $\mu=3$. (See Figure \ref{figA8.3.4}.)
\indent\par We list these paths $P$ while omitting symmetric scenarios. For $i=1,2,$ and $j=1,2$,
\begin{align*}
P^1=&\langle [1],[1]_1\rangle \langle [1],\mathtt{c}_1\rangle \langle \mathtt{c},\mathtt{c}_1\rangle \langle \mathtt{c},[j]_1\rangle \langle [2],[j]_1\rangle.\\
P^2=&\langle [1],[1]_1\rangle \langle [1],\mathtt{c}_1\rangle \langle \mathtt{c},\mathtt{c}_1\rangle \langle \mathtt{c},\mathtt{c}_2\rangle \langle [i],\mathtt{c}_2\rangle \langle [i],[j]_2\rangle.\\
P^3=&\langle [1],[1]_1\rangle \langle [1],\mathtt{c}_1\rangle \langle \mathtt{c},\mathtt{c}_1\rangle \langle \mathtt{c},\mathtt{c}_2\rangle \langle \mathtt{c},[j]_2\rangle.\\
P^4=&\langle [1],[1]_2\rangle \langle [1],\mathtt{c}_2\rangle \langle [1],\mathtt{c}_1\rangle \langle [1],[j]_1\rangle \langle\mathtt{c},[j]_1\rangle \langle [2],[j]_1\rangle.\\
P^5=&\langle [1],[1]_2\rangle \langle\mathtt{c},[1]_2\rangle \langle\mathtt{c},\mathtt{c}_2\rangle \langle[2],\mathtt{c}_2\rangle \langle[2],\mathtt{c}_1\rangle \langle\mathtt{c},\mathtt{c}_1\rangle.\\
P^6=&\langle [1],[1]_2\rangle \langle\mathtt{c},[1]_2\rangle \langle\mathtt{c},\mathtt{c}_2\rangle \langle[2],\mathtt{c}_2\rangle \langle[2],[j]_2\rangle.\\
P^7=&\langle \mathtt{c},[1]_1\rangle \langle \mathtt{c},\mathtt{c}_1\rangle \langle\mathtt{c},\mathtt{c}_2\rangle \langle [i],\mathtt{c}_2\rangle \langle [i],[j]_2\rangle \langle \mathtt{c},[j]_2\rangle.\\
P^8=&\langle \mathtt{c},[1]_2\rangle \langle \mathtt{c},\mathtt{c}_2\rangle \langle [i],\mathtt{c}_2\rangle \langle [i],\mathtt{c}_1\rangle \langle [i],[j]_1\rangle \langle\mathtt{c},[j]_1\rangle.\\
P^9=&\langle \mathtt{c},[1]_2\rangle \langle \mathtt{c},\mathtt{c}_2\rangle \langle [i],\mathtt{c}_2\rangle \langle [i],\mathtt{c}_1\rangle \langle \mathtt{c},\mathtt{c}_1\rangle.
\end{align*}
Subcase 1.2. $\lambda=2$ and $\mu=5$. (See Figure \ref{figA8.3.5}.) 
\indent\par Note that $D_{(2,3)}$ is a subdigraph of $D_{(2,5)}$. Moveover, for any $(p,u)\in V(D_{(2,5)})-V(D_{(2,3)})$, there exists a vertex $(r,x)\in V(D_{(2,3)})$ such that $u\rightsquigarrow x$ or $x\rightsquigarrow u$. Hence, if $(p,u)\in V(D_{(2,5)})-V(D_{(2,3)})$ and $(q,v)\in V(D_{(2,3)})$, then $\max\{d_{D_{(2,5)}}((p,u),(q,v)),$ $d_{D_{(2,5)}} ((q,v),$ $(p,u))\}\le 1+d(D_{(2,3)})\le 7$ for $p,q=1,2$. Similarly, if $(p,u),(q,v)\in V(D_{(2,5)})-V(D_{(2,3)})$, then $d_{D_{(2,5)}}((p,u),(q,v))\le 2+d(D_{(2,3)})\le 7$ for $p,q=1,2$.
\\
\\Subcase 1.3. $\lambda=4$ and $\mu=3$. (See Figure \ref{figA8.3.6}.)
\indent\par Note that $D_{(2,3)}$ is a subdigraph of $D_{(4,3)}$ and this subcase follows by an argument similar to Subcase 1.2.
\\
\\Subcase 1.4. $\lambda=4$ and $\mu=5$. (See Figure \ref{figA8.3.7}.)
\indent\par Note that $D_{(4,3)}$ is a subdigraph of $D_{(4,5)}$ and this subcase follows by an argument similar to Subcase 1.2.
\\
\begin{figure}[h]
\noindent\makebox[\textwidth]{%
\tikzstyle{every node}=[circle, draw, fill=black!100,
                       inner sep=0pt, minimum width=5pt]
\begin{tikzpicture}[thick,scale=0.7]%
\draw(-7,10)node[label={[yshift=0cm, xshift=0cm]90:{\small $[2]_2$}}](1_22){};
\draw(-7,4)node[label={[yshift=0cm, xshift=0cm]270:{\small $[2]_1$}}](1_21){};

\draw(-9,10)node[label={[yshift=0cm, xshift=0cm]90:{\small $[1]_2$}}](1_12){};
\draw(-8,8)node[label={[yshift=0cm, xshift=-0.1cm]180:{\small $\mathtt{c}_2$}}](1_c2){};
\draw(-8,6)node[label={[yshift=0cm, xshift=-0.1cm]180:{\small $\mathtt{c}_1$}}](1_c1){};
\draw(-9,4)node[label={[yshift=0cm, xshift=0cm]270:{\small $[1]_1$}}](1_11){};
\draw(1,10)node[label={[yshift=0cm, xshift=0cm]90:{\small $[2]_2$}}](c_22){};
\draw(1,4)node[label={[yshift=0cm, xshift=0cm]270:{\small $[2]_1$}}](c_21){};

\draw(-1,10)node[label={[yshift=0cm, xshift=0cm]90:{\small $[1]_2$}}](c_12){};
\draw(0,8)node[label={[yshift=-0.4cm, xshift=-0.1cm]180:{\small $\mathtt{c}_2$}}](c_c2){};
\draw(0,6)node[label={[yshift=0.4cm, xshift=-0.1cm]180:{\small $\mathtt{c}_1$}}](c_c1){};
\draw(-1,4)node[label={[yshift=0cm, xshift=0cm]270:{$[1]_1$}}](c_11){};
\draw(9,10)node[label={[yshift=0cm, xshift=0cm]90:{\small $[2]_2$}}](2_22){};
\draw(9,4)node[label={[yshift=0cm, xshift=0cm]270:{\small $[2]_1$}}](2_21){};

\draw(7,10)node[label={[yshift=0cm, xshift=0cm]90:{\small $[1]_2$}}](2_12){};
\draw(8,8)node[label={[yshift=0cm, xshift=0.1cm]0:{\small $\mathtt{c}_2$}}](2_c2){};
\draw(8,6)node[label={[yshift=0cm, xshift=0.1cm]0:{\small $\mathtt{c}_1$}}](2_c1){};
\draw(7,4)node[label={[yshift=0cm, xshift=0cm]270:{\small $[1]_1$}}](2_11){};
\draw [-angle 90, shorten >= 0.1cm, line join=round,decorate, decoration={zigzag, segment length=8,amplitude=1.8,post=lineto,post length=8pt}](1_11)--(1_c1);
\draw [-angle 90, shorten >= 0.1cm, line join=round,decorate, decoration={zigzag, segment length=8,amplitude=1.8,post=lineto,post length=8pt}](1_21)--(1_c1);

\draw [-angle 90, shorten >= 0.1cm, line join=round,decorate, decoration={zigzag, segment length=8,amplitude=1.8,post=lineto,post length=8pt}](1_12)--(1_c2);
\draw [-angle 90, shorten >= 0.1cm, line join=round,decorate, decoration={zigzag, segment length=8,amplitude=1.8,post=lineto,post length=8pt}](1_22)--(1_c2);
\draw [-angle 90, shorten >= 0.1cm, line join=round,decorate, decoration={zigzag, segment length=8,amplitude=1.8,post=lineto,post length=8pt}](c_11)--(c_c1);
\draw [-angle 90, shorten >= 0.1cm, line join=round,decorate, decoration={zigzag, segment length=8,amplitude=1.8,post=lineto,post length=8pt}](c_21)--(c_c1);

\draw [-angle 90, shorten >= 0.1cm, line join=round,decorate, decoration={zigzag, segment length=8,amplitude=1.8,post=lineto,post length=8pt}](c_12)--(c_c2);
\draw [-angle 90, shorten >= 0.1cm, line join=round,decorate, decoration={zigzag, segment length=8,amplitude=1.8,post=lineto,post length=8pt}](c_22)--(c_c2);
\draw [-angle 90, shorten >= 0.1cm, line join=round,decorate, decoration={zigzag, segment length=8,amplitude=1.8,post=lineto,post length=8pt}](2_11)--(2_c1);
\draw [-angle 90, shorten >= 0.1cm, line join=round,decorate, decoration={zigzag, segment length=8,amplitude=1.8,post=lineto,post length=8pt}](2_21)--(2_c1);

\draw [-angle 90, shorten >= 0.1cm, line join=round,decorate, decoration={zigzag, segment length=8,amplitude=1.8,post=lineto,post length=8pt}](2_12)--(2_c2);
\draw [-angle 90, shorten >= 0.1cm, line join=round,decorate, decoration={zigzag, segment length=8,amplitude=1.8,post=lineto,post length=8pt}](2_22)--(2_c2);
\draw [-angle 90, shorten >= 0.1cm, line join=round,decorate, decoration={zigzag, segment length=8,amplitude=1.8,post=lineto,post length=8pt}](-5.5,10)--(-2.5,10);
\draw [-angle 90, shorten >= 0.1cm, line join=round,decorate, decoration={zigzag, segment length=8,amplitude=1.8,post=lineto,post length=8pt}](5.5,10)--(2.5,10);
\draw [-angle 90, shorten >= 0.1cm, line join=round,decorate, decoration={zigzag, segment length=8,amplitude=1.8,post=lineto,post length=8pt}](-5.5,4)--(-2.5,4);
\draw [-angle 90, shorten >= 0.1cm, line join=round,decorate, decoration={zigzag, segment length=8,amplitude=1.8,post=lineto,post length=8pt}](5.5,4)--(2.5,4);
\draw[-angle 90, line width=0.3mm, >=latex, shorten <= 0.2cm, shorten >= 0.15cm](-7.85,8)--(-7.85,6);
\draw[-angle 90, line width=0.3mm, >=latex, shorten <= 0.2cm, shorten >= 0.15cm](-8.15,8)--(-8.15,6);
\draw[-angle 90, line width=0.3mm, >=latex, shorten <= 0.2cm, shorten >= 0.15cm](-0.15,6)--(-0.15,8);
\draw[-angle 90, line width=0.3mm, >=latex, shorten <= 0.2cm, shorten >= 0.15cm](0.15,6)--(0.15,8);
\draw[-angle 90, line width=0.3mm, >=latex, shorten <= 0.2cm, shorten >= 0.15cm](7.85,8)--(7.85,6);
\draw[-angle 90, line width=0.3mm, >=latex, shorten <= 0.2cm, shorten >= 0.15cm](8.15,8)--(8.15,6);
\draw[-angle 90, line width=0.3mm, >=latex, shorten <= 0.2cm, shorten >= 0.15cm](-7.9,5.85)--(-0.1,5.85);
\draw[-angle 90, line width=0.3mm, >=latex, shorten <= 0.2cm, shorten >= 0.15cm](-7.9,6.15)--(-0.1,6.15);
\draw[-angle 90, line width=0.3mm, >=latex, shorten <= 0.2cm, shorten >= 0.15cm](7.9,5.85)--(0.1,5.85);
\draw[-angle 90, line width=0.3mm, >=latex, shorten <= 0.2cm, shorten >= 0.15cm](7.9,6.15)--(0.1,6.15);
\draw[-angle 90, line width=0.3mm, >=latex, shorten <= 0.2cm, shorten >= 0.15cm](-0.1,7.85)--(-7.9,7.85);
\draw[-angle 90, line width=0.3mm, >=latex, shorten <= 0.2cm, shorten >= 0.15cm](-0.1,8.15)--(-7.9,8.15);
\draw[-angle 90, line width=0.3mm, >=latex, shorten <= 0.2cm, shorten >= 0.15cm](0.1,7.85)--(7.9,7.85);
\draw[-angle 90, line width=0.3mm, >=latex, shorten <= 0.2cm, shorten >= 0.15cm](0.1,8.15)--(7.9,8.15);
\node (e1) [E, minimum height=6.3cm, minimum width=2.8cm, label={[yshift=-0.2cm, xshift=0cm]270:{$[1]$}}] at (-8,7){};
\node (e2) [E, minimum height=6.3cm, minimum width=2.8cm, label={[yshift=-0.4cm, xshift=0cm]270:{$\mathtt{c}$}}] at (0,7){};
\node (e3) [E, minimum height=6.3cm, minimum width=2.8cm, label={[yshift=-0.2cm, xshift=0cm]270:{$[2]$}}] at (8,7){};
\end{tikzpicture}
}
{\captionof{figure}{Orientation $D_{(2,3)}$ of $(T_2\times T_3)^{(2)}$, where $d(D_{(2,3)})=5$.}\label{figA8.3.4}}
\caption*{Note: For $u=[1], \mathtt{c}, [2]$, the $u$-copy of $T_3$ is contained in a rectangle and the vertex $\langle u, x\rangle$ is simply labelled as $x$ for clarity. For example, the bottom leftmost vertex is $\langle[1], [1]_1\rangle$. Figures \ref{figA8.3.5} to \ref{figA8.3.12} are labelled similarly.}
\end{figure}
~\\
\noindent\makebox[\textwidth]{%
\tikzstyle{every node}=[circle, draw, fill=black!100,
                       inner sep=0pt, minimum width=5pt]
\begin{tikzpicture}[thick,scale=0.7]
\draw(-7,12)node[label={[yshift=-0.2cm, xshift=0cm]90:{\small $[1,2]_2$}}](1_122){};
\draw(-7,10)node[label={[yshift=0cm, xshift=0cm]315:{\small $[2]_2$}}](1_22){};
\draw(-7,2)node[label={[yshift=0.2cm, xshift=0cm]270:{\small $[1,2]_1$}}](1_121){};
\draw(-7,4)node[label={[yshift=0cm, xshift=0cm]45:{\small $[2]_1$}}](1_21){};

\draw(-9,12)node[label={[yshift=-0.2cm, xshift=0cm]90:{\small $[1,1]_2$}}](1_112){};
\draw(-9,10)node[label={[yshift=0cm, xshift=0cm]215:{\small $[1]_2$}}](1_12){};
\draw(-8,8)node[label={[yshift=0cm, xshift=-0.1cm]180:{\small $\mathtt{c}_2$}}](1_c2){};
\draw(-8,6)node[label={[yshift=0cm, xshift=-0.1cm]180:{\small $\mathtt{c}_1$}}](1_c1){};
\draw(-9,4)node[label={[yshift=0cm, xshift=0cm]135:{\small $[1]_1$}}](1_11){};
\draw(-9,2)node[label={[yshift=0.2cm, xshift=0cm]270:{\small $[1,1]_1$}}](1_111){};
\draw(1,12)node[label={[yshift=-0.2cm, xshift=0cm]90:{\small $[1,2]_2$}}](c_122){};
\draw(1,10)node[label={[yshift=0cm, xshift=0cm]315:{\small $[2]_2$}}](c_22){};
\draw(1,2)node[label={[yshift=0.2cm, xshift=0cm]270:{\small $[1,2]_1$}}](c_121){};
\draw(1,4)node[label={[yshift=0cm, xshift=0cm]45:{\small $[2]_1$}}](c_21){};

\draw(-1,12)node[label={[yshift=-0.2cm, xshift=0cm]90:{\small $[1,1]_2$}}](c_112){};
\draw(-1,10)node[label={[yshift=0cm, xshift=0cm]215:{\small $[1]_2$}}](c_12){};
\draw(0,8)node[label={[yshift=-0.4cm, xshift=-0.1cm]180:{\small $\mathtt{c}_2$}}](c_c2){};
\draw(0,6)node[label={[yshift=0.4cm, xshift=-0.1cm]180:{\small $\mathtt{c}_1$}}](c_c1){};
\draw(-1,4)node[label={[yshift=0cm, xshift=0cm]135:{\small $[1]_1$}}](c_11){};
\draw(-1,2)node[label={[yshift=0.2cm, xshift=0cm]270:{\small $[1,1]_1$}}](c_111){};
\draw(9,12)node[label={[yshift=-0.2cm, xshift=0cm]90:{\small $[1,2]_2$}}](2_122){};
\draw(9,10)node[label={[yshift=0cm, xshift=0cm]315:{\small $[2]_2$}}](2_22){};
\draw(9,2)node[label={[yshift=0.2cm, xshift=0cm]270:{\small $[1,2]_1$}}](2_121){};
\draw(9,4)node[label={[yshift=0cm, xshift=0cm]45:{\small $[2]_1$}}](2_21){};

\draw(7,12)node[label={[yshift=-0.2cm, xshift=0cm]90:{\small $[1,1]_2$}}](2_112){};
\draw(7,10)node[label={[yshift=0cm, xshift=0cm]215:{\small $[1]_2$}}](2_12){};
\draw(8,8)node[label={[yshift=0cm, xshift=0.1cm]0:{\small $\mathtt{c}_2$}}](2_c2){};
\draw(8,6)node[label={[yshift=0cm, xshift=0.1cm]0:{\small $\mathtt{c}_1$}}](2_c1){};
\draw(7,4)node[label={[yshift=0cm, xshift=0cm]135:{\small $[1]_1$}}](2_11){};
\draw(7,2)node[label={[yshift=0.2cm, xshift=0cm]270:{\small $[1,1]_1$}}](2_111){};
\draw [-angle 90, shorten >= 0.1cm, line join=round,decorate, decoration={zigzag, segment length=8,amplitude=1.8,post=lineto,post length=8pt}](1_111)--(1_11);
\draw [-angle 90, shorten >= 0.1cm, line join=round,decorate, decoration={zigzag, segment length=8,amplitude=1.8,post=lineto,post length=8pt}](1_121)--(1_21);
\draw [-angle 90, shorten >= 0.1cm, line join=round,decorate, decoration={zigzag, segment length=8,amplitude=1.8,post=lineto,post length=8pt}](1_11)--(1_c1);
\draw [-angle 90, shorten >= 0.1cm, line join=round,decorate, decoration={zigzag, segment length=8,amplitude=1.8,post=lineto,post length=8pt}](1_21)--(1_c1);

\draw [-angle 90, shorten >= 0.1cm, line join=round,decorate, decoration={zigzag, segment length=8,amplitude=1.8,post=lineto,post length=8pt}](1_112)--(1_12);
\draw [-angle 90, shorten >= 0.1cm, line join=round,decorate, decoration={zigzag, segment length=8,amplitude=1.8,post=lineto,post length=8pt}](1_122)--(1_22);
\draw [-angle 90, shorten >= 0.1cm, line join=round,decorate, decoration={zigzag, segment length=8,amplitude=1.8,post=lineto,post length=8pt}](1_12)--(1_c2);
\draw [-angle 90, shorten >= 0.1cm, line join=round,decorate, decoration={zigzag, segment length=8,amplitude=1.8,post=lineto,post length=8pt}](1_22)--(1_c2);
\draw [-angle 90, shorten >= 0.1cm, line join=round,decorate, decoration={zigzag, segment length=8,amplitude=1.8,post=lineto,post length=8pt}](c_111)--(c_11);
\draw [-angle 90, shorten >= 0.1cm, line join=round,decorate, decoration={zigzag, segment length=8,amplitude=1.8,post=lineto,post length=8pt}](c_121)--(c_21);
\draw [-angle 90, shorten >= 0.1cm, line join=round,decorate, decoration={zigzag, segment length=8,amplitude=1.8,post=lineto,post length=8pt}](c_11)--(c_c1);
\draw [-angle 90, shorten >= 0.1cm, line join=round,decorate, decoration={zigzag, segment length=8,amplitude=1.8,post=lineto,post length=8pt}](c_21)--(c_c1);

\draw [-angle 90, shorten >= 0.1cm, line join=round,decorate, decoration={zigzag, segment length=8,amplitude=1.8,post=lineto,post length=8pt}](c_112)--(c_12);
\draw [-angle 90, shorten >= 0.1cm, line join=round,decorate, decoration={zigzag, segment length=8,amplitude=1.8,post=lineto,post length=8pt}](c_122)--(c_22);
\draw [-angle 90, shorten >= 0.1cm, line join=round,decorate, decoration={zigzag, segment length=8,amplitude=1.8,post=lineto,post length=8pt}](c_12)--(c_c2);
\draw [-angle 90, shorten >= 0.1cm, line join=round,decorate, decoration={zigzag, segment length=8,amplitude=1.8,post=lineto,post length=8pt}](c_22)--(c_c2);
\draw [-angle 90, shorten >= 0.1cm, line join=round,decorate, decoration={zigzag, segment length=8,amplitude=1.8,post=lineto,post length=8pt}](2_111)--(2_11);
\draw [-angle 90, shorten >= 0.1cm, line join=round,decorate, decoration={zigzag, segment length=8,amplitude=1.8,post=lineto,post length=8pt}](2_121)--(2_21);
\draw [-angle 90, shorten >= 0.1cm, line join=round,decorate, decoration={zigzag, segment length=8,amplitude=1.8,post=lineto,post length=8pt}](2_11)--(2_c1);
\draw [-angle 90, shorten >= 0.1cm, line join=round,decorate, decoration={zigzag, segment length=8,amplitude=1.8,post=lineto,post length=8pt}](2_21)--(2_c1);

\draw [-angle 90, shorten >= 0.1cm, line join=round,decorate, decoration={zigzag, segment length=8,amplitude=1.8,post=lineto,post length=8pt}](2_112)--(2_12);
\draw [-angle 90, shorten >= 0.1cm, line join=round,decorate, decoration={zigzag, segment length=8,amplitude=1.8,post=lineto,post length=8pt}](2_122)--(2_22);
\draw [-angle 90, shorten >= 0.1cm, line join=round,decorate, decoration={zigzag, segment length=8,amplitude=1.8,post=lineto,post length=8pt}](2_12)--(2_c2);
\draw [-angle 90, shorten >= 0.1cm, line join=round,decorate, decoration={zigzag, segment length=8,amplitude=1.8,post=lineto,post length=8pt}](2_22)--(2_c2);

\draw[-angle 90, line width=0.3mm, >=latex, shorten <= 0.2cm, shorten >= 0.15cm](-7.85,8)--(-7.85,6);
\draw[-angle 90, line width=0.3mm, >=latex, shorten <= 0.2cm, shorten >= 0.15cm](-8.15,8)--(-8.15,6);
\draw[-angle 90, line width=0.3mm, >=latex, shorten <= 0.2cm, shorten >= 0.15cm](-0.15,6)--(-0.15,8);
\draw[-angle 90, line width=0.3mm, >=latex, shorten <= 0.2cm, shorten >= 0.15cm](0.15,6)--(0.15,8);
\draw[-angle 90, line width=0.3mm, >=latex, shorten <= 0.2cm, shorten >= 0.15cm](7.85,8)--(7.85,6);
\draw[-angle 90, line width=0.3mm, >=latex, shorten <= 0.2cm, shorten >= 0.15cm](8.15,8)--(8.15,6);
\draw [-angle 90, shorten >= 0.1cm, line join=round,decorate, decoration={zigzag, segment length=8,amplitude=1.8,post=lineto,post length=8pt}](-5.5,11)--(-2.5,11);
\draw [-angle 90, shorten >= 0.1cm, line join=round,decorate, decoration={zigzag, segment length=8,amplitude=1.8,post=lineto,post length=8pt}](5.5,11)--(2.5,11);
\draw[-angle 90, line width=0.3mm, >=latex, shorten <= 0.2cm, shorten >= 0.15cm](-7.9,5.85)--(-0.1,5.85);
\draw[-angle 90, line width=0.3mm, >=latex, shorten <= 0.2cm, shorten >= 0.15cm](-7.9,6.15)--(-0.1,6.15);
\draw[-angle 90, line width=0.3mm, >=latex, shorten <= 0.2cm, shorten >= 0.15cm](7.9,5.85)--(0.1,5.85);
\draw[-angle 90, line width=0.3mm, >=latex, shorten <= 0.2cm, shorten >= 0.15cm](7.9,6.15)--(0.1,6.15);
\draw[-angle 90, line width=0.3mm, >=latex, shorten <= 0.2cm, shorten >= 0.15cm](-0.1,7.85)--(-7.9,7.85);
\draw[-angle 90, line width=0.3mm, >=latex, shorten <= 0.2cm, shorten >= 0.15cm](-0.1,8.15)--(-7.9,8.15);
\draw[-angle 90, line width=0.3mm, >=latex, shorten <= 0.2cm, shorten >= 0.15cm](0.1,7.85)--(7.9,7.85);
\draw[-angle 90, line width=0.3mm, >=latex, shorten <= 0.2cm, shorten >= 0.15cm](0.1,8.15)--(7.9,8.15);
\draw [-angle 90, shorten >= 0.1cm, line join=round,decorate, decoration={zigzag, segment length=8,amplitude=1.8,post=lineto,post length=8pt}](-5.5,3)--(-2.5,3);
\draw [-angle 90, shorten >= 0.1cm, line join=round,decorate, decoration={zigzag, segment length=8,amplitude=1.8,post=lineto,post length=8pt}](5.5,3)--(2.5,3);
\node (e1) [E, minimum height=9cm, minimum width=3.2cm, label={[yshift=-0.2cm, xshift=0cm]270:{$[1]$}}] at (-8,7){};
\node (e2) [E, minimum height=9cm, minimum width=3.2cm, label={[yshift=-0.4cm, xshift=0cm]270:{$\mathtt{c}$}}] at (0,7){};
\node (e3) [E, minimum height=9cm, minimum width=3.2cm, label={[yshift=-0.2cm, xshift=0cm]270:{$[2]$}}] at (8,7){};
\end{tikzpicture}
}
{
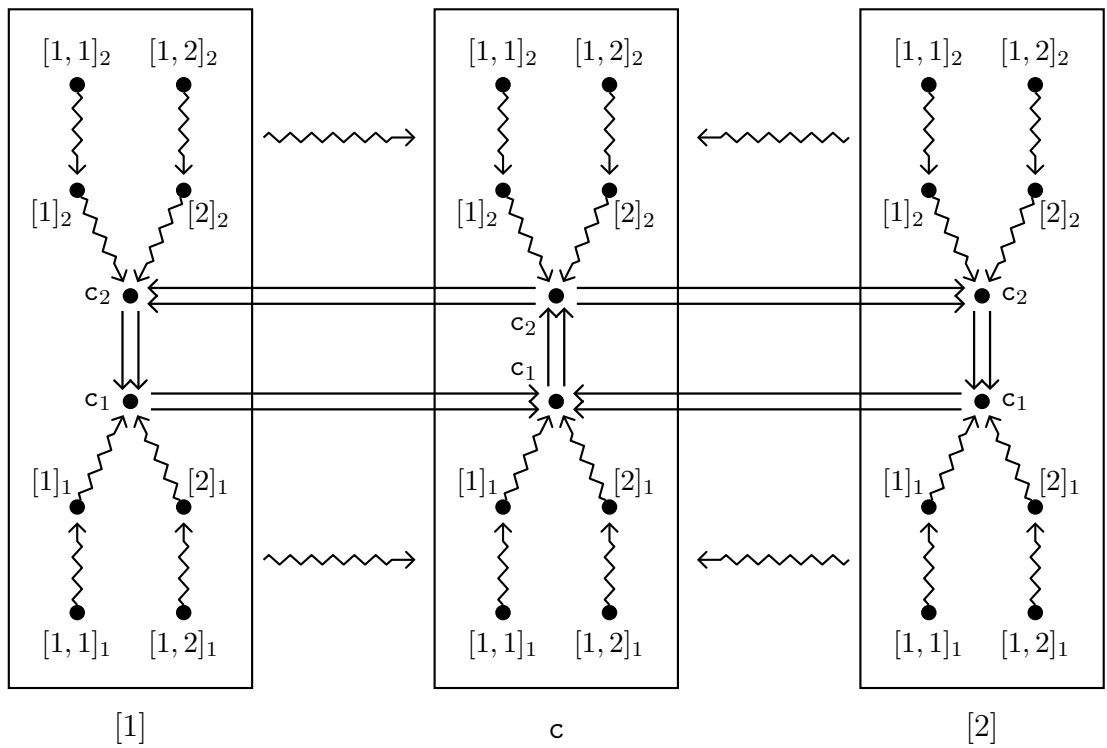
\captionof{figure}{Orientation $D_{(2,5)}$ of $(T_2\times T_5)^{(2)}$, where $d(D_{(2,5)})=7$.}\label{figA8.3.5}}
~\\~\\
\noindent\makebox[\textwidth]{%
\tikzstyle{every node}=[circle, draw, fill=black!100,
                       inner sep=0pt, minimum width=5pt]
\begin{tikzpicture}[thick,scale=0.7]%
\draw(-8,10)node[label={[yshift=0cm, xshift=0cm]90:{\small $[2]_2$}}](11_22){};
\draw(-8,4)node[label={[yshift=0cm, xshift=0cm]270:{\small $[2]_1$}}](11_21){};

\draw(-10,10)node[label={[yshift=0cm, xshift=0cm]90:{\small $[1]_2$}}](11_12){};
\draw(-9,8)node[label={[yshift=0cm, xshift=-0.1cm]180:{\small $\mathtt{c}_2$}}](11_c2){};
\draw(-9,6)node[label={[yshift=0cm, xshift=-0.1cm]180:{\small $\mathtt{c}_1$}}](11_c1){};
\draw(-10,4)node[label={[yshift=0cm, xshift=0cm]270:{\small $[1]_1$}}](11_11){};
\draw(-3.5,10)node[label={[yshift=0cm, xshift=0cm]90:{\small $[2]_2$}}](1_22){};
\draw(-3.5,4)node[label={[yshift=0cm, xshift=0cm]270:{\small $[2]_1$}}](1_21){};

\draw(-5.5,10)node[label={[yshift=0cm, xshift=0cm]90:{\small $[1]_2$}}](1_12){};
\draw(-4.5,8)node[label={[yshift=0cm, xshift=-0.1cm]180:{\small $\mathtt{c}_2$}}](1_c2){};
\draw(-4.5,6)node[label={[yshift=0cm, xshift=-0.1cm]180:{\small $\mathtt{c}_1$}}](1_c1){};
\draw(-5.5,4)node[label={[yshift=0cm, xshift=0cm]270:{\small $[1]_1$}}](1_11){};
\draw(1,10)node[label={[yshift=0cm, xshift=0cm]90:{\small $[2]_2$}}](c_22){};
\draw(1,4)node[label={[yshift=0cm, xshift=0cm]270:{\small $[2]_1$}}](c_21){};

\draw(-1,10)node[label={[yshift=0cm, xshift=0cm]90:{\small $[1]_2$}}](c_12){};
\draw(0,8)node[label={[yshift=-0.4cm, xshift=-0.1cm]180:{\small $\mathtt{c}_2$}}](c_c2){};
\draw(0,6)node[label={[yshift=0.4cm, xshift=-0.1cm]180:{\small $\mathtt{c}_1$}}](c_c1){};
\draw(-1,4)node[label={[yshift=-0cm, xshift=0cm]270:{\small $[1]_1$}}](c_11){};
\draw(5.5,10)node[label={[yshift=0cm, xshift=0cm]90:{\small $[2]_2$}}](2_22){};
\draw(5.5,4)node[label={[yshift=0cm, xshift=0cm]270:{\small $[2]_1$}}](2_21){};

\draw(3.5,10)node[label={[yshift=0cm, xshift=0cm]90:{\small $[1]_2$}}](2_12){};
\draw(4.5,8)node[label={[yshift=0cm, xshift=0.1cm]0:{\small $\mathtt{c}_2$}}](2_c2){};
\draw(4.5,6)node[label={[yshift=0cm, xshift=0.1cm]0:{\small $\mathtt{c}_1$}}](2_c1){};
\draw(3.5,4)node[label={[yshift=0cm, xshift=0cm]270:{\small $[1]_1$}}](2_11){};
\draw(10,10)node[label={[yshift=0cm, xshift=0cm]90:{\small $[2]_2$}}](12_22){};
\draw(10,4)node[label={[yshift=0cm, xshift=0cm]270:{\small $[2]_1$}}](12_21){};

\draw(8,10)node[label={[yshift=0cm, xshift=0cm]90:{\small $[1]_2$}}](12_12){};
\draw(9,8)node[label={[yshift=0cm, xshift=0.1cm]0:{\small $\mathtt{c}_2$}}](12_c2){};
\draw(9,6)node[label={[yshift=0cm, xshift=0.1cm]0:{\small $\mathtt{c}_1$}}](12_c1){};
\draw(8,4)node[label={[yshift=0cm, xshift=0cm]270:{\small $[1]_1$}}](12_11){};
\draw [-angle 90, shorten >= 0.1cm, line join=round,decorate, decoration={zigzag, segment length=8,amplitude=1.8,post=lineto,post length=8pt}](-3,9)--(-1.5,9);
\draw [-angle 90, shorten >= 0.1cm, line join=round,decorate, decoration={zigzag, segment length=8,amplitude=1.8,post=lineto,post length=8pt}](3,9)--(1.5,9);
\draw [-angle 90, shorten >= 0.1cm, line join=round,decorate, decoration={zigzag, segment length=8,amplitude=1.8,post=lineto,post length=8pt}](-7.5,7)--(-6,7);
\draw [-angle 90, shorten >= 0.1cm, line join=round,decorate, decoration={zigzag, segment length=8,amplitude=1.8,post=lineto,post length=8pt}](7.5,7)--(6,7);
\draw [-angle 90, shorten >= 0.1cm, line join=round,decorate, decoration={zigzag, segment length=8,amplitude=1.8,post=lineto,post length=8pt}](-3,5)--(-1.5,5);
\draw [-angle 90, shorten >= 0.1cm, line join=round,decorate, decoration={zigzag, segment length=8,amplitude=1.8,post=lineto,post length=8pt}](3,5)--(1.5,5);
\draw [-angle 90, shorten >= 0.1cm, line join=round,decorate, decoration={zigzag, segment length=8,amplitude=1.8,post=lineto,post length=8pt}](11_12)--(11_c2);
\draw [-angle 90, shorten >= 0.1cm, line join=round,decorate, decoration={zigzag, segment length=8,amplitude=1.8,post=lineto,post length=8pt}](11_22)--(11_c2);
\draw [-angle 90, shorten >= 0.1cm, line join=round,decorate, decoration={zigzag, segment length=8,amplitude=1.8,post=lineto,post length=8pt}](11_c2)--(11_c1);
\draw [-angle 90, shorten >= 0.1cm, line join=round,decorate, decoration={zigzag, segment length=8,amplitude=1.8,post=lineto,post length=8pt}](11_11)--(11_c1);
\draw [-angle 90, shorten >= 0.1cm, line join=round,decorate, decoration={zigzag, segment length=8,amplitude=1.8,post=lineto,post length=8pt}](11_21)--(11_c1);
\draw [-angle 90, shorten >= 0.1cm, line join=round,decorate, decoration={zigzag, segment length=8,amplitude=1.8,post=lineto,post length=8pt}](1_11)--(1_c1);
\draw [-angle 90, shorten >= 0.1cm, line join=round,decorate, decoration={zigzag, segment length=8,amplitude=1.8,post=lineto,post length=8pt}](1_21)--(1_c1);

\draw [-angle 90, shorten >= 0.1cm, line join=round,decorate, decoration={zigzag, segment length=8,amplitude=1.8,post=lineto,post length=8pt}](1_12)--(1_c2);
\draw [-angle 90, shorten >= 0.1cm, line join=round,decorate, decoration={zigzag, segment length=8,amplitude=1.8,post=lineto,post length=8pt}](1_22)--(1_c2);
\draw [-angle 90, shorten >= 0.1cm, line join=round,decorate, decoration={zigzag, segment length=8,amplitude=1.8,post=lineto,post length=8pt}](c_11)--(c_c1);
\draw [-angle 90, shorten >= 0.1cm, line join=round,decorate, decoration={zigzag, segment length=8,amplitude=1.8,post=lineto,post length=8pt}](c_21)--(c_c1);

\draw [-angle 90, shorten >= 0.1cm, line join=round,decorate, decoration={zigzag, segment length=8,amplitude=1.8,post=lineto,post length=8pt}](c_12)--(c_c2);
\draw [-angle 90, shorten >= 0.1cm, line join=round,decorate, decoration={zigzag, segment length=8,amplitude=1.8,post=lineto,post length=8pt}](c_22)--(c_c2);
\draw [-angle 90, shorten >= 0.1cm, line join=round,decorate, decoration={zigzag, segment length=8,amplitude=1.8,post=lineto,post length=8pt}](2_11)--(2_c1);
\draw [-angle 90, shorten >= 0.1cm, line join=round,decorate, decoration={zigzag, segment length=8,amplitude=1.8,post=lineto,post length=8pt}](2_21)--(2_c1);

\draw [-angle 90, shorten >= 0.1cm, line join=round,decorate, decoration={zigzag, segment length=8,amplitude=1.8,post=lineto,post length=8pt}](2_12)--(2_c2);
\draw [-angle 90, shorten >= 0.1cm, line join=round,decorate, decoration={zigzag, segment length=8,amplitude=1.8,post=lineto,post length=8pt}](2_22)--(2_c2);
\draw [-angle 90, shorten >= 0.1cm, line join=round,decorate, decoration={zigzag, segment length=8,amplitude=1.8,post=lineto,post length=8pt}](12_12)--(12_c2);
\draw [-angle 90, shorten >= 0.1cm, line join=round,decorate, decoration={zigzag, segment length=8,amplitude=1.8,post=lineto,post length=8pt}](12_22)--(12_c2);
\draw [-angle 90, shorten >= 0.1cm, line join=round,decorate, decoration={zigzag, segment length=8,amplitude=1.8,post=lineto,post length=8pt}](12_c2)--(12_c1);
\draw [-angle 90, shorten >= 0.1cm, line join=round,decorate, decoration={zigzag, segment length=8,amplitude=1.8,post=lineto,post length=8pt}](12_11)--(12_c1);
\draw [-angle 90, shorten >= 0.1cm, line join=round,decorate, decoration={zigzag, segment length=8,amplitude=1.8,post=lineto,post length=8pt}](12_21)--(12_c1);
\draw[-angle 90, line width=0.3mm, >=latex, shorten <= 0.2cm, shorten >= 0.15cm](-4.35,8)--(-4.35,6);
\draw[-angle 90, line width=0.3mm, >=latex, shorten <= 0.2cm, shorten >= 0.15cm](-4.65,8)--(-4.65,6);
\draw[-angle 90, line width=0.3mm, >=latex, shorten <= 0.2cm, shorten >= 0.15cm](-0.15,6)--(-0.15,8);
\draw[-angle 90, line width=0.3mm, >=latex, shorten <= 0.2cm, shorten >= 0.15cm](0.15,6)--(0.15,8);
\draw[-angle 90, line width=0.3mm, >=latex, shorten <= 0.2cm, shorten >= 0.15cm](4.35,8)--(4.35,6);
\draw[-angle 90, line width=0.3mm, >=latex, shorten <= 0.2cm, shorten >= 0.15cm](4.65,8)--(4.65,6);

\draw[-angle 90, line width=0.3mm, >=latex, shorten <= 0.2cm, shorten >= 0.15cm](-4.5,5.85)--(-0.1,5.85);
\draw[-angle 90, line width=0.3mm, >=latex, shorten <= 0.2cm, shorten >= 0.15cm](-4.5,6.15)--(-0.1,6.15);
\draw[-angle 90, line width=0.3mm, >=latex, shorten <= 0.2cm, shorten >= 0.15cm](4.5,5.85)--(0.1,5.85);
\draw[-angle 90, line width=0.3mm, >=latex, shorten <= 0.2cm, shorten >= 0.15cm](4.5,6.15)--(0.1,6.15);
\draw[-angle 90, line width=0.3mm, >=latex, shorten <= 0.2cm, shorten >= 0.15cm](-0.1,7.85)--(-4.5,7.85);
\draw[-angle 90, line width=0.3mm, >=latex, shorten <= 0.2cm, shorten >= 0.15cm](-0.1,8.15)--(-4.5,8.15);
\draw[-angle 90, line width=0.3mm, >=latex, shorten <= 0.2cm, shorten >= 0.15cm](0.1,7.85)--(4.5,7.85);
\draw[-angle 90, line width=0.3mm, >=latex, shorten <= 0.2cm, shorten >= 0.15cm](0.1,8.15)--(4.5,8.15);
\node (e1) [E, minimum height=6.3cm, minimum width=2.8cm, label={[yshift=0cm, xshift=0cm]270:{$[1,1]$}}] at (-9,7){};
\node (e2) [E, minimum height=6.3cm, minimum width=2.8cm, label={[yshift=-0.2cm, xshift=0cm]270:{$[1]$}}] at (-4.5,7){};
\node (e3) [E, minimum height=6.3cm, minimum width=2.8cm, label={[yshift=-0.4cm, xshift=0cm]270:{$\mathtt{c}$}}] at (0,7){};
\node (e4) [E, minimum height=6.3cm, minimum width=2.8cm, label={[yshift=-0.2cm, xshift=0cm]270:{$[2]$}}] at (4.5,7){};
\node (e5) [E, minimum height=6.3cm, minimum width=2.8cm, label={[yshift=0cm, xshift=0cm]270:{$[1,2]$}}] at (9,7){};
\end{tikzpicture}
}
{
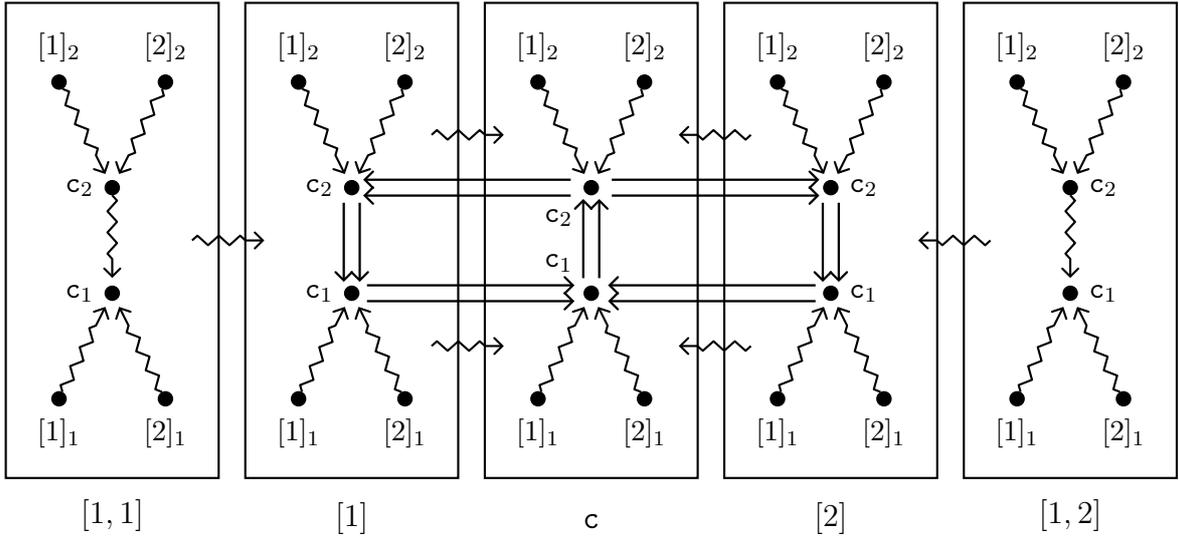
\captionof{figure}{Orientation $D_{(4,3)}$ of $(T_4\times T_3)^{(2)}$, where $d(D_{(4,3)})=7$.}\label{figA8.3.6}}
~\\
\noindent\makebox[\textwidth]{%
\tikzstyle{every node}=[circle, draw, fill=black!100,
                       inner sep=0pt, minimum width=5pt]
\begin{tikzpicture}[thick,scale=0.7]
\draw(-8,11)node[, label={[yshift=-0.2cm, xshift=0cm]90:{\small $[1,2]_2$}}](11_122){};
\draw(-8,9)node[label={[yshift=0cm, xshift=0cm]315:{\small $[2]_2$}}](11_22){};
\draw(-8,1)node[label={[yshift=0.2cm, xshift=0cm]270:{\small $[1,2]_1$}}](11_121){};
\draw(-8,3)node[label={[yshift=0cm, xshift=0cm]45:{\small $[2]_1$}}](11_21){};

\draw(-10,11)node[label={[yshift=-0.2cm, xshift=0cm]90:{\small $[1,1]_2$}}](11_112){};
\draw(-10,9)node[label={[yshift=0cm, xshift=0cm]225:{\small $[1]_2$}}](11_12){};
\draw(-9,7)node[label={[yshift=0cm, xshift=0cm]180:{\small $\mathtt{c}_2$}}](11_c2){};
\draw(-9,5)node[label={[yshift=0cm, xshift=0cm]180:{\small $\mathtt{c}_1$}}](11_c1){};
\draw(-10,3)node[label={[yshift=0cm, xshift=0cm]135:{\small $[1]_1$}}](11_11){};
\draw(-10,1)node[label={[yshift=0.2cm, xshift=0cm]270:{\small $[1,1]_1$}}](11_111){};
\draw(-3.5,11)node[label={[yshift=-0.2cm, xshift=0cm]90:{\small $[1,2]_2$}}](1_122){};
\draw(-3.5,9)node[label={[yshift=0cm, xshift=0cm]315:{\small $[2]_2$}}](1_22){};
\draw(-3.5,1)node[label={[yshift=0.2cm, xshift=0cm]270:{\small $[1,2]_1$}}](1_121){};
\draw(-3.5,3)node[label={[yshift=0cm, xshift=0cm]45:{\small $[2]_1$}}](1_21){};

\draw(-5.5,11)node[label={[yshift=-0.2cm, xshift=0cm]90:{\small $[1,1]_2$}}](1_112){};
\draw(-5.5,9)node[label={[yshift=0cm, xshift=0cm]225:{\small $[1]_2$}}](1_12){};
\draw(-4.5,7)node[label={[yshift=0cm, xshift=0cm]180:{\small $\mathtt{c}_2$}}](1_c2){};
\draw(-4.5,5)node[label={[yshift=0cm, xshift=0cm]180:{\small $\mathtt{c}_1$}}](1_c1){};
\draw(-5.5,3)node[label={[yshift=0cm, xshift=0cm]135:{\small $[1]_1$}}](1_11){};
\draw(-5.5,1)node[label={[yshift=0.2cm, xshift=0cm]270:{\small $[1,1]_1$}}](1_111){};
\draw(1,11)node[label={[yshift=-0.2cm, xshift=0cm]90:{\small $[1,2]_2$}}](c_122){};
\draw(1,9)node[label={[yshift=0cm, xshift=0cm]315:{\small $[2]_2$}}](c_22){};
\draw(1,1)node[label={[yshift=0.2cm, xshift=0cm]270:{\small $[1,2]_1$}}](c_121){};
\draw(1,3)node[label={[yshift=0cm, xshift=0cm]45:{\small $[2]_1$}}](c_21){};

\draw(-1,11)node[label={[yshift=-0.2cm, xshift=0cm]90:{\small $[1,1]_2$}}](c_112){};
\draw(-1,9)node[label={[yshift=0cm, xshift=0cm]225:{\small $[1]_2$}}](c_12){};
\draw(0,7)node[label={[yshift=-0.4cm, xshift=-0.1cm]180:{\small $\mathtt{c}_2$}}](c_c2){};
\draw(0,5)node[label={[yshift=0.4cm, xshift=-0.1cm]180:{\small $\mathtt{c}_1$}}](c_c1){};
\draw(-1,3)node[label={[yshift=0cm, xshift=0cm]135:{\small $[1]_1$}}](c_11){};
\draw(-1,1)node[label={[yshift=0.2cm, xshift=0cm]270:{\small $[1,1]_1$}}](c_111){};
\draw(5.5,11)node[label={[yshift=-0.2cm, xshift=0cm]90:{\small $[1,2]_2$}}](2_122){};
\draw(5.5,9)node[label={[yshift=0cm, xshift=0cm]315:{\small $[2]_2$}}](2_22){};
\draw(5.5,1)node[label={[yshift=0.2cm, xshift=0cm]270:{\small $[1,2]_1$}}](2_121){};
\draw(5.5,3)node[label={[yshift=0cm, xshift=0cm]45:{\small $[2]_1$}}](2_21){};

\draw(3.5,11)node[label={[yshift=-0.2cm, xshift=0cm]90:{\small $[1,1]_2$}}](2_112){};
\draw(3.5,9)node[label={[yshift=0cm, xshift=0cm]225:{\small $[1]_2$}}](2_12){};
\draw(4.5,7)node[label={[yshift=0cm, xshift=0.05cm]0:{\small $\mathtt{c}_2$}}](2_c2){};
\draw(4.5,5)node[label={[yshift=0cm, xshift=0.05cm]0:{\small $\mathtt{c}_1$}}](2_c1){};
\draw(3.5,3)node[label={[yshift=0cm, xshift=0cm]135:{\small $[1]_1$}}](2_11){};
\draw(3.55,1)node[label={[yshift=0.2cm, xshift=0cm]270:{\small $[1,1]_1$}}](2_111){};
\draw(10,11)node[label={[yshift=-0.2cm, xshift=0cm]90:{\small $[1,2]_2$}}](12_122){};
\draw(10,9)node[label={[yshift=0cm, xshift=0cm]315:{\small $[2]_2$}}](12_22){};
\draw(10,1)node[label={[yshift=0.2cm, xshift=0cm]270:{\small $[1,2]_1$}}](12_121){};
\draw(10,3)node[label={[yshift=0cm, xshift=0cm]45:{\small $[2]_1$}}](12_21){};

\draw(8,11)node[label={[yshift=-0.2cm, xshift=0cm]90:{\small $[1,1]_2$}}](12_112){};
\draw(8,9)node[label={[yshift=0cm, xshift=0cm]225:{\small $[1]_2$}}](12_12){};
\draw(9,7)node[label={[yshift=0cm, xshift=0.05cm]0:{\small $\mathtt{c}_2$}}](12_c2){};
\draw(9,5)node[label={[yshift=0cm, xshift=0.05cm]0:{\small $\mathtt{c}_1$}}](12_c1){};
\draw(8,3)node[label={[yshift=0cm, xshift=0cm]135:{\small $[1]_1$}}](12_11){};
\draw(8,1)node[label={[yshift=0.2cm, xshift=0cm]270:{\small $[1,1]_1$}}](12_111){};
\draw [-angle 90, shorten >= 0.1cm, line join=round,decorate, decoration={zigzag, segment length=8,amplitude=1.8,post=lineto,post length=8pt}](11_112)--(11_12);
\draw [-angle 90, shorten >= 0.1cm, line join=round,decorate, decoration={zigzag, segment length=8,amplitude=1.8,post=lineto,post length=8pt}](11_122)--(11_22);
\draw [-angle 90, shorten >= 0.1cm, line join=round,decorate, decoration={zigzag, segment length=8,amplitude=1.8,post=lineto,post length=8pt}](11_21)--(11_c1);
\draw [-angle 90, shorten >= 0.1cm, line join=round,decorate, decoration={zigzag, segment length=8,amplitude=1.8,post=lineto,post length=8pt}](11_22)--(11_c2);
\draw [-angle 90, shorten >= 0.1cm, line join=round,decorate, decoration={zigzag, segment length=8,amplitude=1.8,post=lineto,post length=8pt}](11_c2)--(11_c1);
\draw [-angle 90, shorten >= 0.1cm, line join=round,decorate, decoration={zigzag, segment length=8,amplitude=1.8,post=lineto,post length=8pt}](11_111)--(11_11);
\draw [-angle 90, shorten >= 0.1cm, line join=round,decorate, decoration={zigzag, segment length=8,amplitude=1.8,post=lineto,post length=8pt}](11_121)--(11_21);
\draw [-angle 90, shorten >= 0.1cm, line join=round,decorate, decoration={zigzag, segment length=8,amplitude=1.8,post=lineto,post length=8pt}](11_11)--(11_c1);
\draw [-angle 90, shorten >= 0.1cm, line join=round,decorate, decoration={zigzag, segment length=8,amplitude=1.8,post=lineto,post length=8pt}](11_12)--(11_c2);
\draw [-angle 90, shorten >= 0.1cm, line join=round,decorate, decoration={zigzag, segment length=8,amplitude=1.8,post=lineto,post length=8pt}](1_112)--(1_12);
\draw [-angle 90, shorten >= 0.1cm, line join=round,decorate, decoration={zigzag, segment length=8,amplitude=1.8,post=lineto,post length=8pt}](1_122)--(1_22);
\draw [-angle 90, shorten >= 0.1cm, line join=round,decorate, decoration={zigzag, segment length=8,amplitude=1.8,post=lineto,post length=8pt}](1_21)--(1_c1);
\draw [-angle 90, shorten >= 0.1cm, line join=round,decorate, decoration={zigzag, segment length=8,amplitude=1.8,post=lineto,post length=8pt}](1_22)--(1_c2);
\draw [-angle 90, shorten >= 0.1cm, line join=round,decorate, decoration={zigzag, segment length=8,amplitude=1.8,post=lineto,post length=8pt}](1_111)--(1_11);
\draw [-angle 90, shorten >= 0.1cm, line join=round,decorate, decoration={zigzag, segment length=8,amplitude=1.8,post=lineto,post length=8pt}](1_121)--(1_21);
\draw [-angle 90, shorten >= 0.1cm, line join=round,decorate, decoration={zigzag, segment length=8,amplitude=1.8,post=lineto,post length=8pt}](1_11)--(1_c1);
\draw [-angle 90, shorten >= 0.1cm, line join=round,decorate, decoration={zigzag, segment length=8,amplitude=1.8,post=lineto,post length=8pt}](1_12)--(1_c2);
\draw [-angle 90, shorten >= 0.1cm, line join=round,decorate, decoration={zigzag, segment length=8,amplitude=1.8,post=lineto,post length=8pt}](c_112)--(c_12);
\draw [-angle 90, shorten >= 0.1cm, line join=round,decorate, decoration={zigzag, segment length=8,amplitude=1.8,post=lineto,post length=8pt}](c_122)--(c_22);
\draw [-angle 90, shorten >= 0.1cm, line join=round,decorate, decoration={zigzag, segment length=8,amplitude=1.8,post=lineto,post length=8pt}](c_21)--(c_c1);
\draw [-angle 90, shorten >= 0.1cm, line join=round,decorate, decoration={zigzag, segment length=8,amplitude=1.8,post=lineto,post length=8pt}](c_22)--(c_c2);
\draw [-angle 90, shorten >= 0.1cm, line join=round,decorate, decoration={zigzag, segment length=8,amplitude=1.8,post=lineto,post length=8pt}](c_111)--(c_11);
\draw [-angle 90, shorten >= 0.1cm, line join=round,decorate, decoration={zigzag, segment length=8,amplitude=1.8,post=lineto,post length=8pt}](c_121)--(c_21);
\draw [-angle 90, shorten >= 0.1cm, line join=round,decorate, decoration={zigzag, segment length=8,amplitude=1.8,post=lineto,post length=8pt}](c_11)--(c_c1);
\draw [-angle 90, shorten >= 0.1cm, line join=round,decorate, decoration={zigzag, segment length=8,amplitude=1.8,post=lineto,post length=8pt}](c_12)--(c_c2);
\draw [-angle 90, shorten >= 0.1cm, line join=round,decorate, decoration={zigzag, segment length=8,amplitude=1.8,post=lineto,post length=8pt}](2_112)--(2_12);
\draw [-angle 90, shorten >= 0.1cm, line join=round,decorate, decoration={zigzag, segment length=8,amplitude=1.8,post=lineto,post length=8pt}](2_122)--(2_22);
\draw [-angle 90, shorten >= 0.1cm, line join=round,decorate, decoration={zigzag, segment length=8,amplitude=1.8,post=lineto,post length=8pt}](2_21)--(2_c1);
\draw [-angle 90, shorten >= 0.1cm, line join=round,decorate, decoration={zigzag, segment length=8,amplitude=1.8,post=lineto,post length=8pt}](2_22)--(2_c2);
\draw [-angle 90, shorten >= 0.1cm, line join=round,decorate, decoration={zigzag, segment length=8,amplitude=1.8,post=lineto,post length=8pt}](2_111)--(2_11);
\draw [-angle 90, shorten >= 0.1cm, line join=round,decorate, decoration={zigzag, segment length=8,amplitude=1.8,post=lineto,post length=8pt}](2_121)--(2_21);
\draw [-angle 90, shorten >= 0.1cm, line join=round,decorate, decoration={zigzag, segment length=8,amplitude=1.8,post=lineto,post length=8pt}](2_11)--(2_c1);
\draw [-angle 90, shorten >= 0.1cm, line join=round,decorate, decoration={zigzag, segment length=8,amplitude=1.8,post=lineto,post length=8pt}](2_12)--(2_c2);
\draw [-angle 90, shorten >= 0.1cm, line join=round,decorate, decoration={zigzag, segment length=8,amplitude=1.8,post=lineto,post length=8pt}](12_112)--(12_12);
\draw [-angle 90, shorten >= 0.1cm, line join=round,decorate, decoration={zigzag, segment length=8,amplitude=1.8,post=lineto,post length=8pt}](12_122)--(12_22);
\draw [-angle 90, shorten >= 0.1cm, line join=round,decorate, decoration={zigzag, segment length=8,amplitude=1.8,post=lineto,post length=8pt}](12_21)--(12_c1);
\draw [-angle 90, shorten >= 0.1cm, line join=round,decorate, decoration={zigzag, segment length=8,amplitude=1.8,post=lineto,post length=8pt}](12_22)--(12_c2);
\draw [-angle 90, shorten >= 0.1cm, line join=round,decorate, decoration={zigzag, segment length=8,amplitude=1.8,post=lineto,post length=8pt}](12_c2)--(12_c1);
\draw [-angle 90, shorten >= 0.1cm, line join=round,decorate, decoration={zigzag, segment length=8,amplitude=1.8,post=lineto,post length=8pt}](12_111)--(12_11);
\draw [-angle 90, shorten >= 0.1cm, line join=round,decorate, decoration={zigzag, segment length=8,amplitude=1.8,post=lineto,post length=8pt}](12_121)--(12_21);
\draw [-angle 90, shorten >= 0.1cm, line join=round,decorate, decoration={zigzag, segment length=8,amplitude=1.8,post=lineto,post length=8pt}](12_11)--(12_c1);
\draw [-angle 90, shorten >= 0.1cm, line join=round,decorate, decoration={zigzag, segment length=8,amplitude=1.8,post=lineto,post length=8pt}](12_12)--(12_c2);
\draw [-angle 90, shorten >= 0.1cm, line join=round,decorate, decoration={zigzag, segment length=8,amplitude=1.8,post=lineto,post length=8pt}](-3,10)--(-1.5,10);
\draw [-angle 90, shorten >= 0.1cm, line join=round,decorate, decoration={zigzag, segment length=8,amplitude=1.8,post=lineto,post length=8pt}](3,10)--(1.5,10);
\draw [-angle 90, shorten >= 0.1cm, line join=round,decorate, decoration={zigzag, segment length=8,amplitude=1.8,post=lineto,post length=8pt}](-7.5,6)--(-6,6);
\draw [-angle 90, shorten >= 0.1cm, line join=round,decorate, decoration={zigzag, segment length=8,amplitude=1.8,post=lineto,post length=8pt}](7.5,6)--(6,6);
\draw [-angle 90, shorten >= 0.1cm, line join=round,decorate, decoration={zigzag, segment length=8,amplitude=1.8,post=lineto,post length=8pt}](-3,2)--(-1.5,2);
\draw [-angle 90, shorten >= 0.1cm, line join=round,decorate, decoration={zigzag, segment length=8,amplitude=1.8,post=lineto,post length=8pt}](3,2)--(1.5,2);
\draw[-angle 90, line width=0.3mm, >=latex, shorten <= 0.2cm, shorten >= 0.15cm](-0.1,6.85)--(-4.4,6.85);
\draw[-angle 90, line width=0.3mm, >=latex, shorten <= 0.2cm, shorten >= 0.15cm](-0.1,7.15)--(-4.4,7.15);
\draw[-angle 90, line width=0.3mm, >=latex, shorten <= 0.2cm, shorten >= 0.15cm](0.1,6.85)--(4.4,6.85);
\draw[-angle 90, line width=0.3mm, >=latex, shorten <= 0.2cm, shorten >= 0.15cm](0.1,7.15)--(4.4,7.15);
\draw[-angle 90, line width=0.3mm, >=latex, shorten <= 0.2cm, shorten >= 0.15cm](-4.65,6.9)--(-4.65,5);
\draw[-angle 90, line width=0.3mm, >=latex, shorten <= 0.2cm, shorten >= 0.15cm](-4.35,6.9)--(-4.35,5);
\draw[-angle 90, line width=0.3mm, >=latex, shorten <= 0.2cm, shorten >= 0.15cm](-0.15,5.1)--(-0.15,6.9);
\draw[-angle 90, line width=0.3mm, >=latex, shorten <= 0.2cm, shorten >= 0.15cm](0.15,5.1)--(0.15,6.9);
\draw[-angle 90, line width=0.3mm, >=latex, shorten <= 0.2cm, shorten >= 0.15cm](4.35,6.9)--(4.355,5);
\draw[-angle 90, line width=0.3mm, >=latex, shorten <= 0.2cm, shorten >= 0.15cm](4.65,6.9)--(4.65,5);
\draw[-angle 90, line width=0.3mm, >=latex, shorten <= 0.2cm, shorten >= 0.15cm](-4.4,5.15)--(-0.1,5.15);
\draw[-angle 90, line width=0.3mm, >=latex, shorten <= 0.2cm, shorten >= 0.15cm](-4.4,4.85)--(-0.1,4.85);
\draw[-angle 90, line width=0.3mm, >=latex, shorten <= 0.2cm, shorten >= 0.15cm](4.4,5.15)--(0.1,5.15);
\draw[-angle 90, line width=0.3mm, >=latex, shorten <= 0.2cm, shorten >= 0.15cm](4.4,4.85)--(0.1,4.85);
\node (e1) [E, minimum height=9cm, minimum width=2.8cm, label={[yshift=0cm, xshift=0cm]270:{$[1,1]$}}] at (-9,6){};
\node (e2) [E, minimum height=9cm, minimum width=2.8cm,  label={[yshift=-0.2cm, xshift=0cm]270:{$[1]$}}] at (-4.5,6){};
\node (e3) [E, minimum height=9cm, minimum width=2.8cm,  label={[yshift=-0.4cm, xshift=0cm]270:{$\mathtt{c}$}}] at (0,6){};
\node (e4) [E, minimum height=9cm, minimum width=2.8cm,  label={[yshift=-0.2cm, xshift=0cm]270:{$[2]$}}] at (4.5,6){};
\node (e5) [E, minimum height=9cm, minimum width=2.8cm,  label={[yshift=0cm, xshift=0cm]270:{$[1,2]$}}] at (9,6){};
\end{tikzpicture}
}
{
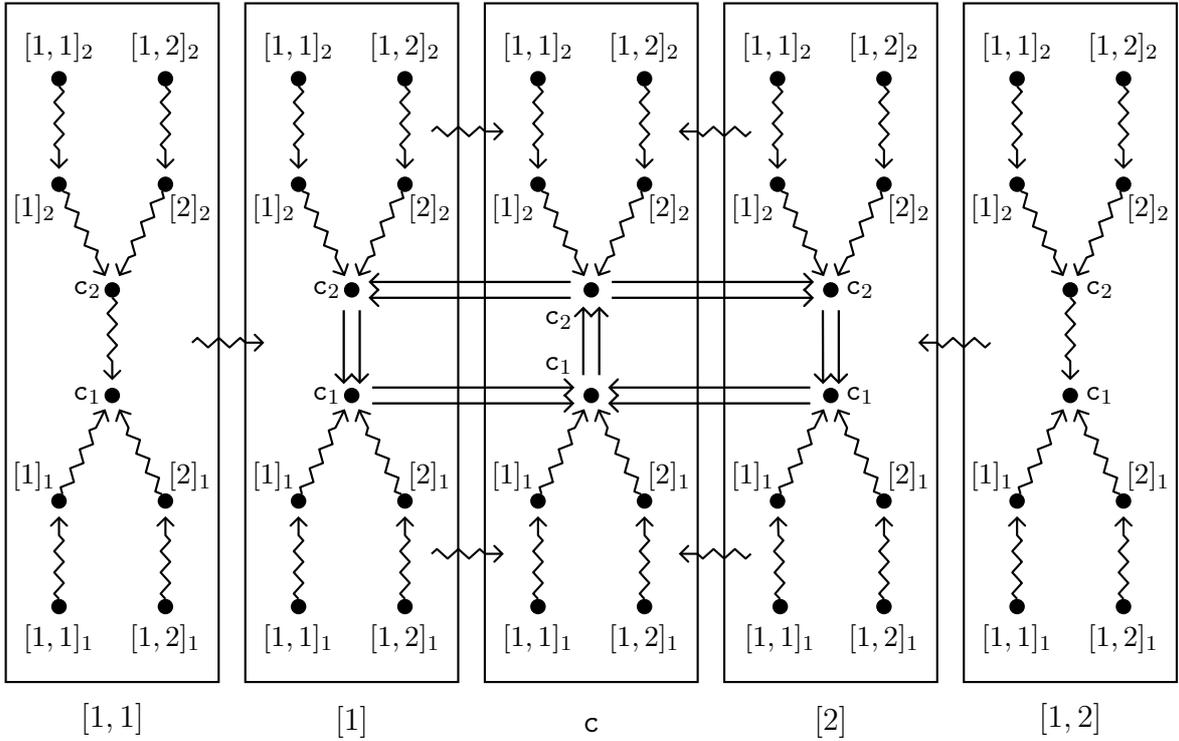
\captionof{figure}{Orientation $D_{(4,5)}$ of $(T_4\times T_5)^{(2)}$, where $d(D_{(4,5)})=9$.}\label{figA8.3.7}}

\newpage
\noindent Case 2. $\lambda$ and $\mu$ are both even, i.e. $\lambda=2,4$ and $\mu=4$.
\indent\par For each $[i]\in N_{T_\lambda}(\mathtt{c})$, and each $\alpha=1,2,\ldots, deg_{T_\lambda}([i])-1$ and each $[j]\in N_{T_\mu}(\mathtt{c})$, and each $\beta=1,2,\ldots, deg_{T_\mu}([j])-1$,
\begin{align}
&\langle [i],\mathtt{c}\rangle \rightrightarrows \langle \mathtt{c},\mathtt{c}\rangle \text{ and }\langle\mathtt{c},[j]\rangle \rightrightarrows \langle [i],[j]\rangle;
\label{eqA8.3.6}
\end{align}
excluding the edges defined above,
\begin{align}
&\langle [\alpha,i],y\rangle\rightsquigarrow \langle [i],y\rangle \rightsquigarrow \langle \mathtt{c},y\rangle \label{eqA8.3.7}
\end{align}
for all $y\in V(T_\mu)$, and 
\begin{align}
\langle x,[\beta,j]\rangle \rightsquigarrow \langle x,[j]\rangle \rightsquigarrow \langle x,\mathtt{c}\rangle \label{eqA8.3.8}
\end{align}
for all $x\in V(T_\lambda)$. (See Figures \ref{figA8.3.8} and \ref{figA8.3.9}.)
\indent\par We use the same strategy as before to prove $d(D_{(\lambda,\mu)})=\lambda+\mu$. Let $u,v\in V(G)$ and $P: u=w_0 w_1 \ldots  w_k=v$ be a shortest $u-v$ path in $G$. If $d_G(u,v)\le d(G)-2$ and $P$ satisfies (\ref{eqA8.3.4}), then $d_{D_{(\lambda,\mu)}}((p,u),(q,v))\le d_G(u,v)+2\le d(G)$ for $p,q=1,2$. Particularly, this holds for $u=\langle [\alpha_1,k],y_1\rangle$, $v=\langle[ \alpha_2,k],y_2\rangle$ with $\alpha_1\neq \alpha_2$. So, by symmetry of (\ref{eqA8.3.6})-(\ref{eqA8.3.8}), we may assume $\mathtt{c}$ has two eccentric vertices in $T_\lambda$, i.e. $T_\lambda=P_3$ if $\lambda=2$, and $T_\lambda=P_5$ if $\lambda=4$. Furthermore, by symmetry of (\ref{eqA8.3.8}), we may assume $\mathtt{c}$ has two eccentric vertices $[1,1]$ and $[1,2]$ in $T_\mu$.
\indent\par For the pairs of $u,v$ that do not satisfy (\ref{eqA8.3.4}), we claim that there exists a path $P$ with length at most $d(G)$ and satisfies (\ref{eqA8.3.5}). Hence, $d_{D_{(\lambda,\mu)}}((p,u),(q,v))\le d(G)$.
\\
\\Subcase 2.1. $\lambda=2$ and $\mu=4$. (See Figure \ref{figA8.3.8}.)
\indent\par We list these paths $P$ in each of the respective subcases, omitting symmetric scenarios. For $i=1,2$, and $j=2,3$,
\begin{align*}
P^1=&\langle [1],[1,1]\rangle \langle [1],[1]\rangle \langle [1],\mathtt{c}\rangle \langle \mathtt{c},\mathtt{c}\rangle \langle\mathtt{c},[j]\rangle \langle\mathtt{c},[1,j]\rangle \langle [2],[1,j]\rangle.\\
P^2=&\langle [1],[1,1]\rangle \langle [1],[1]\rangle \langle [1],\mathtt{c}\rangle \langle \mathtt{c},\mathtt{c}\rangle \langle\mathtt{c},[j]\rangle \langle [2],[j]\rangle.\\
P^3=&\langle\mathtt{c},[1,1]\rangle \langle\mathtt{c},[1]\rangle \langle\mathtt{c},\mathtt{c}\rangle \langle\mathtt{c},[j]\rangle \langle [i],[j]\rangle \langle[i],[1,j]\rangle.
\end{align*}
Subcase 2.2. $\lambda=4$ and $\mu=4$. (See Figure \ref{figA8.3.9}.)
\indent\par Note that $D_{(2,4)}$ is a subdigraph of $D_{(4,4)}$ and this subcase follows by an argument similar to Subcase 1.2.
\\
\\
\noindent\makebox[\textwidth]{%
\tikzstyle{every node}=[circle, draw, fill=black!100,
                       inner sep=0pt, minimum width=5pt]
\begin{tikzpicture}[thick,scale=0.7]%
\draw(-8,10)node[label={[yshift=-0.1cm, xshift=0cm]90:{\small $[1,2]$}}](1_12){};
\draw(-8,8)node[label={[yshift=0cm, xshift=0cm]180:{\small $[2]$}}](1_2){};
\draw(-8,6)node[label={[yshift=0cm, xshift=-0.15cm]180:{\small $\mathtt{c}$}}](1_c){};
\draw(-8,4)node[label={[yshift=0cm, xshift=0cm]180:{\small $[1]$}}](1_1){};
\draw(-8,2)node[label={[yshift=0.1cm, xshift=0cm]270:{\small $[1,1]$}}](1_11){};
\draw(0,10)node[label={[yshift=-0.1cm, xshift=0cm]90:{\small $[1,2]$}}](c_12){};
\draw(0,8)node[label={[yshift=-0.4cm, xshift=0cm]180:{\small $[2]$}}](c_2){};
\draw(0,6)node[label={[yshift=-0.4cm, xshift=-0.15cm]180:{\small $\mathtt{c}$}}](c_c){};
\draw(0,4)node[label={[yshift=-0.4cm, xshift=0cm]180:{\small $[1]$}}](c_1){};
\draw(0,2)node[label={[yshift=0.1cm, xshift=0cm]270:{\small $[1,1]$}}](c_11){};
\draw(8,10)node[label={[yshift=-0.1cm, xshift=0cm]90:{\small $[1,2]$}}](2_12){};
\draw(8,8)node[label={[yshift=0cm, xshift=0cm]0:{\small $[2]$}}](2_2){};
\draw(8,6)node[label={[yshift=0cm, xshift=0.15cm]0:{\small $\mathtt{c}$}}](2_c){};
\draw(8,4)node[label={[yshift=0cm, xshift=0cm]0:{\small $[1]$}}](2_1){};
\draw(8,2)node[label={[yshift=0.1cm, xshift=0cm]270:{\small $[1,1]$}}](2_11){};
\draw [-angle 90, shorten >= 0.1cm, line join=round,decorate, decoration={zigzag, segment length=8,amplitude=1.8,post=lineto,post length=8pt}](1_11)--(1_1);
\draw [-angle 90, shorten >= 0.1cm, line join=round,decorate, decoration={zigzag, segment length=8,amplitude=1.8,post=lineto,post length=8pt}](1_1)--(1_c);

\draw [-angle 90, shorten >= 0.1cm, line join=round,decorate, decoration={zigzag, segment length=8,amplitude=1.8,post=lineto,post length=8pt}](1_12)--(1_2);
\draw [-angle 90, shorten >= 0.1cm, line join=round,decorate, decoration={zigzag, segment length=8,amplitude=1.8,post=lineto,post length=8pt}](1_2)--(1_c);
\draw [-angle 90, shorten >= 0.1cm, line join=round,decorate, decoration={zigzag, segment length=8,amplitude=1.8,post=lineto,post length=8pt}](c_11)--(c_1);
\draw [-angle 90, shorten >= 0.1cm, line join=round,decorate, decoration={zigzag, segment length=8,amplitude=1.8,post=lineto,post length=8pt}](c_1)--(c_c);

\draw [-angle 90, shorten >= 0.1cm, line join=round,decorate, decoration={zigzag, segment length=8,amplitude=1.8,post=lineto,post length=8pt}](c_12)--(c_2);
\draw [-angle 90, shorten >= 0.1cm, line join=round,decorate, decoration={zigzag, segment length=8,amplitude=1.8,post=lineto,post length=8pt}](c_2)--(c_c);
\draw [-angle 90, shorten >= 0.1cm, line join=round,decorate, decoration={zigzag, segment length=8,amplitude=1.8,post=lineto,post length=8pt}](2_11)--(2_1);
\draw [-angle 90, shorten >= 0.1cm, line join=round,decorate, decoration={zigzag, segment length=8,amplitude=1.8,post=lineto,post length=8pt}](2_1)--(2_c);

\draw [-angle 90, shorten >= 0.1cm, line join=round,decorate, decoration={zigzag, segment length=8,amplitude=1.8,post=lineto,post length=8pt}](2_12)--(2_2);

\draw [-angle 90, shorten >= 0.1cm, line join=round,decorate, decoration={zigzag, segment length=8,amplitude=1.8,post=lineto,post length=8pt}](2_2)--(2_c);
\draw[-angle 90, line width=0.3mm, >=latex, shorten <= 0.2cm, shorten >= 0.15cm](-8,6.15)--(0,6.15);
\draw[-angle 90, line width=0.3mm, >=latex, shorten <= 0.2cm, shorten >= 0.15cm](-8,5.85)--(0,5.85);

\draw [-angle 90, shorten >= 0.1cm, line join=round,decorate, decoration={zigzag, segment length=8,amplitude=1.8,post=lineto,post length=8pt}](-6,10)--(-2,10);
\draw [-angle 90, shorten >= 0.1cm, line join=round,decorate, decoration={zigzag, segment length=8,amplitude=1.8,post=lineto,post length=8pt}](-6,2)--(-2,2);
\draw[-angle 90, line width=0.3mm, >=latex, shorten <= 0.2cm, shorten >= 0.15cm](8,6.15)--(0,6.15);
\draw[-angle 90, line width=0.3mm, >=latex, shorten <= 0.2cm, shorten >= 0.15cm](8,5.85)--(0,5.85);

\draw [-angle 90, shorten >= 0.1cm, line join=round,decorate, decoration={zigzag, segment length=8,amplitude=1.8,post=lineto,post length=8pt}](6,10)--(2,10);
\draw [-angle 90, shorten >= 0.1cm, line join=round,decorate, decoration={zigzag, segment length=8,amplitude=1.8,post=lineto,post length=8pt}](6,2)--(2,2);
\draw[-angle 90, line width=0.3mm, >=latex, shorten <= 0.2cm, shorten >= 0.15cm](0,8.15)--(-8,8.15);
\draw[-angle 90, line width=0.3mm, >=latex, shorten <= 0.2cm, shorten >= 0.15cm](0,7.85)--(-8,7.85);
\draw[-angle 90, line width=0.3mm, >=latex, shorten <= 0.2cm, shorten >= 0.15cm](0,4.15)--(-8,4.15);
\draw[-angle 90, line width=0.3mm, >=latex, shorten <= 0.2cm, shorten >= 0.15cm](0,3.85)--(-8,3.85);
\draw[-angle 90, line width=0.3mm, >=latex, shorten <= 0.2cm, shorten >= 0.15cm](0,8.15)--(8,8.15);
\draw[-angle 90, line width=0.3mm, >=latex, shorten <= 0.2cm, shorten >= 0.15cm](0,7.85)--(8,7.85);
\draw[-angle 90, line width=0.3mm, >=latex, shorten <= 0.2cm, shorten >= 0.15cm](0,4.15)--(8,4.15);
\draw[-angle 90, line width=0.3mm, >=latex, shorten <= 0.2cm, shorten >= 0.15cm](0,3.85)--(8,3.85);

\node (e1) [E, minimum height=7.6cm, minimum width=2.6cm, label={[yshift=-0.2cm, xshift=0cm]270:{$[1]$}}] at (-8,6){};
\node (e2) [E, minimum height=7.7cm, minimum width=2.6cm, label={[yshift=-0.4cm, xshift=0cm]270:{$\mathtt{c}$}}] at (0,6){};
\node (e3) [E, minimum height=7.8cm, minimum width=2.6cm, label={[yshift=-0.2cm, xshift=0cm]270:{$[2]$}}] at (8,6){};
\end{tikzpicture}
}
{
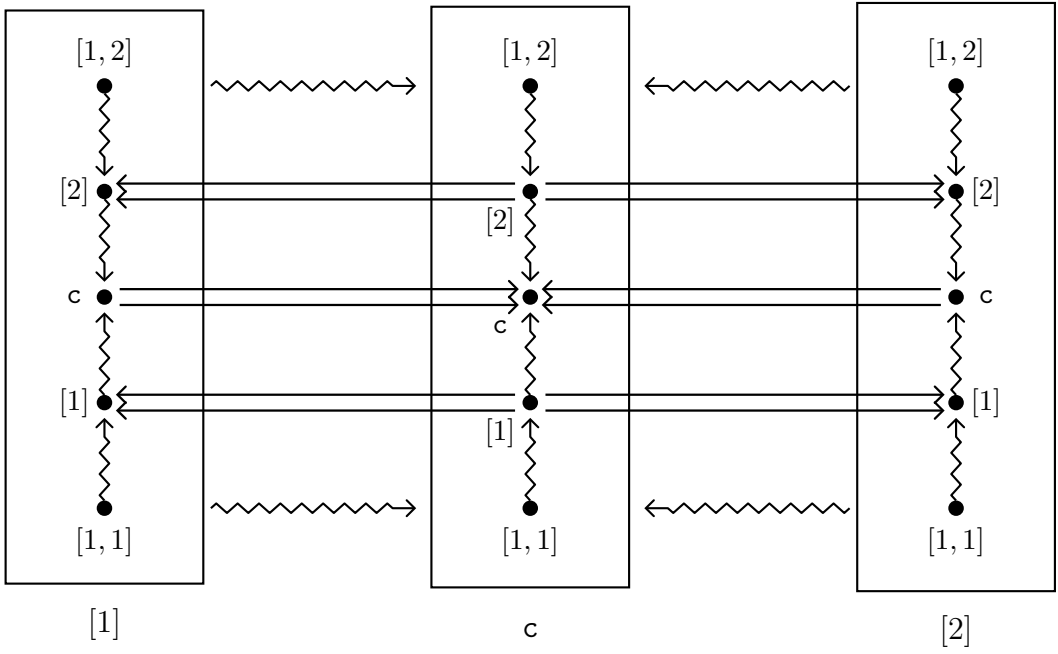
\captionof{figure}{Orientation $D_{(2,4)}$ of $(T_2\times T_4)^{(2)}$, where $d(D_{(2,4)})=6$.}\label{figA8.3.8}}
~\\
\noindent\makebox[\textwidth]{%
\tikzstyle{every node}=[circle, draw, fill=black!100,
                       inner sep=0pt, minimum width=5pt]
\begin{tikzpicture}[thick,scale=0.7]%
\draw(-9,10)node[label={[yshift=-0.1cm, xshift=0cm]90:{\small $[1,2]$}}](11_12){};
\draw(-9,8)node[label={[yshift=0cm, xshift=0cm]180:{\small $[2]$}}](11_2){};
\draw(-9,6)node[label={[yshift=0cm, xshift=-0.15cm]180:{\small $\mathtt{c}$}}](11_c){};
\draw(-9,4)node[label={[yshift=0cm, xshift=0cm]180:{\small $[1]$}}](11_1){};
\draw(-9,2)node[label={[yshift=0.1cm, xshift=0cm]270:{\small $[1,1]$}}](11_11){};
\draw(-4.5,10)node[label={[yshift=-0.1cm, xshift=0cm]90:{\small $[1,2]$}}](1_12){};
\draw(-4.5,8)node[label={[yshift=0cm, xshift=0cm]180:{\small $[2]$}}](1_2){};
\draw(-4.5,6)node[label={[yshift=0cm, xshift=-0.15cm]180:{\small $\mathtt{c}$}}](1_c){};
\draw(-4.5,4)node[label={[yshift=0cm, xshift=0cm]180:{\small $[1]$}}](1_1){};
\draw(-4.5,2)node[label={[yshift=0.1cm, xshift=0cm]270:{\small $[1,1]$}}](1_11){};
\draw(0,10)node[label={[yshift=-0.1cm, xshift=0cm]90:{\small $[1,2]$}}](c_12){};
\draw(0,8)node[label={[yshift=-0.4cm, xshift=0cm]180:{\small $[2]$}}](c_2){};
\draw(-0,6)node[label={[yshift=-0.4cm, xshift=-0.15cm]180:{\small $\mathtt{c}$}}](c_c){};
\draw(0,4)node[label={[yshift=-0.4cm, xshift=0cm]180:{\small $[1]$}}](c_1){};
\draw(0,2)node[label={[yshift=0.1cm, xshift=0cm]270:{\small $[1,1]$}}](c_11){};
\draw(4.5,10)node[label={[yshift=-0.1cm, xshift=0cm]90:{\small $[1,2]$}}](2_12){};
\draw(4.5,8)node[label={[yshift=0cm, xshift=0cm]0:{\small $[2]$}}](2_2){};
\draw(4.5,6)node[label={[yshift=0cm, xshift=0.15cm]0:{\small $\mathtt{c}$}}](2_c){};
\draw(4.5,4)node[label={[yshift=0cm, xshift=0cm]0:{\small $[1]$}}](2_1){};
\draw(4.5,2)node[label={[yshift=0.1cm, xshift=0cm]270:{\small $[1,1]$}}](2_11){};
\draw(9,10)node[label={[yshift=-0.1cm, xshift=0cm]90:{\small $[1,2]$}}](12_12){};
\draw(9,8)node[label={[yshift=0cm, xshift=0cm]0:{\small $[2]$}}](12_2){};
\draw(9,6)node[label={[yshift=0cm, xshift=0.15cm]0:{\small $\mathtt{c}$}}](12_c){};
\draw(9,4)node[label={[yshift=0cm, xshift=0cm]0:{\small $[1]$}}](12_1){};
\draw(9,2)node[label={[yshift=0.1cm, xshift=0cm]270:{\small $[1,1]$}}](12_11){};
\draw [-angle 90, shorten >= 0.1cm, line join=round,decorate, decoration={zigzag, segment length=8,amplitude=1.8,post=lineto,post length=8pt}](11_11)--(11_1);
\draw [-angle 90, shorten >= 0.1cm, line join=round,decorate, decoration={zigzag, segment length=8,amplitude=1.8,post=lineto,post length=8pt}](11_1)--(11_c);

\draw [-angle 90, shorten >= 0.1cm, line join=round,decorate, decoration={zigzag, segment length=8,amplitude=1.8,post=lineto,post length=8pt}](11_12)--(11_2);
\draw [-angle 90, shorten >= 0.1cm, line join=round,decorate, decoration={zigzag, segment length=8,amplitude=1.8,post=lineto,post length=8pt}](11_2)--(11_c);
\draw [-angle 90, shorten >= 0.1cm, line join=round,decorate, decoration={zigzag, segment length=8,amplitude=1.8,post=lineto,post length=8pt}](1_11)--(1_1);
\draw [-angle 90, shorten >= 0.1cm, line join=round,decorate, decoration={zigzag, segment length=8,amplitude=1.8,post=lineto,post length=8pt}](1_1)--(1_c);

\draw [-angle 90, shorten >= 0.1cm, line join=round,decorate, decoration={zigzag, segment length=8,amplitude=1.8,post=lineto,post length=8pt}](1_12)--(1_2);
\draw [-angle 90, shorten >= 0.1cm, line join=round,decorate, decoration={zigzag, segment length=8,amplitude=1.8,post=lineto,post length=8pt}](1_2)--(1_c);
\draw [-angle 90, shorten >= 0.1cm, line join=round,decorate, decoration={zigzag, segment length=8,amplitude=1.8,post=lineto,post length=8pt}](c_11)--(c_1);
\draw [-angle 90, shorten >= 0.1cm, line join=round,decorate, decoration={zigzag, segment length=8,amplitude=1.8,post=lineto,post length=8pt}](c_1)--(c_c);

\draw [-angle 90, shorten >= 0.1cm, line join=round,decorate, decoration={zigzag, segment length=8,amplitude=1.8,post=lineto,post length=8pt}](c_12)--(c_2);
\draw [-angle 90, shorten >= 0.1cm, line join=round,decorate, decoration={zigzag, segment length=8,amplitude=1.8,post=lineto,post length=8pt}](c_2)--(c_c);
\draw [-angle 90, shorten >= 0.1cm, line join=round,decorate, decoration={zigzag, segment length=8,amplitude=1.8,post=lineto,post length=8pt}](2_11)--(2_1);
\draw [-angle 90, shorten >= 0.1cm, line join=round,decorate, decoration={zigzag, segment length=8,amplitude=1.8,post=lineto,post length=8pt}](2_1)--(2_c);

\draw [-angle 90, shorten >= 0.1cm, line join=round,decorate, decoration={zigzag, segment length=8,amplitude=1.8,post=lineto,post length=8pt}](2_12)--(2_2);
\draw [-angle 90, shorten >= 0.1cm, line join=round,decorate, decoration={zigzag, segment length=8,amplitude=1.8,post=lineto,post length=8pt}](2_2)--(2_c);
\draw [-angle 90, shorten >= 0.1cm, line join=round,decorate, decoration={zigzag, segment length=8,amplitude=1.8,post=lineto,post length=8pt}](12_11)--(12_1);
\draw [-angle 90, shorten >= 0.1cm, line join=round,decorate, decoration={zigzag, segment length=8,amplitude=1.8,post=lineto,post length=8pt}](12_1)--(12_c);

\draw [-angle 90, shorten >= 0.1cm, line join=round,decorate, decoration={zigzag, segment length=8,amplitude=1.8,post=lineto,post length=8pt}](12_12)--(12_2);
\draw [-angle 90, shorten >= 0.1cm, line join=round,decorate, decoration={zigzag, segment length=8,amplitude=1.8,post=lineto,post length=8pt}](12_2)--(12_c);
\draw [-angle 90, shorten >= 0.1cm, line join=round,decorate, decoration={zigzag, segment length=8,amplitude=1.8,post=lineto,post length=8pt}](-7.5,6)--(-6,6);
\draw [-angle 90, shorten >= 0.1cm, line join=round,decorate, decoration={zigzag, segment length=8,amplitude=1.8,post=lineto,post length=8pt}](7.5,6)--(6,6);
\draw[-angle 90, line width=0.3mm, >=latex, shorten <= 0.2cm, shorten >= 0.15cm](-4.5,6.15)--(0,6.15);
\draw[-angle 90, line width=0.3mm, >=latex, shorten <= 0.2cm, shorten >= 0.15cm](-4.5,5.85)--(0,5.85);

\draw [-angle 90, shorten >= 0.1cm, line join=round,decorate, decoration={zigzag, segment length=8,amplitude=1.8,post=lineto,post length=8pt}](-3,10)--(-1.5,10);
\draw [-angle 90, shorten >= 0.1cm, line join=round,decorate, decoration={zigzag, segment length=8,amplitude=1.8,post=lineto,post length=8pt}](-3,2)--(-1.5,2);
\draw[-angle 90, line width=0.3mm, >=latex, shorten <= 0.2cm, shorten >= 0.15cm](4.5,6.15)--(0,6.15);
\draw[-angle 90, line width=0.3mm, >=latex, shorten <= 0.2cm, shorten >= 0.15cm](4.5,5.85)--(0,5.85);

\draw [-angle 90, shorten >= 0.1cm, line join=round,decorate, decoration={zigzag, segment length=8,amplitude=1.8,post=lineto,post length=8pt}](3,10)--(1.5,10);
\draw [-angle 90, shorten >= 0.1cm, line join=round,decorate, decoration={zigzag, segment length=8,amplitude=1.8,post=lineto,post length=8pt}](3,2)--(1.5,2);
\draw[-angle 90, line width=0.3mm, >=latex, shorten <= 0.2cm, shorten >= 0.15cm](0,8.15)--(-4.5,8.15);
\draw[-angle 90, line width=0.3mm, >=latex, shorten <= 0.2cm, shorten >= 0.15cm](0,7.85)--(-4.5,7.85);

\draw[-angle 90, line width=0.3mm, >=latex, shorten <= 0.2cm, shorten >= 0.15cm](0,4.15)--(-4.5,4.15);
\draw[-angle 90, line width=0.3mm, >=latex, shorten <= 0.2cm, shorten >= 0.15cm](0,3.85)--(-4.5,3.85);
\draw[-angle 90, line width=0.3mm, >=latex, shorten <= 0.2cm, shorten >= 0.15cm](0.2,8.15)--(4.5,8.15);
\draw[-angle 90, line width=0.3mm, >=latex, shorten <= 0.2cm, shorten >= 0.15cm](0.2, 7.85)--(4.5,7.85);

\draw[-angle 90, line width=0.3mm, >=latex, shorten <= 0.2cm, shorten >= 0.15cm](0,4.15)--(4.5,4.15);
\draw[-angle 90, line width=0.3mm, >=latex, shorten <= 0.2cm, shorten >= 0.15cm](0,3.85)--(4.5,3.85);
\node (e1) [E, minimum height=7.6cm, minimum width=2.6cm, label={[yshift=0cm, xshift=0cm]270:{$[1,1]$}}] at (-9,6){};
\node (e2) [E, minimum height=7.6cm, minimum width=2.6cm,  label={[yshift=-0.2cm, xshift=0cm]270:{$[1]$}}] at (-4.5,6){};
\node (e3) [E, minimum height=7.6cm, minimum width=2.6cm,  label={[yshift=-0.4cm, xshift=0cm]270:{$\mathtt{c}$}}] at (0,6){};
\node (e4) [E, minimum height=7.6cm, minimum width=2.6cm,  label={[yshift=-0.2cm, xshift=0cm]270:{$[2]$}}] at (4.5,6){};
\node (e5) [E, minimum height=7.6cm, minimum width=2.6cm,  label={[yshift=0cm, xshift=0cm]270:{$[1,2]$}}] at (9,6){};
\end{tikzpicture}
}
{
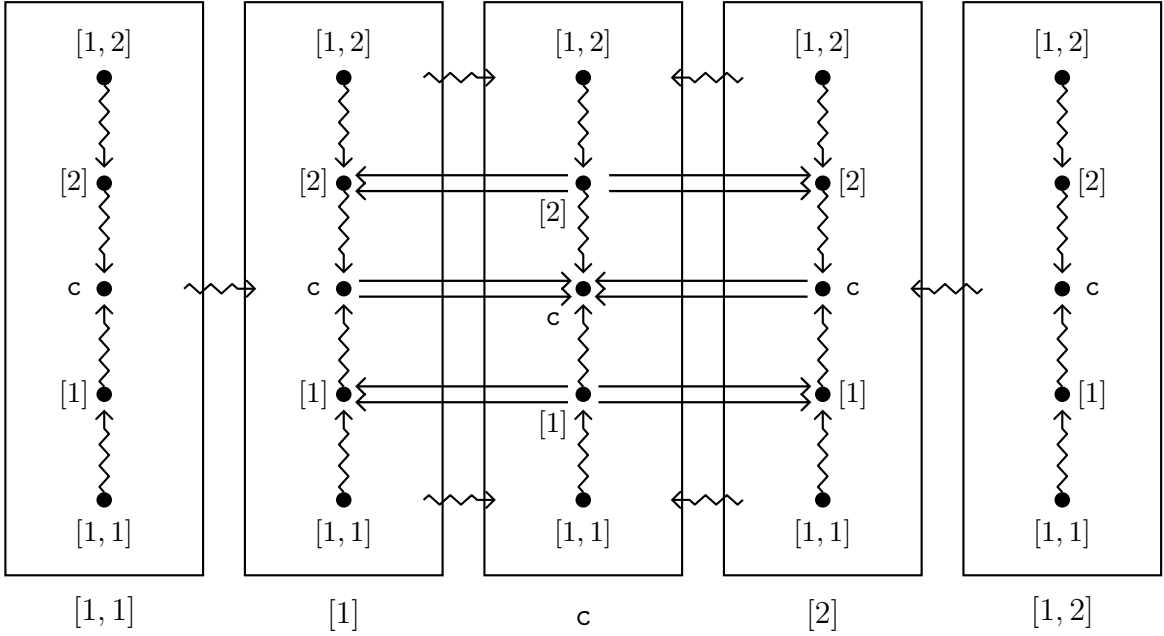
\captionof{figure}{Orientation $D_{(4,4)}$ of $(T_4\times T_4)^{(2)}$, where $d(D_{(4,4)})=8$.}\label{figA8.3.9}}

~\\Case 3. $\lambda$ and $\mu$ are both odd, i.e. $\lambda,\mu=3,5$.
\indent\par For each $[i]_1\in N_{T_\mu}(\mathtt{c}_1)-\{\mathtt{c}_2\}$, each $\alpha=1,2,\ldots,deg_{T_\mu}([i]_1)-1$, each $[j]_2\in N_{T_\mu}(\mathtt{c}_2)-\{\mathtt{c}_1\}$, each $\beta=1,2,\ldots,deg_{T_\mu}([j]_2)-1$,
\begin{equation}
\left. \begin{array}{@{}ll@{}}
&\langle\mathtt{c}_1,[i]_1\rangle 
\rightrightarrows \langle \mathtt{c}_1,\mathtt{c}_1\rangle
\rightrightarrows \langle\mathtt{c}_2,\mathtt{c}_1\rangle
\rightrightarrows \langle\mathtt{c}_2,[i]_1\rangle
\rightrightarrows \langle\mathtt{c}_1,[i]_1\rangle,\\
&\langle\mathtt{c}_1,[j]_2\rangle 
\rightrightarrows \langle \mathtt{c}_1,\mathtt{c}_2\rangle
\rightrightarrows \langle\mathtt{c}_2,\mathtt{c}_2\rangle
\rightrightarrows \langle\mathtt{c}_2,[j]_2\rangle
\rightrightarrows \langle\mathtt{c}_1,[j]_2\rangle;
  \end{array}\right\}
\label{eqA8.3.9}
\end{equation}
excluding the edges defined above,
\begin{align}
\langle x,[\alpha,i]_1\rangle \rightsquigarrow \langle x,[i]_1\rangle \rightsquigarrow \langle x,\mathtt{c}_1\rangle,
\langle x,[\beta,j]_2\rangle \rightsquigarrow \langle x,[j]_2\rangle \rightsquigarrow \langle x,\mathtt{c}_2\rangle\rightsquigarrow \langle x,\mathtt{c}_1\rangle \label{eqA8.3.10}
\end{align}
for all $x\in V(T_\lambda)$, and
\begin{align}
&\langle [\gamma,k]_1,y\rangle\rightsquigarrow \langle [k]_1,y\rangle \rightsquigarrow \langle \mathtt{c}_1,y\rangle,
\langle [\theta,l]_2,y\rangle\rightsquigarrow \langle [l]_2,y\rangle \rightsquigarrow \langle \mathtt{c}_2,y\rangle \rightsquigarrow \langle \mathtt{c}_1,y\rangle, \label{eqA8.3.11}
\end{align}
for each $[k]_1\in N_{T_\lambda}(\mathtt{c}_1)-\{\mathtt{c}_2\}$, each $\gamma=1,2,\ldots, \deg_{T_\lambda}([k]_1)-1$, each $[l]_2\in N_{T_\lambda}(\mathtt{c}_2)-\{\mathtt{c}_1\}$, each $\theta=1,2,\ldots, \deg_{T_\lambda}([l]_2)-1$ and each $y\in V(T_\mu)$. (See Figures \ref{figA8.3.10}-\ref{figA8.3.12}.)
\indent\par We use the same strategy as before to prove $d(D_{(\lambda,\mu)})=\lambda+\mu$. Let $u,v\in V(G)$ and $P: u=w_0 w_1 \ldots w_k=v$ be a shortest $u-v$ path in $G$. If $d_G(u,v)\le d(G)-2$ and $P$ satisfies (\ref{eqA8.3.4}), then $d_{D_{(\lambda,\mu)}}((p,u),(q,v))\le d_G(u,v)+2\le d(G)$ for $p,q=1,2$.
\indent\par By symmetry of (\ref{eqA8.3.9})-(\ref{eqA8.3.11}), we may assume for $i=1,2$, $N_{T_\lambda}(\mathtt{c}_i)=\{[1]_i,[2]_i\}$, and $N_{T_\lambda}([j]_i)=\{[1,j]_i\}$ for $j=1,2$. Furthermore, by symmetry of (\ref{eqA8.3.10}), we may assume the same holds for $T_\mu$. For the pairs of $u,v$ that do not satisfy (\ref{eqA8.3.4}), we claim that there exists a path $P$ with length at most $d(G)$ and satisfies (\ref{eqA8.3.5}). Hence, $d_{D_{(\lambda,\mu)}}((p,u),(q,v))\le d(G)$.
\\
\\Subcase 3.1. $\lambda=3$ and $\mu=3$. (See Figure \ref{figA8.3.10}.)
\indent\par We list these paths $P$ while omitting symmetric scenarios. For $j=1,2$,
\begin{align*}
P^1=&\langle [1]_1,[1]_1\rangle \langle [1]_1,\mathtt{c}_1\rangle \langle \mathtt{c}_1,\mathtt{c}_1\rangle \langle \mathtt{c}_2,\mathtt{c}_1\rangle \langle \mathtt{c}_2,[j]_1\rangle \langle [1]_2,[j]_1\rangle \langle [1]_2,\mathtt{c}_1\rangle.\\
P^2=&\langle [1]_1,[1]_1\rangle \langle [1]_1,\mathtt{c}_1\rangle \langle \mathtt{c}_1,\mathtt{c}_1\rangle \langle \mathtt{c}_2,\mathtt{c}_1\rangle \langle \mathtt{c}_2,\mathtt{c}_2\rangle \langle \mathtt{c}_2,[j]_2\rangle \langle [1]_1,[j]_2\rangle.\\
P^3=&\langle [1]_1,[1]_1\rangle \langle [1]_1,\mathtt{c}_1\rangle \langle \mathtt{c}_1,\mathtt{c}_1\rangle \langle \mathtt{c}_2,\mathtt{c}_1\rangle \langle \mathtt{c}_2,\mathtt{c}_2\rangle \langle [1]_2,\mathtt{c}_2\rangle.\\
P^4=&\langle [1]_1,[1]_1\rangle \langle \mathtt{c},[1]_1\rangle \langle \mathtt{c}_1,\mathtt{c}_1\rangle \langle\mathtt{c}_1,\mathtt{c}_2\rangle \langle[2]_1,\mathtt{c}_2\rangle \langle [2]_1,[j]_2\rangle.\\
P^5=&\langle [1]_2,[1]_2\rangle \langle [1]_2,\mathtt{c}_2\rangle \langle [1]_2,\mathtt{c}_1\rangle \langle [1]_2,[j]_1\rangle\langle \mathtt{c}_2,[j]_1\rangle \langle \mathtt{c}_1,[j]_1\rangle \langle [1]_1,[j]_1\rangle,\\
P^6=&\langle [1]_2,[1]_2\rangle \langle \mathtt{c}_2,[1]_2\rangle \langle \mathtt{c}_1,[1]_2\rangle \langle [1]_1,[j]_2\rangle \langle [1]_1,\mathtt{c}_2\rangle \langle [1]_1,\mathtt{c}_1\rangle.\\
P^7=&\langle [1]_2,[1]_2\rangle \langle \mathtt{c}_2,[1]_2\rangle \langle \mathtt{c}_1,[1]_2\rangle \langle \mathtt{c}_1,\mathtt{c}_2\rangle \langle \mathtt{c}_1,\mathtt{c}_1\rangle.\\
P^8=&\langle [1]_2,[1]_2\rangle \langle [1]_2,\mathtt{c}_2\rangle \langle \mathtt{c}_2,\mathtt{c}_2\rangle \langle \mathtt{c}_2,\mathtt{c}_1\rangle \langle \mathtt{c}_2,[j]_1\rangle \langle [2]_2,[j]_1\rangle.\\
P^9=&\langle[1]_2,\mathtt{c}_2\rangle \langle [1]_2,[1]_2\rangle \langle \mathtt{c}_2,[1]_2\rangle \langle \mathtt{c}_1,[1]_2\rangle \langle \mathtt{c}_1,\mathtt{c}_2\rangle \langle \mathtt{c}_1,\mathtt{c}_1\rangle.\\
P^{10}=&\langle[1]_2,\mathtt{c}_2\rangle \langle [1]_2,[1]_2\rangle \langle \mathtt{c}_2,[j]_2\rangle \langle \mathtt{c}_1,[j]_2\rangle \langle [1]_1,[j]_2\rangle \langle [1]_1,\mathtt{c}_2\rangle \langle [1]_1,\mathtt{c}_1\rangle.\\
P^{11}=&\langle[1]_2,\mathtt{c}_2\rangle \langle [1]_2,\mathtt{c}_1\rangle \langle [1]_2,[j]_1\rangle \langle\mathtt{c}_2,[j]_1\rangle \langle\mathtt{c}_1,[j]_1\rangle \langle [1]_1,[j]_1\rangle.\\
P^{12}=&\langle \mathtt{c}_1,[1]_1\rangle \langle \mathtt{c}_1,\mathtt{c}_1\rangle \langle \mathtt{c}_2,\mathtt{c}_1\rangle \langle \mathtt{c}_2,\mathtt{c}_2\rangle \langle \mathtt{c}_2,[j]_2\rangle \langle \mathtt{c}_1,[j]_2\rangle.\\
P^{13}=&\langle \mathtt{c}_1,[1]_1\rangle \langle \mathtt{c}_1,\mathtt{c}_1\rangle \langle \mathtt{c}_2,\mathtt{c}_1\rangle \langle \mathtt{c}_2,[j]_1\rangle \langle \mathtt{c}_1,[j]_1\rangle.\\
P^{14}=&\langle \mathtt{c}_1,[1]_1\rangle \langle \mathtt{c}_1,\mathtt{c}_1\rangle \langle \mathtt{c}_2,\mathtt{c}_1\rangle \langle \mathtt{c}_2,[j]_1\rangle \langle [1]_2,[j]_1\rangle.\\
P^{15}=&\langle \mathtt{c}_1,[1]_1\rangle \langle \mathtt{c}_1,\mathtt{c}_1\rangle \langle \mathtt{c}_2,\mathtt{c}_1\rangle \langle [1]_2,\mathtt{c}_1\rangle \langle [1]_2,\mathtt{c}_2\rangle \langle [1]_2,[j]_2\rangle.\\
P^{16}=&\langle \mathtt{c}_2,\mathtt{c}_1\rangle \langle \mathtt{c}_2,[j]_1\rangle \langle \mathtt{c}_1,[j]_1\rangle \langle [1]_1,[j]_1\rangle \langle [1]_1,\mathtt{c}_1\rangle \langle [1]_1,\mathtt{c}_2\rangle \langle [1]_1,[j]_2\rangle.\\
P^{17}=&\langle \mathtt{c}_2,\mathtt{c}_1\rangle \langle \mathtt{c}_2,[j]_1\rangle \langle \mathtt{c}_1,[j]_1\rangle \langle \mathtt{c}_1,\mathtt{c}_1\rangle \langle \mathtt{c}_1,\mathtt{c}_2\rangle.\\
P^{18}=&\langle \mathtt{c}_2,\mathtt{c}_1\rangle \langle \mathtt{c}_2,\mathtt{c}_2\rangle \langle \mathtt{c}_2,[j]_2\rangle \langle \mathtt{c}_1,[j]_2\rangle.\\
P^{19}=&\langle \mathtt{c}_2,[1]_1\rangle \langle \mathtt{c}_1,[1]_1\rangle \langle [1]_1,[1]_1\rangle \langle [1]_1,\mathtt{c}_1\rangle \langle [1]_1,[2]_1\rangle.\\
P^{20}=&\langle \mathtt{c}_2,[1]_1\rangle \langle \mathtt{c}_1,[1]_1\rangle \langle [1]_1,[1]_1\rangle \langle [1]_1,\mathtt{c}_1\rangle \langle [1]_1,\mathtt{c}_2\rangle \langle [1]_1,[j]_2\rangle \langle \mathtt{c}_1,[j]_2\rangle.\\
P^{21}=&\langle \mathtt{c}_2,[1]_1\rangle \langle \mathtt{c}_1,[1]_1\rangle \langle \mathtt{c}_1,\mathtt{c}_1\rangle \langle \mathtt{c}_2,\mathtt{c}_1\rangle \langle \mathtt{c}_2,[j]_1\rangle \langle \mathtt{c}_1,[j]_1\rangle.\\
P^{22}=&\langle \mathtt{c}_2,[1]_1\rangle \langle \mathtt{c}_1,[1]_1\rangle \langle \mathtt{c}_1,\mathtt{c}_1\rangle \langle \mathtt{c}_2,\mathtt{c}_1\rangle \langle \mathtt{c}_2,\mathtt{c}_2\rangle \langle \mathtt{c}_2,[j]_2\rangle. \\
P^{23}=&\langle \mathtt{c}_2,[1]_1\rangle \langle \mathtt{c}_1,[1]_1\rangle \langle \mathtt{c}_1,\mathtt{c}_1\rangle \langle \mathtt{c}_1,\mathtt{c}_2\rangle.
\end{align*}
Subcase 3.2. $\lambda=3$ and $\mu=5$. (See Figure \ref{figA8.3.11}.)
\indent\par Note that $D_{(3,3)}$ is a subdigraph of $D_{(3,5)}$ and this subcase follows by an argument similar to Subcase 1.2.
\\
\\Subcase 3.3. $\lambda=5$ and $\mu=5$. (See Figure \ref{figA8.3.12}.)
\indent\par Note that $D_{(3,5)}$ is a subdigraph of $D_{(5,5)}$ and this subcase follows by an argument similar to Subcase 1.2. \qed
~\\~\\
\noindent\makebox[\textwidth]{%
\tikzstyle{every node}=[circle, draw, fill=black!100,
                       inner sep=0pt, minimum width=5pt]
\begin{tikzpicture}[thick,scale=0.65]%
\draw(-6.5,10)node[label={[yshift=0cm, xshift=0cm]90:{\small $[2]_2$}}](11_22){};
\draw(-6.5,4)node[label={[yshift=0cm, xshift=0cm]270:{\small $[2]_1$}}](11_21){};

\draw(-8.5,10)node[label={[yshift=0cm, xshift=0cm]90:{\small $[1]_2$}}](11_12){};
\draw(-7.5,8)node[label={[yshift=0cm, xshift=-0.1cm]180:{\small $\mathtt{c}_2$}}](11_c2){};
\draw(-7.5,6)node[label={[yshift=0cm, xshift=-0.1cm]180:{\small $\mathtt{c}_1$}}](11_c1){};
\draw(-8.5,4)node[label={[yshift=0cm, xshift=0cm]270:{\small $[1]_1$}}](11_11){};
\draw(-1.5,9.5)node[label={[yshift=-0.1cm, xshift=0cm]90:{\small $[2]_2$}}](c1_22){};
\draw(-1.5,3.5)node[label={[yshift=0cm, xshift=0cm]270:{\small $[2]_1$}}](c1_21){};

\draw(-3.5,10.5)node[label={[yshift=0cm, xshift=0cm]90:{\small $[1]_2$}}](c1_12){};
\draw(-2.5,8)node[label={[yshift=0cm, xshift=-0.1cm]180:{\small $\mathtt{c}_2$}}](c1_c2){};
\draw(-2.5,6)node[label={[yshift=0cm, xshift=-0.1cm]180:{\small $\mathtt{c}_1$}}](c1_c1){};
\draw(-3.5,4.5)node[label={[yshift=0.1cm, xshift=0cm]270:{\small $[1]_1$}}](c1_11){};
\draw(3.5,9.5)node[label={[yshift=-0.1cm, xshift=0cm]90:{\small $[2]_2$}}](c2_22){};
\draw(3.5,3.5)node[label={[yshift=0cm, xshift=0cm]270:{\small $[2]_1$}}](c2_21){};

\draw(1.5,10.5)node[label={[yshift=0cm, xshift=0cm]90:{\small $[1]_2$}}](c2_12){};
\draw(2.5,8)node[label={[yshift=0cm, xshift=0.1cm]0:{\small $\mathtt{c}_2$}}](c2_c2){};
\draw(2.5,6)node[label={[yshift=0cm, xshift=0.1cm]0:{\small $\mathtt{c}_1$}}](c2_c1){};
\draw(1.5,4.5)node[label={[yshift=0.1cm, xshift=0cm]270:{\small $[1]_1$}}](c2_11){};
\draw(8.5,10)node[label={[yshift=0cm, xshift=0cm]90:{\small $[2]_2$}}](12_22){};
\draw(8.5,4)node[label={[yshift=0cm, xshift=0cm]270:{\small $[2]_1$}}](12_21){};

\draw(6.5,10)node[label={[yshift=0cm, xshift=0cm]90:{\small $[1]_2$}}](12_12){};
\draw(7.5,8)node[label={[yshift=0cm, xshift=0.1cm]0:{\small $\mathtt{c}_2$}}](12_c2){};
\draw(7.5,6)node[label={[yshift=0cm, xshift=0.1cm]0:{\small $\mathtt{c}_1$}}](12_c1){};
\draw(6.5,4)node[label={[yshift=0cm, xshift=0cm]270:{\small $[1]_1$}}](12_11){};
\draw [-angle 90, shorten >= 0.1cm, line join=round,decorate, decoration={zigzag, segment length=8,amplitude=1.8,post=lineto,post length=8pt}](11_11)--(11_c1);
\draw [-angle 90, shorten >= 0.1cm, line join=round,decorate, decoration={zigzag, segment length=8,amplitude=1.8,post=lineto,post length=8pt}](11_21)--(11_c1);

\draw [-angle 90, shorten >= 0.1cm, line join=round,decorate, decoration={zigzag, segment length=8,amplitude=1.8,post=lineto,post length=8pt}](11_12)--(11_c2);
\draw [-angle 90, shorten >= 0.1cm, line join=round,decorate, decoration={zigzag, segment length=8,amplitude=1.8,post=lineto,post length=8pt}](11_22)--(11_c2);

\draw [-angle 90, shorten >= 0.1cm, line join=round,decorate, decoration={zigzag, segment length=8,amplitude=1.8,post=lineto,post length=8pt}](11_c2)--(11_c1);
\draw [-angle 90, shorten >= 0.1cm, line join=round,decorate, decoration={zigzag, segment length=8,amplitude=1.8,post=lineto,post length=8pt}](c1_c2)--(c1_c1);
\draw [-angle 90, shorten >= 0.1cm, line join=round,decorate, decoration={zigzag, segment length=8,amplitude=1.8,post=lineto,post length=8pt}](c2_c2)--(c2_c1);
\draw [-angle 90, shorten >= 0.1cm, line join=round,decorate, decoration={zigzag, segment length=8,amplitude=1.8,post=lineto,post length=8pt}](12_11)--(12_c1);
\draw [-angle 90, shorten >= 0.1cm, line join=round,decorate, decoration={zigzag, segment length=8,amplitude=1.8,post=lineto,post length=8pt}](12_21)--(12_c1);

\draw [-angle 90, shorten >= 0.1cm, line join=round,decorate, decoration={zigzag, segment length=8,amplitude=1.8,post=lineto,post length=8pt}](12_12)--(12_c2);
\draw [-angle 90, shorten >= 0.1cm, line join=round,decorate, decoration={zigzag, segment length=8,amplitude=1.8,post=lineto,post length=8pt}](12_22)--(12_c2);

\draw [-angle 90, shorten >= 0.1cm, line join=round,decorate, decoration={zigzag, segment length=8,amplitude=1.8,post=lineto,post length=8pt}](12_c2)--(12_c1);
\draw [-angle 90, shorten >= 0.1cm, line join=round,decorate, decoration={zigzag, segment length=8,amplitude=1.8,post=lineto,post length=8pt}](-6,7)--(-4,7);
\draw [-angle 90, shorten >= 0.1cm, line join=round,decorate, decoration={zigzag, segment length=8,amplitude=1.8,post=lineto,post length=8pt}](6,7)--(4,7);
\draw[-angle 90, line width=0.3mm, >=latex, shorten <= 0.2cm, shorten >= 0.15cm](1.5,4.65)--(-3.4,4.65);
\draw[-angle 90, line width=0.3mm, >=latex, shorten <= 0.2cm, shorten >= 0.15cm](1.5,4.35)--(-3.4,4.35);
\draw[-angle 90, line width=0.3mm, >=latex, shorten <= 0.2cm, shorten >= 0.15cm](3.5,3.65)--(-1.5,3.65);
\draw[-angle 90, line width=0.3mm, >=latex, shorten <= 0.2cm, shorten >= 0.15cm](3.5,3.35)--(-1.5,3.35);
\draw[-angle 90, line width=0.3mm, >=latex, shorten <= 0.2cm, shorten >= 0.15cm](-3.35,4.5)--(-2.45,5.8);
\draw[-angle 90, line width=0.3mm, >=latex, shorten <= 0.2cm, shorten >= 0.15cm](-3.65,4.5)--(-2.75,5.8);
\draw[-angle 90, line width=0.3mm, >=latex, shorten <= 0.2cm, shorten >= 0.15cm](-1.2,3.5)--(-2.25,5.8);
\draw[-angle 90, line width=0.3mm, >=latex, shorten <= 0.2cm, shorten >= 0.15cm](-1.5,3.5)--(-2.55,5.8);
\draw[-angle 90, line width=0.3mm, >=latex, shorten <= 0.2cm, shorten >= 0.15cm](2.25,5.8)--(1.25,4.5);
\draw[-angle 90, line width=0.3mm, >=latex, shorten <= 0.2cm, shorten >= 0.15cm](2.55,5.8)--(1.55,4.5);
\draw[-angle 90, line width=0.3mm, >=latex, shorten <= 0.2cm, shorten >= 0.15cm](2.45,5.8)--(3.45,3.5);
\draw[-angle 90, line width=0.3mm, >=latex, shorten <= 0.2cm, shorten >= 0.15cm](2.75,5.8)--(3.75,3.5);
\draw[-angle 90, line width=0.3mm, >=latex, shorten <= 0.2cm, shorten >= 0.15cm](-2.5,6.15)--(2.5,6.15);
\draw[-angle 90, line width=0.3mm, >=latex, shorten <= 0.2cm, shorten >= 0.15cm](-2.5,5.85)--(2.5,5.85);
\draw[-angle 90, line width=0.3mm, >=latex, shorten <= 0.2cm, shorten >= 0.15cm](-2.5,8.15)--(2.5,8.15);
\draw[-angle 90, line width=0.3mm, >=latex, shorten <= 0.2cm, shorten >= 0.15cm](-2.5,7.85)--(2.5,7.85);
\draw[-angle 90, line width=0.3mm, >=latex, shorten <= 0.2cm, shorten >= 0.15cm](3.5,9.65)--(-1.5,9.65);
\draw[-angle 90, line width=0.3mm, >=latex, shorten <= 0.2cm, shorten >= 0.15cm](3.5,9.35)--(-1.5,9.35);
\draw[-angle 90, line width=0.3mm, >=latex, shorten <= 0.2cm, shorten >= 0.15cm](1.5,10.65)--(-3.5,10.65);
\draw[-angle 90, line width=0.3mm, >=latex, shorten <= 0.2cm, shorten >= 0.15cm](1.5,10.35)--(-3.5,10.35);
\draw[-angle 90, line width=0.3mm, >=latex, shorten <= 0.2cm, shorten >= 0.15cm](-1.2,9.5)--(-2.25,8.1);
\draw[-angle 90, line width=0.3mm, >=latex, shorten <= 0.2cm, shorten >= 0.15cm](-1.5,9.5)--(-2.55,8.1);
\draw[-angle 90, line width=0.3mm, >=latex, shorten <= 0.2cm, shorten >= 0.15cm](-3.35,10.35)--(-2.45,8.1);
\draw[-angle 90, line width=0.3mm, >=latex, shorten <= 0.2cm, shorten >= 0.15cm](-3.65,10.35)--(-2.75,8.1);
\draw[-angle 90, line width=0.3mm, >=latex, shorten <= 0.2cm, shorten >= 0.15cm](2.25,8.2)--(1.35,10.35);
\draw[-angle 90, line width=0.3mm, >=latex, shorten <= 0.2cm, shorten >= 0.15cm](2.55,8.2)--(1.65,10.35);
\draw[-angle 90, line width=0.3mm, >=latex, shorten <= 0.2cm, shorten >= 0.15cm](2.45,8.2)--(3.55,9.45);
\draw[-angle 90, line width=0.3mm, >=latex, shorten <= 0.2cm, shorten >= 0.15cm](2.75,8.2)--(3.85,9.45);

\node (e1) [E, minimum height=6.6cm, minimum width=2.8cm,  label={[yshift=-0.2cm, xshift=0cm]270:{$[1]_1$}}] at (-7.5,7){};
\node (e2) [E, minimum height=6.6cm, minimum width=2.8cm,  label={[yshift=-0.4cm, xshift=0cm]270:{$\mathtt{c}_1$}}] at (-2.5,7){};
\node (e3) [E, minimum height=6.6cm, minimum width=2.8cm,  label={[yshift=-0.4cm, xshift=0cm]270:{$\mathtt{c}_2$}}] at (2.5,7){};
\node (e4) [E, minimum height=6.6cm, minimum width=2.8cm,  label={[yshift=-0.2cm, xshift=0cm]270:{$[1]_2$}}] at (7.5,7){};
\end{tikzpicture}
}
{
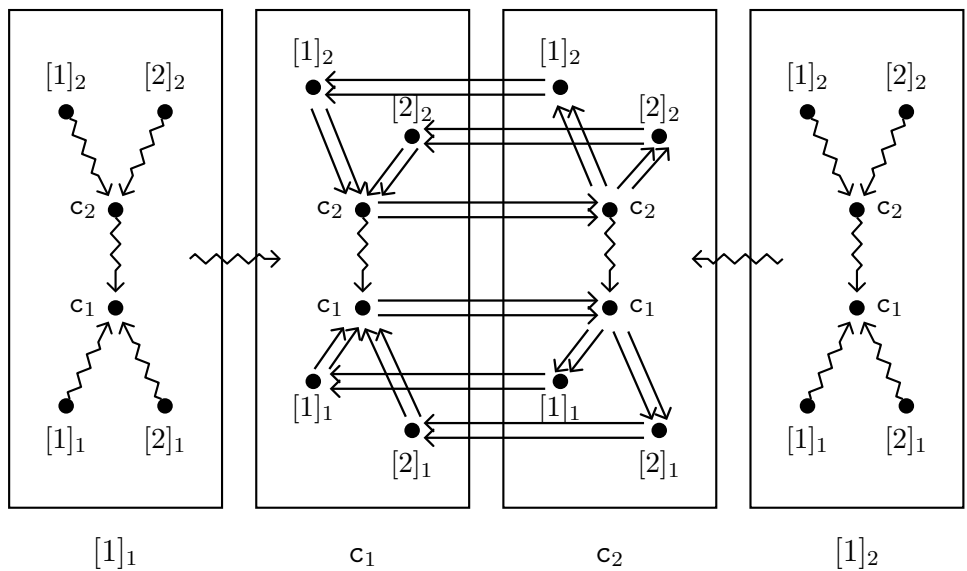
\captionof{figure}{Orientation $D_{(3,3)}$ of $(T_3\times T_3)^{(2)}$, where $d(D_{(3,3)})=6$.}\label{figA8.3.10}}
~\\
\noindent\makebox[\textwidth]{%
\tikzstyle{every node}=[circle, draw, fill=black!100,
                       inner sep=0pt, minimum width=5pt]
\begin{tikzpicture}[thick,scale=0.65]%
\draw(-6.5,12)node[label={[yshift=-0.2cm, xshift=0cm]90:{\small $[1,2]_2$}}](11_122){};
\draw(-6.5,10)node[label={[yshift=0cm, xshift=0cm]315:{\small $[2]_2$}}](11_22){};
\draw(-6.5,2)node[label={[yshift=0.2cm, xshift=0cm]270:{\small $[1,2]_1$}}](11_121){};
\draw(-6.5,4)node[label={[yshift=0cm, xshift=0cm]45:{\small $[2]_1$}}](11_21){};

\draw(-8.5,12)node[label={[yshift=-0.2cm, xshift=0cm]90:{\small $[1,1]_2$}}](11_112){};
\draw(-8.5,10)node[label={[yshift=0cm, xshift=0cm]215:{\small $[1]_2$}}](11_12){};
\draw(-7.5,8)node[label={[yshift=0cm, xshift=-0.1cm]180:{\small $\mathtt{c}_2$}}](11_c2){};
\draw(-7.5,6)node[label={[yshift=0cm, xshift=-0.1cm]180:{\small $\mathtt{c}_1$}}](11_c1){};
\draw(-8.5,4)node[label={[yshift=0cm, xshift=0cm]135:{\small $[1]_1$}}](11_11){};
\draw(-8.5,2)node[label={[yshift=0.2cm, xshift=0cm]270:{\small $[1,1]_1$}}](11_111){};
\draw(-1.5,12)node[label={[yshift=-0.2cm, xshift=0cm]90:{\small $[1,2]_2$}}](c1_122){};
\draw(-1.5,9.5)node[label={[yshift=-0.5cm, xshift=0cm]0:{\small $[2]_2$}}](c1_22){};
\draw(-1.5,2)node[label={[yshift=0.2cm, xshift=0cm]270:{\small $[1,2]_1$}}](c1_121){};
\draw(-1.5,3.5)node[label={[yshift=-0.5cm, xshift=0cm]0:{\small $[2]_1$}}](c1_21){};

\draw(-3.5,12)node[label={[yshift=-0.2cm, xshift=0cm]90:{\small $[1,1]_2$}}](c1_112){};
\draw(-3.5,10.5)node[label={[yshift=0cm, xshift=0cm]180:{\small $[1]_2$}}](c1_12){};
\draw(-2.5,8)node[label={[yshift=0cm, xshift=-0.1cm]180:{\small $\mathtt{c}_2$}}](c1_c2){};
\draw(-2.5,6)node[label={[yshift=0cm, xshift=-0.1cm]180:{\small $\mathtt{c}_1$}}](c1_c1){};
\draw(-3.5,4.5)node[label={[yshift=0cm, xshift=0cm]180:{\small $[1]_1$}}](c1_11){};
\draw(-3.5,2)node[label={[yshift=0.2cm, xshift=0cm]270:{\small $[1,1]_1$}}](c1_111){};
\draw(3.5,12)node[label={[yshift=-0.2cm, xshift=0cm]90:{\small $[1,2]_2$}}](c2_122){};
\draw(3.5,9.5)node[label={[yshift=0cm, xshift=0cm]0:{\small $[2]_2$}}](c2_22){};
\draw(3.5,2)node[label={[yshift=0.2cm, xshift=0cm]270:{\small $[1,2]_1$}}](c2_121){};
\draw(3.5,3.5)node[label={[yshift=0cm, xshift=0cm]0:{\small $[2]_1$}}](c2_21){};

\draw(1.5,12)node[label={[yshift=-0.2cm, xshift=0cm]90:{\small $[1,1]_2$}}](c2_112){};
\draw(1.5,10.5)node[label={[yshift=0.4cm, xshift=0cm]180:{\small $[1]_2$}}](c2_12){};
\draw(2.5,8)node[label={[yshift=0cm, xshift=0.1cm]0:{\small $\mathtt{c}_2$}}](c2_c2){};
\draw(2.5,6)node[label={[yshift=0cm, xshift=0.1cm]0:{\small $\mathtt{c}_1$}}](c2_c1){};
\draw(1.5,4.5)node[label={[yshift=0.4cm, xshift=0cm]180:{\small $[1]_1$}}](c2_11){};
\draw(1.5,2)node[label={[yshift=0.2cm, xshift=0cm]270:{\small $[1,1]_1$}}](c2_111){};
\draw(8.5,12)node[label={[yshift=-0.2cm, xshift=0cm]90:{\small $[1,2]_2$}}](12_122){};
\draw(8.5,10)node[label={[yshift=0cm, xshift=0cm]315:{\small $[2]_2$}}](12_22){};
\draw(8.5,2)node[label={[yshift=0.2cm, xshift=0cm]270:{\small $[1,2]_1$}}](12_121){};
\draw(8.5,4)node[label={[yshift=0cm, xshift=0cm]45:{\small $[2]_1$}}](12_21){};

\draw(6.5,12)node[label={[yshift=-0.2cm, xshift=0cm]90:{\small $[1,1]_2$}}](12_112){};
\draw(6.5,10)node[label={[yshift=0cm, xshift=0cm]215:{\small $[1]_2$}}](12_12){};
\draw(7.5,8)node[label={[yshift=0cm, xshift=0.1cm]0:{\small $\mathtt{c}_2$}}](12_c2){};
\draw(7.5,6)node[label={[yshift=0cm, xshift=0.1cm]0:{\small $\mathtt{c}_1$}}](12_c1){};
\draw(6.5,4)node[label={[yshift=0cm, xshift=0cm]135:{\small $[1]_1$}}](12_11){};
\draw(6.5,2)node[label={[yshift=0.2cm, xshift=0cm]270:{\small $[1,1]_1$}}](12_111){};
\draw [-angle 90, shorten >= 0.1cm, line join=round,decorate, decoration={zigzag, segment length=8,amplitude=1.8,post=lineto,post length=8pt}](11_111)--(11_11);
\draw [-angle 90, shorten >= 0.1cm, line join=round,decorate, decoration={zigzag, segment length=8,amplitude=1.8,post=lineto,post length=8pt}](11_121)--(11_21);
\draw [-angle 90, shorten >= 0.1cm, line join=round,decorate, decoration={zigzag, segment length=8,amplitude=1.8,post=lineto,post length=8pt}](11_11)--(11_c1);
\draw [-angle 90, shorten >= 0.1cm, line join=round,decorate, decoration={zigzag, segment length=8,amplitude=1.8,post=lineto,post length=8pt}](11_21)--(11_c1);

\draw [-angle 90, shorten >= 0.1cm, line join=round,decorate, decoration={zigzag, segment length=8,amplitude=1.8,post=lineto,post length=8pt}](11_112)--(11_12);
\draw [-angle 90, shorten >= 0.1cm, line join=round,decorate, decoration={zigzag, segment length=8,amplitude=1.8,post=lineto,post length=8pt}](11_122)--(11_22);
\draw [-angle 90, shorten >= 0.1cm, line join=round,decorate, decoration={zigzag, segment length=8,amplitude=1.8,post=lineto,post length=8pt}](11_12)--(11_c2);
\draw [-angle 90, shorten >= 0.1cm, line join=round,decorate, decoration={zigzag, segment length=8,amplitude=1.8,post=lineto,post length=8pt}](11_22)--(11_c2);

\draw [-angle 90, shorten >= 0.1cm, line join=round,decorate, decoration={zigzag, segment length=8,amplitude=1.8,post=lineto,post length=8pt}](11_c2)--(11_c1);
\draw [-angle 90, shorten >= 0.1cm, line join=round,decorate, decoration={zigzag, segment length=8,amplitude=1.8,post=lineto,post length=8pt}](c1_111)--(c1_11);
\draw [-angle 90, shorten >= 0.1cm, line join=round,decorate, decoration={zigzag, segment length=8,amplitude=1.8,post=lineto,post length=8pt}](c1_121)--(c1_21);

\draw [-angle 90, shorten >= 0.1cm, line join=round,decorate, decoration={zigzag, segment length=8,amplitude=1.8,post=lineto,post length=8pt}](c1_112)--(c1_12);
\draw [-angle 90, shorten >= 0.1cm, line join=round,decorate, decoration={zigzag, segment length=8,amplitude=1.8,post=lineto,post length=8pt}](c1_122)--(c1_22);

\draw [-angle 90, shorten >= 0.1cm, line join=round,decorate, decoration={zigzag, segment length=8,amplitude=1.8,post=lineto,post length=8pt}](c1_c2)--(c1_c1);
\draw [-angle 90, shorten >= 0.1cm, line join=round,decorate, decoration={zigzag, segment length=8,amplitude=1.8,post=lineto,post length=8pt}](c2_111)--(c2_11);
\draw [-angle 90, shorten >= 0.1cm, line join=round,decorate, decoration={zigzag, segment length=8,amplitude=1.8,post=lineto,post length=8pt}](c2_121)--(c2_21);

\draw [-angle 90, shorten >= 0.1cm, line join=round,decorate, decoration={zigzag, segment length=8,amplitude=1.8,post=lineto,post length=8pt}](c2_112)--(c2_12);
\draw [-angle 90, shorten >= 0.1cm, line join=round,decorate, decoration={zigzag, segment length=8,amplitude=1.8,post=lineto,post length=8pt}](c2_122)--(c2_22);

\draw [-angle 90, shorten >= 0.1cm, line join=round,decorate, decoration={zigzag, segment length=8,amplitude=1.8,post=lineto,post length=8pt}](c2_c2)--(c2_c1);
\draw [-angle 90, shorten >= 0.1cm, line join=round,decorate, decoration={zigzag, segment length=8,amplitude=1.8,post=lineto,post length=8pt}](12_111)--(12_11);
\draw [-angle 90, shorten >= 0.1cm, line join=round,decorate, decoration={zigzag, segment length=8,amplitude=1.8,post=lineto,post length=8pt}](12_121)--(12_21);
\draw [-angle 90, shorten >= 0.1cm, line join=round,decorate, decoration={zigzag, segment length=8,amplitude=1.8,post=lineto,post length=8pt}](12_11)--(12_c1);
\draw [-angle 90, shorten >= 0.1cm, line join=round,decorate, decoration={zigzag, segment length=8,amplitude=1.8,post=lineto,post length=8pt}](12_21)--(12_c1);

\draw [-angle 90, shorten >= 0.1cm, line join=round,decorate, decoration={zigzag, segment length=8,amplitude=1.8,post=lineto,post length=8pt}](12_112)--(12_12);
\draw [-angle 90, shorten >= 0.1cm, line join=round,decorate, decoration={zigzag, segment length=8,amplitude=1.8,post=lineto,post length=8pt}](12_122)--(12_22);
\draw [-angle 90, shorten >= 0.1cm, line join=round,decorate, decoration={zigzag, segment length=8,amplitude=1.8,post=lineto,post length=8pt}](12_12)--(12_c2);
\draw [-angle 90, shorten >= 0.1cm, line join=round,decorate, decoration={zigzag, segment length=8,amplitude=1.8,post=lineto,post length=8pt}](12_22)--(12_c2);

\draw [-angle 90, shorten >= 0.1cm, line join=round,decorate, decoration={zigzag, segment length=8,amplitude=1.8,post=lineto,post length=8pt}](12_c2)--(12_c1);
\draw [-angle 90, shorten >= 0.1cm, line join=round,decorate, decoration={zigzag, segment length=8,amplitude=1.8,post=lineto,post length=8pt}](1,12)--(-1,12);
\draw [-angle 90, shorten >= 0.1cm, line join=round,decorate, decoration={zigzag, segment length=8,amplitude=1.8,post=lineto,post length=8pt}](-6,7)--(-4,7);
\draw [-angle 90, shorten >= 0.1cm, line join=round,decorate, decoration={zigzag, segment length=8,amplitude=1.8,post=lineto,post length=8pt}](6,7)--(4,7);
\draw [-angle 90, shorten >= 0.1cm, line join=round,decorate, decoration={zigzag, segment length=8,amplitude=1.8,post=lineto,post length=8pt}](1,2)--(-1,2);
\draw[-angle 90, line width=0.3mm, >=latex, shorten <= 0.2cm, shorten >= 0.15cm](1.5,4.65)--(-3.4,4.65);
\draw[-angle 90, line width=0.3mm, >=latex, shorten <= 0.2cm, shorten >= 0.15cm](1.5,4.35)--(-3.4,4.35);
\draw[-angle 90, line width=0.3mm, >=latex, shorten <= 0.2cm, shorten >= 0.15cm](3.5,3.65)--(-1.5,3.65);
\draw[-angle 90, line width=0.3mm, >=latex, shorten <= 0.2cm, shorten >= 0.15cm](3.5,3.35)--(-1.5,3.35);
\draw[-angle 90, line width=0.3mm, >=latex, shorten <= 0.2cm, shorten >= 0.15cm](-3.35,4.5)--(-2.45,5.8);
\draw[-angle 90, line width=0.3mm, >=latex, shorten <= 0.2cm, shorten >= 0.15cm](-3.65,4.5)--(-2.75,5.8);
\draw[-angle 90, line width=0.3mm, >=latex, shorten <= 0.2cm, shorten >= 0.15cm](-1.2,3.5)--(-2.25,5.8);
\draw[-angle 90, line width=0.3mm, >=latex, shorten <= 0.2cm, shorten >= 0.15cm](-1.5,3.5)--(-2.55,5.8);
\draw[-angle 90, line width=0.3mm, >=latex, shorten <= 0.2cm, shorten >= 0.15cm](2.25,5.8)--(1.25,4.5);
\draw[-angle 90, line width=0.3mm, >=latex, shorten <= 0.2cm, shorten >= 0.15cm](2.55,5.8)--(1.55,4.5);
\draw[-angle 90, line width=0.3mm, >=latex, shorten <= 0.2cm, shorten >= 0.15cm](2.45,5.8)--(3.45,3.5);
\draw[-angle 90, line width=0.3mm, >=latex, shorten <= 0.2cm, shorten >= 0.15cm](2.75,5.8)--(3.75,3.5);
\draw[-angle 90, line width=0.3mm, >=latex, shorten <= 0.2cm, shorten >= 0.15cm](-2.5,6.15)--(2.5,6.15);
\draw[-angle 90, line width=0.3mm, >=latex, shorten <= 0.2cm, shorten >= 0.15cm](-2.5,5.85)--(2.5,5.85);
\draw[-angle 90, line width=0.3mm, >=latex, shorten <= 0.2cm, shorten >= 0.15cm](-2.5,8.15)--(2.5,8.15);
\draw[-angle 90, line width=0.3mm, >=latex, shorten <= 0.2cm, shorten >= 0.15cm](-2.5,7.85)--(2.5,7.85);
\draw[-angle 90, line width=0.3mm, >=latex, shorten <= 0.2cm, shorten >= 0.15cm](3.5,9.65)--(-1.5,9.65);
\draw[-angle 90, line width=0.3mm, >=latex, shorten <= 0.2cm, shorten >= 0.15cm](3.5,9.35)--(-1.5,9.35);
\draw[-angle 90, line width=0.3mm, >=latex, shorten <= 0.2cm, shorten >= 0.15cm](1.5,10.65)--(-3.5,10.65);
\draw[-angle 90, line width=0.3mm, >=latex, shorten <= 0.2cm, shorten >= 0.15cm](1.5,10.35)--(-3.5,10.35);
\draw[-angle 90, line width=0.3mm, >=latex, shorten <= 0.2cm, shorten >= 0.15cm](-1.2,9.5)--(-2.25,8.1);
\draw[-angle 90, line width=0.3mm, >=latex, shorten <= 0.2cm, shorten >= 0.15cm](-1.5,9.5)--(-2.55,8.1);
\draw[-angle 90, line width=0.3mm, >=latex, shorten <= 0.2cm, shorten >= 0.15cm](-3.35,10.35)--(-2.45,8.1);
\draw[-angle 90, line width=0.3mm, >=latex, shorten <= 0.2cm, shorten >= 0.15cm](-3.65,10.35)--(-2.75,8.1);
\draw[-angle 90, line width=0.3mm, >=latex, shorten <= 0.2cm, shorten >= 0.15cm](2.25,8.2)--(1.35,10.35);
\draw[-angle 90, line width=0.3mm, >=latex, shorten <= 0.2cm, shorten >= 0.15cm](2.55,8.2)--(1.65,10.35);
\draw[-angle 90, line width=0.3mm, >=latex, shorten <= 0.2cm, shorten >= 0.15cm](2.45,8.2)--(3.55,9.45);
\draw[-angle 90, line width=0.3mm, >=latex, shorten <= 0.2cm, shorten >= 0.15cm](2.75,8.2)--(3.85,9.45);

\node (e1) [E, minimum height=8.4cm, minimum width=3cm,  label={[yshift=-0.2cm, xshift=0cm]270:{$[1]_1$}}] at (-7.5,7){};
\node (e2) [E, minimum height=8.4cm, minimum width=3cm,  label={[yshift=-0.4cm, xshift=0cm]270:{$\mathtt{c}_1$}}] at (-2.5,7){};
\node (e3) [E, minimum height=8.4cm, minimum width=3cm,  label={[yshift=-0.3cm, xshift=0cm]270:{$\mathtt{c}_2$}}] at (2.5,7){};
\node (e4) [E, minimum height=8.4cm, minimum width=3cm,  label={[yshift=-0.2cm, xshift=0cm]270:{$[1]_2$}}] at (7.5,7){};
\end{tikzpicture}
}
{\captionof{figure}{Orientation $D_{(3,5)}$ of $(T_3\times T_5)^{(2)}$, where $d(D_{(3,5)})=8$.}\label{figA8.3.11}}
\enlargethispage{1\baselineskip}
\begin{sidewaysfigure}
\centering
\noindent\makebox[\textwidth]{%
\tikzstyle{every node}=[circle, draw, fill=black!100,
                       inner sep=0pt, minimum width=5pt]
\begin{tikzpicture}[thick,scale=0.75]%
\draw(-11.5,12)node[label={[yshift=-0.2cm, xshift=0cm]90:{\small $[1,2]_2$}}](111_122){};
\draw(-11.5,10)node[label={[yshift=0cm, xshift=0cm]315:{\small $[2]_2$}}](111_22){};
\draw(-11.5,2)node[label={[yshift=0.2cm, xshift=0cm]270:{\small $[1,2]_1$}}](111_121){};
\draw(-11.5,4)node[label={[yshift=0cm, xshift=0cm]45:{\small $[2]_1$}}](111_21){};

\draw(-13.5,12)node[label={[yshift=-0.2cm, xshift=0cm]90:{\small $[1,1]_2$}}](111_112){};
\draw(-13.5,10)node[label={[yshift=0cm, xshift=0cm]215:{\small $[1]_2$}}](111_12){};
\draw(-12.5,8)node[label={[yshift=0cm, xshift=-0.1cm]180:{\small $\mathtt{c}_2$}}](111_c2){};
\draw(-12.5,6)node[label={[yshift=0cm, xshift=-0.1cm]180:{\small $\mathtt{c}_1$}}](111_c1){};
\draw(-13.5,4)node[label={[yshift=0cm, xshift=0cm]135:{\small $[1]_1$}}](111_11){};
\draw(-13.5,2)node[label={[yshift=0.2cm, xshift=0cm]270:{\small $[1,1]_1$}}](111_111){};
\draw(-6.5,12)node[label={[yshift=-0.2cm, xshift=0cm]90:{\small $[1,2]_2$}}](11_122){};
\draw(-6.5,10)node[label={[yshift=0cm, xshift=0cm]315:{\small $[2]_2$}}](11_22){};
\draw(-6.5,2)node[label={[yshift=0.2cm, xshift=0cm]270:{\small $[1,2]_1$}}](11_121){};
\draw(-6.5,4)node[label={[yshift=0cm, xshift=0cm]45:{\small $[2]_1$}}](11_21){};

\draw(-8.5,12)node[label={[yshift=-0.2cm, xshift=0cm]90:{\small $[1,1]_2$}}](11_112){};
\draw(-8.5,10)node[label={[yshift=0cm, xshift=0cm]215:{\small $[1]_2$}}](11_12){};
\draw(-7.5,8)node[label={[yshift=0cm, xshift=-0.1cm]180:{\small $\mathtt{c}_2$}}](11_c2){};
\draw(-7.5,6)node[label={[yshift=0cm, xshift=-0.1cm]180:{\small $\mathtt{c}_1$}}](11_c1){};
\draw(-8.5,4)node[label={[yshift=0cm, xshift=0cm]135:{\small $[1]_1$}}](11_11){};
\draw(-8.5,2)node[label={[yshift=0.2cm, xshift=0cm]270:{\small $[1,1]_1$}}](11_111){};
\draw(-1.5,12)node[label={[yshift=-0.2cm, xshift=0cm]90:{\small $[1,2]_2$}}](c1_122){};
\draw(-1.5,9.5)node[label={[yshift=0.4cm, xshift=0.12cm]0:{\small $[2]_2$}}](c1_22){};
\draw(-1.5,2)node[label={[yshift=0.2cm, xshift=0cm]270:{\small $[1,2]_1$}}](c1_121){};
\draw(-1.5,3.5)node[label={[yshift=0.4cm, xshift=0cm]0:{\small $[2]_1$}}](c1_21){};

\draw(-3.5,12)node[label={[yshift=-0.2cm, xshift=0cm]90:{\small $[1,1]_2$}}](c1_112){};
\draw(-3.5,10.5)node[label={[yshift=0cm, xshift=0cm]180:{\small $[1]_2$}}](c1_12){};
\draw(-2.5,8)node[label={[yshift=0cm, xshift=-0.1cm]180:{\small $\mathtt{c}_2$}}](c1_c2){};
\draw(-2.5,6)node[label={[yshift=0cm, xshift=-0.1cm]180:{\small $\mathtt{c}_1$}}](c1_c1){};
\draw(-3.5,4.5)node[label={[yshift=0cm, xshift=0cm]180:{\small $[1]_1$}}](c1_11){};
\draw(-3.5,2)node[label={[yshift=0.2cm, xshift=0cm]270:{\small $[1,1]_1$}}](c1_111){};
\draw(3.5,12)node[label={[yshift=-0.2cm, xshift=0cm]90:{\small $[1,2]_2$}}](c2_122){};
\draw(3.5,9.5)node[label={[yshift=0cm, xshift=0cm]0:{\small $[2]_2$}}](c2_22){};
\draw(3.5,2)node[label={[yshift=0.2cm, xshift=0cm]270:{\small $[1,2]_1$}}](c2_121){};
\draw(3.5,3.5)node[label={[yshift=0cm, xshift=0cm]0:{\small $[2]_1$}}](c2_21){};

\draw(1.5,12)node[label={[yshift=-0.2cm, xshift=0cm]90:{\small $[1,1]_2$}}](c2_112){};
\draw(1.5,10.5)node[label={[yshift=0.4cm, xshift=-0.07cm]180:{\small $[1]_2$}}](c2_12){};
\draw(2.5,8)node[label={[yshift=0cm, xshift=0.1cm]0:{\small $\mathtt{c}_2$}}](c2_c2){};
\draw(2.5,6)node[label={[yshift=0cm, xshift=0.1cm]0:{\small $\mathtt{c}_1$}}](c2_c1){};
\draw(1.5,4.5)node[label={[yshift=0.4cm, xshift=0cm]180:{\small $[1]_1$}}](c2_11){};
\draw(1.5,2)node[label={[yshift=0.2cm, xshift=0cm]270:{\small $[1,1]_1$}}](c2_111){};
\draw(8.5,12)node[label={[yshift=-0.2cm, xshift=0cm]90:{\small $[1,2]_2$}}](12_122){};
\draw(8.5,10)node[label={[yshift=0cm, xshift=0cm]315:{\small $[2]_2$}}](12_22){};
\draw(8.5,2)node[label={[yshift=0.2cm, xshift=0cm]270:{\small $[1,2]_1$}}](12_121){};
\draw(8.5,4)node[label={[yshift=0cm, xshift=0cm]45:{\small $[2]_1$}}](12_21){};

\draw(6.5,12)node[label={[yshift=-0.2cm, xshift=0cm]90:{\small $[1,1]_2$}}](12_112){};
\draw(6.5,10)node[label={[yshift=0cm, xshift=0cm]215:{\small $[1]_2$}}](12_12){};
\draw(7.5,8)node[label={[yshift=0cm, xshift=0.1cm]0:{\small $\mathtt{c}_2$}}](12_c2){};
\draw(7.5,6)node[label={[yshift=0cm, xshift=0.1cm]0:{\small $\mathtt{c}_1$}}](12_c1){};
\draw(6.5,4)node[label={[yshift=0cm, xshift=0cm]135:{\small $[1]_1$}}](12_11){};
\draw(6.5,2)node[label={[yshift=0.2cm, xshift=0cm]270:{\small $[1,1]_1$}}](12_111){};
\draw(13.5,12)node[label={[yshift=-0.2cm, xshift=0cm]90:{\small $[1,2]_2$}}](112_122){};
\draw(13.5,10)node[label={[yshift=0cm, xshift=0cm]315:{\small $[2]_2$}}](112_22){};
\draw(13.5,2)node[label={[yshift=0.2cm, xshift=0cm]270:{\small $[1,2]_1$}}](112_121){};
\draw(13.5,4)node[label={[yshift=0cm, xshift=0cm]45:{\small $[2]_1$}}](112_21){};

\draw(11.5,12)node[label={[yshift=-0.2cm, xshift=0cm]90:{\small $[1,1]_2$}}](112_112){};
\draw(11.5,10)node[label={[yshift=0cm, xshift=0cm]215:{\small $[1]_2$}}](112_12){};
\draw(12.5,8)node[label={[yshift=0cm, xshift=0.1cm]0:{\small $\mathtt{c}_2$}}](112_c2){};
\draw(12.5,6)node[label={[yshift=0cm, xshift=0.1cm]0:{\small $\mathtt{c}_1$}}](112_c1){};
\draw(11.5,4)node[label={[yshift=0cm, xshift=0cm]135:{\small $[1]_1$}}](112_11){};
\draw(11.5,2)node[label={[yshift=0.2cm, xshift=0cm]270:{\small $[1,1]_1$}}](112_111){};
\draw [-angle 90, shorten >= 0.1cm, line join=round,decorate, decoration={zigzag, segment length=8,amplitude=1.8,post=lineto,post length=8pt}](111_112)--(111_12);
\draw [-angle 90, shorten >= 0.1cm, line join=round,decorate, decoration={zigzag, segment length=8,amplitude=1.8,post=lineto,post length=8pt}](111_122)--(111_22);
\draw [-angle 90, shorten >= 0.1cm, line join=round,decorate, decoration={zigzag, segment length=8,amplitude=1.8,post=lineto,post length=8pt}](111_21)--(111_c1);
\draw [-angle 90, shorten >= 0.1cm, line join=round,decorate, decoration={zigzag, segment length=8,amplitude=1.8,post=lineto,post length=8pt}](111_22)--(111_c2);
\draw [-angle 90, shorten >= 0.1cm, line join=round,decorate, decoration={zigzag, segment length=8,amplitude=1.8,post=lineto,post length=8pt}](111_c2)--(111_c1);
\draw [-angle 90, shorten >= 0.1cm, line join=round,decorate, decoration={zigzag, segment length=8,amplitude=1.8,post=lineto,post length=8pt}](111_111)--(111_11);
\draw [-angle 90, shorten >= 0.1cm, line join=round,decorate, decoration={zigzag, segment length=8,amplitude=1.8,post=lineto,post length=8pt}](111_121)--(111_21);
\draw [-angle 90, shorten >= 0.1cm, line join=round,decorate, decoration={zigzag, segment length=8,amplitude=1.8,post=lineto,post length=8pt}](111_11)--(111_c1);
\draw [-angle 90, shorten >= 0.1cm, line join=round,decorate, decoration={zigzag, segment length=8,amplitude=1.8,post=lineto,post length=8pt}](111_12)--(111_c2);
\draw [-angle 90, shorten >= 0.1cm, line join=round,decorate, decoration={zigzag, segment length=8,amplitude=1.8,post=lineto,post length=8pt}](11_112)--(11_12);
\draw [-angle 90, shorten >= 0.1cm, line join=round,decorate, decoration={zigzag, segment length=8,amplitude=1.8,post=lineto,post length=8pt}](11_122)--(11_22);
\draw [-angle 90, shorten >= 0.1cm, line join=round,decorate, decoration={zigzag, segment length=8,amplitude=1.8,post=lineto,post length=8pt}](11_21)--(11_c1);
\draw [-angle 90, shorten >= 0.1cm, line join=round,decorate, decoration={zigzag, segment length=8,amplitude=1.8,post=lineto,post length=8pt}](11_22)--(11_c2);
\draw [-angle 90, shorten >= 0.1cm, line join=round,decorate, decoration={zigzag, segment length=8,amplitude=1.8,post=lineto,post length=8pt}](11_c2)--(11_c1);
\draw [-angle 90, shorten >= 0.1cm, line join=round,decorate, decoration={zigzag, segment length=8,amplitude=1.8,post=lineto,post length=8pt}](11_111)--(11_11);
\draw [-angle 90, shorten >= 0.1cm, line join=round,decorate, decoration={zigzag, segment length=8,amplitude=1.8,post=lineto,post length=8pt}](11_121)--(11_21);
\draw [-angle 90, shorten >= 0.1cm, line join=round,decorate, decoration={zigzag, segment length=8,amplitude=1.8,post=lineto,post length=8pt}](11_11)--(11_c1);
\draw [-angle 90, shorten >= 0.1cm, line join=round,decorate, decoration={zigzag, segment length=8,amplitude=1.8,post=lineto,post length=8pt}](11_12)--(11_c2);
\draw [-angle 90, shorten >= 0.1cm, line join=round,decorate, decoration={zigzag, segment length=8,amplitude=1.8,post=lineto,post length=8pt}](c1_112)--(c1_12);
\draw [-angle 90, shorten >= 0.1cm, line join=round,decorate, decoration={zigzag, segment length=8,amplitude=1.8,post=lineto,post length=8pt}](c1_122)--(c1_22);
\draw [-angle 90, shorten >= 0.1cm, line join=round,decorate, decoration={zigzag, segment length=8,amplitude=1.8,post=lineto,post length=8pt}](c1_c2)--(c1_c1);
\draw [-angle 90, shorten >= 0.1cm, line join=round,decorate, decoration={zigzag, segment length=8,amplitude=1.8,post=lineto,post length=8pt}](c1_111)--(c1_11);
\draw [-angle 90, shorten >= 0.1cm, line join=round,decorate, decoration={zigzag, segment length=8,amplitude=1.8,post=lineto,post length=8pt}](c1_121)--(c1_21);
\draw [-angle 90, shorten >= 0.1cm, line join=round,decorate, decoration={zigzag, segment length=8,amplitude=1.8,post=lineto,post length=8pt}](c2_112)--(c2_12);
\draw [-angle 90, shorten >= 0.1cm, line join=round,decorate, decoration={zigzag, segment length=8,amplitude=1.8,post=lineto,post length=8pt}](c2_122)--(c2_22);
\draw [-angle 90, shorten >= 0.1cm, line join=round,decorate, decoration={zigzag, segment length=8,amplitude=1.8,post=lineto,post length=8pt}](c2_c2)--(c2_c1);
\draw [-angle 90, shorten >= 0.1cm, line join=round,decorate, decoration={zigzag, segment length=8,amplitude=1.8,post=lineto,post length=8pt}](c2_111)--(c2_11);
\draw [-angle 90, shorten >= 0.1cm, line join=round,decorate, decoration={zigzag, segment length=8,amplitude=1.8,post=lineto,post length=8pt}](c2_121)--(c2_21);
\draw [-angle 90, shorten >= 0.1cm, line join=round,decorate, decoration={zigzag, segment length=8,amplitude=1.8,post=lineto,post length=8pt}](12_112)--(12_12);
\draw [-angle 90, shorten >= 0.1cm, line join=round,decorate, decoration={zigzag, segment length=8,amplitude=1.8,post=lineto,post length=8pt}](12_122)--(12_22);
\draw [-angle 90, shorten >= 0.1cm, line join=round,decorate, decoration={zigzag, segment length=8,amplitude=1.8,post=lineto,post length=8pt}](12_21)--(12_c1);
\draw [-angle 90, shorten >= 0.1cm, line join=round,decorate, decoration={zigzag, segment length=8,amplitude=1.8,post=lineto,post length=8pt}](12_22)--(12_c2);
\draw [-angle 90, shorten >= 0.1cm, line join=round,decorate, decoration={zigzag, segment length=8,amplitude=1.8,post=lineto,post length=8pt}](12_c2)--(12_c1);
\draw [-angle 90, shorten >= 0.1cm, line join=round,decorate, decoration={zigzag, segment length=8,amplitude=1.8,post=lineto,post length=8pt}](12_111)--(12_11);
\draw [-angle 90, shorten >= 0.1cm, line join=round,decorate, decoration={zigzag, segment length=8,amplitude=1.8,post=lineto,post length=8pt}](12_121)--(12_21);
\draw [-angle 90, shorten >= 0.1cm, line join=round,decorate, decoration={zigzag, segment length=8,amplitude=1.8,post=lineto,post length=8pt}](12_11)--(12_c1);
\draw [-angle 90, shorten >= 0.1cm, line join=round,decorate, decoration={zigzag, segment length=8,amplitude=1.8,post=lineto,post length=8pt}](12_12)--(12_c2);
\draw [-angle 90, shorten >= 0.1cm, line join=round,decorate, decoration={zigzag, segment length=8,amplitude=1.8,post=lineto,post length=8pt}](112_112)--(112_12);
\draw [-angle 90, shorten >= 0.1cm, line join=round,decorate, decoration={zigzag, segment length=8,amplitude=1.8,post=lineto,post length=8pt}](112_122)--(112_22);
\draw [-angle 90, shorten >= 0.1cm, line join=round,decorate, decoration={zigzag, segment length=8,amplitude=1.8,post=lineto,post length=8pt}](112_21)--(112_c1);
\draw [-angle 90, shorten >= 0.1cm, line join=round,decorate, decoration={zigzag, segment length=8,amplitude=1.8,post=lineto,post length=8pt}](112_22)--(112_c2);
\draw [-angle 90, shorten >= 0.1cm, line join=round,decorate, decoration={zigzag, segment length=8,amplitude=1.8,post=lineto,post length=8pt}](112_c2)--(112_c1);
\draw [-angle 90, shorten >= 0.1cm, line join=round,decorate, decoration={zigzag, segment length=8,amplitude=1.8,post=lineto,post length=8pt}](112_111)--(112_11);
\draw [-angle 90, shorten >= 0.1cm, line join=round,decorate, decoration={zigzag, segment length=8,amplitude=1.8,post=lineto,post length=8pt}](112_121)--(112_21);
\draw [-angle 90, shorten >= 0.1cm, line join=round,decorate, decoration={zigzag, segment length=8,amplitude=1.8,post=lineto,post length=8pt}](112_11)--(112_c1);
\draw [-angle 90, shorten >= 0.1cm, line join=round,decorate, decoration={zigzag, segment length=8,amplitude=1.8,post=lineto,post length=8pt}](112_12)--(112_c2);
\draw [-angle 90, shorten >= 0.1cm, line join=round,decorate, decoration={zigzag, segment length=8,amplitude=1.8,post=lineto,post length=8pt}](-11,7)--(-9,7);
\draw [-angle 90, shorten >= 0.1cm, line join=round,decorate, decoration={zigzag, segment length=8,amplitude=1.8,post=lineto,post length=8pt}](11,7)--(9,7);
\draw [-angle 90, shorten >= 0.1cm, line join=round,decorate, decoration={zigzag, segment length=8,amplitude=1.8,post=lineto,post length=8pt}](-6,7)--(-4,7);
\draw [-angle 90, shorten >= 0.1cm, line join=round,decorate, decoration={zigzag, segment length=8,amplitude=1.8,post=lineto,post length=8pt}](6,7)--(4,7);
\draw [-angle 90, shorten >= 0.1cm, line join=round,decorate, decoration={zigzag, segment length=8,amplitude=1.8,post=lineto,post length=8pt}](1,12)--(-1,12);
\draw [-angle 90, shorten >= 0.1cm, line join=round,decorate, decoration={zigzag, segment length=8,amplitude=1.8,post=lineto,post length=8pt}](1,2)--(-1,2);
\draw[-angle 90, line width=0.3mm, >=latex, shorten <= 0.2cm, shorten >= 0.15cm](1.5,4.65)--(-3.4,4.65);
\draw[-angle 90, line width=0.3mm, >=latex, shorten <= 0.2cm, shorten >= 0.15cm](1.5,4.35)--(-3.4,4.35);
\draw[-angle 90, line width=0.3mm, >=latex, shorten <= 0.2cm, shorten >= 0.15cm](3.5,3.65)--(-1.5,3.65);
\draw[-angle 90, line width=0.3mm, >=latex, shorten <= 0.2cm, shorten >= 0.15cm](3.5,3.35)--(-1.5,3.35);
\draw[-angle 90, line width=0.3mm, >=latex, shorten <= 0.2cm, shorten >= 0.15cm](-3.35,4.5)--(-2.45,5.8);
\draw[-angle 90, line width=0.3mm, >=latex, shorten <= 0.2cm, shorten >= 0.15cm](-3.65,4.5)--(-2.75,5.8);
\draw[-angle 90, line width=0.3mm, >=latex, shorten <= 0.2cm, shorten >= 0.15cm](-1.2,3.5)--(-2.25,5.8);
\draw[-angle 90, line width=0.3mm, >=latex, shorten <= 0.2cm, shorten >= 0.15cm](-1.5,3.5)--(-2.55,5.8);
\draw[-angle 90, line width=0.3mm, >=latex, shorten <= 0.2cm, shorten >= 0.15cm](2.25,5.8)--(1.25,4.5);
\draw[-angle 90, line width=0.3mm, >=latex, shorten <= 0.2cm, shorten >= 0.15cm](2.55,5.8)--(1.55,4.5);
\draw[-angle 90, line width=0.3mm, >=latex, shorten <= 0.2cm, shorten >= 0.15cm](2.45,5.8)--(3.45,3.5);
\draw[-angle 90, line width=0.3mm, >=latex, shorten <= 0.2cm, shorten >= 0.15cm](2.75,5.8)--(3.75,3.5);
\draw[-angle 90, line width=0.3mm, >=latex, shorten <= 0.2cm, shorten >= 0.15cm](-2.5,6.15)--(2.5,6.15);
\draw[-angle 90, line width=0.3mm, >=latex, shorten <= 0.2cm, shorten >= 0.15cm](-2.5,5.85)--(2.5,5.85);
\draw[-angle 90, line width=0.3mm, >=latex, shorten <= 0.2cm, shorten >= 0.15cm](-2.5,8.15)--(2.5,8.15);
\draw[-angle 90, line width=0.3mm, >=latex, shorten <= 0.2cm, shorten >= 0.15cm](-2.5,7.85)--(2.5,7.85);
\draw[-angle 90, line width=0.3mm, >=latex, shorten <= 0.2cm, shorten >= 0.15cm](3.5,9.65)--(-1.5,9.65);
\draw[-angle 90, line width=0.3mm, >=latex, shorten <= 0.2cm, shorten >= 0.15cm](3.5,9.35)--(-1.5,9.35);
\draw[-angle 90, line width=0.3mm, >=latex, shorten <= 0.2cm, shorten >= 0.15cm](1.5,10.65)--(-3.5,10.65);
\draw[-angle 90, line width=0.3mm, >=latex, shorten <= 0.2cm, shorten >= 0.15cm](1.5,10.35)--(-3.5,10.35);
\draw[-angle 90, line width=0.3mm, >=latex, shorten <= 0.2cm, shorten >= 0.15cm](-1.2,9.5)--(-2.25,8.1);
\draw[-angle 90, line width=0.3mm, >=latex, shorten <= 0.2cm, shorten >= 0.15cm](-1.5,9.5)--(-2.55,8.1);
\draw[-angle 90, line width=0.3mm, >=latex, shorten <= 0.2cm, shorten >= 0.15cm](-3.35,10.35)--(-2.45,8.1);
\draw[-angle 90, line width=0.3mm, >=latex, shorten <= 0.2cm, shorten >= 0.15cm](-3.65,10.35)--(-2.75,8.1);
\draw[-angle 90, line width=0.3mm, >=latex, shorten <= 0.2cm, shorten >= 0.15cm](2.25,8.2)--(1.35,10.35);
\draw[-angle 90, line width=0.3mm, >=latex, shorten <= 0.2cm, shorten >= 0.15cm](2.55,8.2)--(1.65,10.35);
\draw[-angle 90, line width=0.3mm, >=latex, shorten <= 0.2cm, shorten >= 0.15cm](2.45,8.2)--(3.55,9.45);
\draw[-angle 90, line width=0.3mm, >=latex, shorten <= 0.2cm, shorten >= 0.15cm](2.75,8.2)--(3.85,9.45);
\node (e1) [E, minimum height=10.2cm, minimum width=3.5cm,  label={[yshift=-0.2cm, xshift=0cm]270:{$[1,1]_1$}}] at (-12.5,7){};
\node (e2) [E, minimum height=10.2cm, minimum width=3.5cm,  label={[yshift=-0.4cm, xshift=0cm]270:{$[1]_1$}}] at (-7.5,7){};
\node (e3) [E, minimum height=10.2cm, minimum width=3.5cm,  label={[yshift=-0.6cm, xshift=0cm]270:{$\mathtt{c}_1$}}] at (-2.5,7){};
\node (e4) [E, minimum height=10.2cm, minimum width=3.5cm,  label={[yshift=-0.6cm, xshift=0cm]270:{$\mathtt{c}_2$}}] at (2.5,7){};
\node (e5) [E, minimum height=10.2cm, minimum width=3.5cm,  label={[yshift=-0.4cm, xshift=0cm]270:{$[1]_2$}}] at (7.5,7){};
\node (e6) [E, minimum height=10.2cm, minimum width=3.5cm,  label={[yshift=-0.2cm, xshift=0cm]270:{$[1,1]_2$}}] at (12.5,7){};
\end{tikzpicture}
}
{\captionof{figure}{Orientation $D_{(5,5)}$ of $(T_5\times T_5)^{(2)}$, where $d(D_{(5,5)})=10$.}\label{figA8.3.12}}
\end{sidewaysfigure}
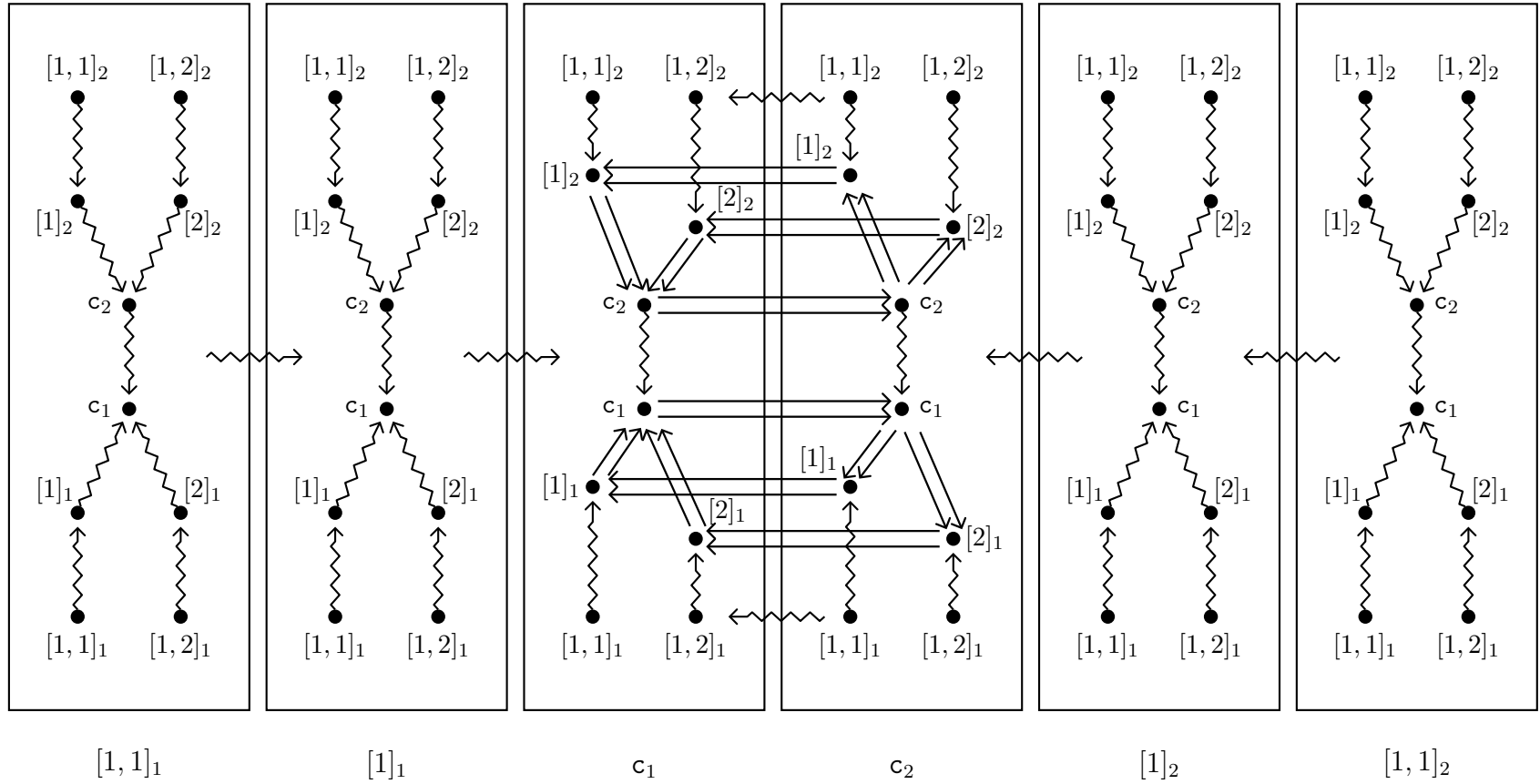
\newpage
\indent\par Next, we shall prove two lemmas for the investigation of the rectangular grid $P_\lambda\times P_\mu$. For $P_n$ ($C_n$ resp.), we shall use the natural labelling of vertices where $E(P_n)=\{ (i,i+1)\mid i=1,2,\ldots,n-1\}$ ($E(C_n)=E(P_n)\cup \{(n,1)\}$ resp.).
\begin{lem}\label{lemA8.3.1}
Let $G$ be a graph and $D$ be an orientation of $G^{(2)}$. If $u_0 u_1 u_2$ is a unique shortest $u_0-u_2$ path in $G$ and $d_D((p,u_0),(q,u_2))=d_D((p,u_2),(q,u_0))=2$ for all $p,q=1,2$, then $u_0\overset{1}\twoheadrightarrow u_1\overset{2}\twoheadleftarrow u_2$ or $u_0\overset{2}\twoheadrightarrow u_1\overset{1}\twoheadleftarrow u_2$.
\end{lem}
\noindent\textit{Proof}: Suppose $(1,u_1)\rightarrow (1,u_2)$. Now, for $p=1,2$, since $d_D((1,u_2),(p,u_0))=2$, it follows that $(1,u_2)\rightarrow (2,u_1)\rightarrow (p,u_0)$. Since $d_D((p,u_0),(q,u_2))=2$ for $p,q=1,2$, it follows that $(p,u_0)\rightarrow (1,u_1)\rightarrow (q,u_2)$ must hold. It is now necessary from $d_D((2,u_2),(1,u_0))=2$ that $(2,u_2)\rightarrow (2,u_1)$. Thus, $u_0\overset{1}\twoheadrightarrow u_1\overset{2}\twoheadleftarrow u_2$. Similarly, an argument reversing all arcs will give $u_0\overset{2}\twoheadrightarrow u_1\overset{1}\twoheadleftarrow u_2$ if $(1,u_2)\rightarrow (1,u_1)$.
\qed
\begin{lem}\label{lemA8.3.2}
Let $G$ be a graph and $D$ be an orientation of $G^{(2)}$. Suppose $v_0v_1 \ldots v_k$, $k\ge 2$, is a shortest $v_0-v_k$ path of length $k$ in $G$ and $D$ satisfies 
\\(a) $v_{i}\overset{1}\twoheadrightarrow v_{i+1}\overset{2}\twoheadleftarrow v_{i+2}$ for some $i$, $0\le i \le k-2$, and
\\(b) if $j\not\in\{i, i+1\}$, then either $v_j\rightsquigarrow v_{j+1}$ or $v_{j+1}\rightsquigarrow v_{j}$.
\\Then, $d_D((p,v_0),(q,v_k))=d_D((p,v_k),(q,v_0))=k$ for $p,q=1,2$.
\end{lem}
\noindent\textit{Proof}: Assume $v_j\rightsquigarrow v_{j+1}$ for all $j\not\in\{i,i+1\}$; the proof is similar otherwise. Note that $(p,v_0)\rightarrow (p,v_1)\rightarrow \cdots \rightarrow (p,v_i)$, $\{(1,v_{i}),(2,v_{i})\}\rightarrow (1,v_{i+1})\rightarrow \{(1,v_{i+2}),(2,v_{i+2})\}$ and $(p,v_{i+2})\rightarrow (p,v_{i+3})\rightarrow \cdots \rightarrow (p,v_k)$ for all $p=1,2$. Thus, for $p,q=1,2$, $d_D((p,v_0),(q,v_k))=d_G(v_0,v_k)=k$. By symmetry, we have $d_D((p,v_k),(q,v_0))=d_G(v_k,v_0)$ $=k$ for $p,q=1,2$.
\qed
\\
\\\textit{Proof of Theorem \ref{thmA8.1.8}}: Let $G=P_\lambda\times P_\mu$ and note that $d(G)=\lambda+\mu-2$.
\\ Case 1. $\lambda=3$ and $\mu=2$.
\indent\par We first prove $\bar{d}(G^{(2)})=4$. Suppose there exists an orientation $D$ of $G^{(2)}$ such that $d(D)=3$. Since $d_D((p,\langle 1,2\rangle),(q,\langle 3,2\rangle))=d_D((q,\langle 3,2\rangle),(p,\langle 1,2\rangle))=2$ for all $p,q=1,2$, we may assume from Lemma \ref{lemA8.3.1} that $\langle 1,2\rangle\overset{1}\twoheadrightarrow \langle 2,2\rangle\overset{2}\twoheadleftarrow \langle 3,2\rangle$. Similarly, we assume $\langle 1,1\rangle \overset{1}\twoheadrightarrow \langle 2,1\rangle \overset{2}\twoheadleftarrow \langle 3,1\rangle$ (the case $\langle 1,1\rangle \overset{2}\twoheadrightarrow \langle 2,1\rangle \overset{1}\twoheadleftarrow \langle 3,1\rangle$ is similar). Since $d_D((1,\langle 1,1\rangle),(2,\langle 2,2\rangle))\le 3$, it follows that $(1,\langle 2,1\rangle)\rightarrow (2,\langle 2,2\rangle)$. However, we have $d_D((1,\langle 3,2\rangle),(1,\langle 2,1\rangle))>3$, which contradicts $d(D)=3$. Hence, $\bar{d}(G^{(2)})\ge 4$.
\indent\par Define an orientation $D_{(3,2)}$ for $G^{(2)}$ as follows. (See Figure \ref{figA8.3.13}.)
\begin{align*}
&\langle 1,j\rangle\overset{1}\twoheadrightarrow \langle 2,j\rangle\overset{2}\twoheadleftarrow \langle 3,j\rangle \text{ for }j=1,2, \text{ and }\langle i,1\rangle\rightsquigarrow \langle i,2\rangle \text{ for }i=1,2,3.
\end{align*}
It is easy to verify $d(D_{(3,2)})=4$. Hence, $G^{(2)}\in\mathscr{C}_1$ and we are done for (a).

\begin{figure}[h]
\begin{center}
\tikzstyle{every node}=[circle, draw, fill=black!100,
                       inner sep=0pt, minimum width=5pt]
\begin{tikzpicture}[thick,scale=0.6]%
\draw(-10,0)node[star,star points=4, label={[yshift=0.5cm] 270:{}}](2_11){};
\draw(-6,0)node[star,star points=4, label={[yshift=0.5cm] 270:{}}](2_21){};
\draw(-2,0)node[star,star points=4, label={[yshift=0.5cm]270:{}}](2_31){};
\draw(-8,2)node[label={[yshift=-0.1cm, xshift=0.1cm]45:{\small $\langle 1,1\rangle$}}](1_11){};
\draw(-4,2)node[label={[yshift=-0.1cm, xshift=0.1cm]45:{\small $\langle 2,1\rangle$}}](1_21){};
\draw(0,2)node[label={[yshift=0.15cm, xshift=0cm]270:{\small $\langle 3,1\rangle$}}](1_31){};

\draw(-10,4)node[star,star points=4, label={[yshift=-0.3cm, xshift=0.3cm] 90:{}}](2_12){};
\draw(-6,4)node[star,star points=4, label={[yshift=-0.3cm, xshift=0.3cm] 90:{}}](2_22){};
\draw(-2,4)node[star,star points=4, label={[yshift=-0.3cm, xshift=0.3cm] 90:{}}](2_32){};
\draw(-8,6)node[label={[yshift=-0.2cm]90:{\small $\langle 1,2\rangle$}}](1_12){};
\draw(-4,6)node[label={[yshift=-0.2cm]90:{\small $\langle 2,2\rangle$}}](1_22){};
\draw(0,6)node[label={[yshift=-0.2cm]90:{\small $\langle 3,2\rangle$}}](1_32){};

\draw[->, line width=0.3mm, >=latex, shorten <= 0.2cm, shorten >= 0.15cm](1_11)--(1_21);
\draw[->, line width=0.3mm, >=latex, shorten <= 0.2cm, shorten >= 0.15cm](2_11)--(1_21);
\draw[->, line width=0.3mm, >=latex, shorten <= 0.2cm, shorten >= 0.15cm](1_21)--(1_31);
\draw[->, line width=0.3mm, >=latex, shorten <= 0.2cm, shorten >= 0.15cm](1_21)--(2_31);
\draw[->, line width=0.3mm, >=latex, shorten <= 0.2cm, shorten >= 0.15cm](2_21)--(1_11);
\draw[->, line width=0.3mm, >=latex, shorten <= 0.2cm, shorten >= 0.15cm](2_21)--(2_11);
\draw[->, line width=0.3mm, >=latex, shorten <= 0.2cm, shorten >= 0.15cm](2_31)--(2_21);
\draw[->, line width=0.3mm, >=latex, shorten <= 0.2cm, shorten >= 0.15cm](1_31)--(2_21);

\draw[->, line width=0.3mm, >=latex, shorten <= 0.2cm, shorten >= 0.15cm](1_12)--(1_22);
\draw[->, line width=0.3mm, >=latex, shorten <= 0.2cm, shorten >= 0.15cm](2_12)--(1_22);
\draw[->, line width=0.3mm, >=latex, shorten <= 0.2cm, shorten >= 0.15cm](1_22)--(1_32);
\draw[->, line width=0.3mm, >=latex, shorten <= 0.2cm, shorten >= 0.15cm](1_22)--(2_32);
\draw[->, line width=0.3mm, >=latex, shorten <= 0.2cm, shorten >= 0.15cm](2_22)--(1_12);
\draw[->, line width=0.3mm, >=latex, shorten <= 0.2cm, shorten >= 0.15cm](2_22)--(2_12);
\draw[->, line width=0.3mm, >=latex, shorten <= 0.2cm, shorten >= 0.15cm](2_32)--(2_22);
\draw[->, line width=0.3mm, >=latex, shorten <= 0.2cm, shorten >= 0.15cm](1_32)--(2_22);

\draw[dashed, ->, line width=0.3mm, >=latex, shorten <= 0.2cm, shorten >= 0.15cm](1_11)--(1_12);
\draw[dashed, ->, line width=0.3mm, >=latex, shorten <= 0.2cm, shorten >= 0.15cm](2_11)--(2_12);
\draw[dashed, ->, line width=0.3mm, >=latex, shorten <= 0.2cm, shorten >= 0.15cm](1_21)--(1_22);
\draw[dashed, ->, line width=0.3mm, >=latex, shorten <= 0.2cm, shorten >= 0.15cm](2_21)--(2_22);
\draw[dashed, ->, line width=0.3mm, >=latex, shorten <= 0.2cm, shorten >= 0.15cm](1_31)--(1_32);
\draw[dashed, ->, line width=0.3mm, >=latex, shorten <= 0.2cm, shorten >= 0.15cm](2_31)--(2_32);

\draw[dashed, ->, line width=0.3mm, >=latex, shorten <= 0.2cm, shorten >= 0.15cm](1_12)--(2_11);
\draw[dashed, ->, line width=0.3mm, >=latex, shorten <= 0.2cm, shorten >= 0.15cm](2_12)--(1_11);
\draw[dashed, ->, line width=0.3mm, >=latex, shorten <= 0.2cm, shorten >= 0.15cm](1_32)--(2_31);
\draw[dashed, ->, line width=0.3mm, >=latex, shorten <= 0.2cm, shorten >= 0.15cm](2_32)--(1_31);
\draw[dashed, ->, line width=0.3mm, >=latex, shorten <= 0.2cm, shorten >= 0.15cm](1_22)--(2_21);
\draw[dashed, ->, line width=0.3mm, >=latex, shorten <= 0.2cm, shorten >= 0.15cm](2_22)--(1_21);
\end{tikzpicture}
{\captionof{figure}{Orientation $D_{(3,2)}$ of $(P_3\times P_2)^{(2)}$, where $d(D_{(3,2)})=4$.}\label{figA8.3.13}}
\caption*{Note: The vertices $(p,\langle u,x\rangle)$, for $p=1,2$, are represented by $\bullet$ and $\blackdiamond$ respectively. The vertex $(1,\langle u,x\rangle)$ is simply labelled as $\langle u,x\rangle$ for clarity. For example, the bottom leftmost $\bullet$ and $\blackdiamond$ are $(1, \langle 1,1\rangle)$ and $(2, \langle 1,1\rangle)$ respectively. The same simplification is used for Figures \ref{figA8.3.14} to \ref{figA8.3.15} and \ref{figA8.3.19} to \ref{figA8.5.24}.}
\end{center}
\end{figure}

\noindent Case 2. $\lambda\ge 4$ and $\mu=2$.
\indent\par Define an orientation $D_{(\lambda,2)}$ for $G^{(2)}$ as follows. (See Figure \ref{figA8.3.14} when $\lambda=4$.)
\begin{align*}
&\langle 1,2\rangle\overset{1}\twoheadrightarrow \langle 2,2\rangle\overset{2}\twoheadleftarrow \langle 3,2\rangle,\  \langle \lambda-2,1\rangle\overset{1}\twoheadrightarrow \langle \lambda-1,1\rangle \overset{2}\twoheadleftarrow \langle \lambda,1\rangle,\\
&\langle i,1\rangle\rightsquigarrow \langle i,2\rangle \text{ for } i=1,2,\ldots, \lambda,\\
&\langle j,1\rangle\rightsquigarrow \langle j+1,1\rangle \text{ for }j=1,2,3,\ldots, \lambda-3, \text{ and }\\
&\langle k,2\rangle \rightsquigarrow \langle k+1,2\rangle \text{ for }k=3,4,\ldots, \lambda-1.
\end{align*}

\begin{center}
\tikzstyle{every node}=[circle, draw, fill=black!100,
                       inner sep=0pt, minimum width=5pt]
\begin{tikzpicture}[thick,scale=0.6]%
\draw(-10,0)node[star,star points=4, label={[yshift=0.5cm] 270:{}}](2_11){};
\draw(-6,0)node[star,star points=4, label={[yshift=0.5cm] 270:{}}](2_21){};
\draw(-2,0)node[star,star points=4, label={[yshift=0.5cm]270:{}}](2_31){};
\draw(2,0)node[star,star points=4, label={[yshift=0.5cm]270:{}}](2_41){};
\draw(-8,2)node[label={[yshift=-0.1cm, xshift=0.1cm]45:{\small $\langle 1,1\rangle$}}](1_11){};
\draw(-4,2)node[label={[yshift=-0.1cm, xshift=0.1cm]45:{\small $\langle 2,1\rangle$}}](1_21){};
\draw(0,2)node[label={[yshift=-0.1cm, xshift=0.1cm]45:{\small $\langle 3,1\rangle$}}](1_31){};
\draw(4,2)node[label={[yshift=0.15cm, xshift=0cm]270:{\small $\langle 4,1\rangle$}}](1_41){};

\draw(-10,4)node[star,star points=4, label={[yshift=-0.35cm, xshift=0.3cm] 90:{}}](2_12){};
\draw(-6,4)node[star,star points=4, label={[yshift=-0.35cm, xshift=0.3cm] 90:{}}](2_22){};
\draw(-2,4)node[star,star points=4, label={[yshift=-0.35cm, xshift=0.3cm] 90:{}}](2_32){};
\draw(2,4)node[star,star points=4, label={[yshift=-0.35cm, xshift=0.3cm] 90:{}}](2_42){};
\draw(-8,6)node[label={[yshift=-0.2cm]90:{\small $\langle 1,2\rangle$}}](1_12){};
\draw(-4,6)node[label={[yshift=-0.2cm]90:{\small $\langle 2,2\rangle$}}](1_22){};
\draw(0,6)node[label={[yshift=-0.2cm]90:{\small $\langle 3,2\rangle$}}](1_32){};
\draw(4,6)node[label={[yshift=-0.2cm]90:{\small $\langle 4,2\rangle$}}](1_42){};

\draw[dashed, ->, line width=0.3mm, >=latex, shorten <= 0.2cm, shorten >= 0.15cm](1_11)--(1_21);
\draw[dashed, ->, line width=0.3mm, >=latex, shorten <= 0.2cm, shorten >= 0.15cm](2_11)--(2_21);
\draw[dashed, ->, line width=0.3mm, >=latex, shorten <= 0.2cm, shorten >= 0.15cm](2_21)--(1_11);
\draw[dashed, ->, line width=0.3mm, >=latex, shorten <= 0.2cm, shorten >= 0.15cm](1_21)--(2_11);

\draw[->, line width=0.3mm, >=latex, shorten <= 0.2cm, shorten >= 0.15cm](1_21)--(1_31);
\draw[->, line width=0.3mm, >=latex, shorten <= 0.2cm, shorten >= 0.15cm](2_31)--(1_21);
\draw[->, line width=0.3mm, >=latex, shorten <= 0.2cm, shorten >= 0.15cm](2_21)--(1_31);
\draw[->, line width=0.3mm, >=latex, shorten <= 0.2cm, shorten >= 0.15cm](2_31)--(2_21);
\draw[->, line width=0.3mm, >=latex, shorten <= 0.2cm, shorten >= 0.15cm](1_31)--(1_41);
\draw[->, line width=0.3mm, >=latex, shorten <= 0.2cm, shorten >= 0.15cm](1_31)--(2_41);
\draw[->, line width=0.3mm, >=latex, shorten <= 0.2cm, shorten >= 0.15cm](1_41)--(2_31);
\draw[->, line width=0.3mm, >=latex, shorten <= 0.2cm, shorten >= 0.15cm](2_41)--(2_31);

\draw[->, line width=0.3mm, >=latex, shorten <= 0.2cm, shorten >= 0.15cm](1_12)--(1_22);
\draw[->, line width=0.3mm, >=latex, shorten <= 0.2cm, shorten >= 0.15cm](2_12)--(1_22);
\draw[->, line width=0.3mm, >=latex, shorten <= 0.2cm, shorten >= 0.15cm](1_22)--(1_32);
\draw[->, line width=0.3mm, >=latex, shorten <= 0.2cm, shorten >= 0.15cm](1_22)--(2_32);
\draw[->, line width=0.3mm, >=latex, shorten <= 0.2cm, shorten >= 0.15cm](2_22)--(1_12);
\draw[->, line width=0.3mm, >=latex, shorten <= 0.2cm, shorten >= 0.15cm](2_22)--(2_12);
\draw[->, line width=0.3mm, >=latex, shorten <= 0.2cm, shorten >= 0.15cm](2_32)--(2_22);
\draw[->, line width=0.3mm, >=latex, shorten <= 0.2cm, shorten >= 0.15cm](1_32)--(2_22);

\draw[dashed, ->, line width=0.3mm, >=latex, shorten <= 0.2cm, shorten >= 0.15cm](1_12)--(2_11);
\draw[dashed, ->, line width=0.3mm, >=latex, shorten <= 0.2cm, shorten >= 0.15cm](2_12)--(1_11);
\draw[dashed, ->, line width=0.3mm, >=latex, shorten <= 0.2cm, shorten >= 0.15cm](1_11)--(1_12);
\draw[dashed, ->, line width=0.3mm, >=latex, shorten <= 0.2cm, shorten >= 0.15cm](2_11)--(2_12);

\draw[dashed, ->, line width=0.3mm, >=latex, shorten <= 0.2cm, shorten >= 0.15cm](1_21)--(1_22);
\draw[dashed, ->, line width=0.3mm, >=latex, shorten <= 0.2cm, shorten >= 0.15cm](2_21)--(2_22);
\draw[dashed, ->, line width=0.3mm, >=latex, shorten <= 0.2cm, shorten >= 0.15cm](1_22)--(2_21);
\draw[dashed, ->, line width=0.3mm, >=latex, shorten <= 0.2cm, shorten >= 0.15cm](2_22)--(1_21);

\draw[dashed, ->, line width=0.3mm, >=latex, shorten <= 0.2cm, shorten >= 0.15cm](1_31)--(1_32);
\draw[dashed, ->, line width=0.3mm, >=latex, shorten <= 0.2cm, shorten >= 0.15cm](2_31)--(2_32);

\draw[dashed, ->, line width=0.3mm, >=latex, shorten <= 0.2cm, shorten >= 0.15cm](1_32)--(2_31);
\draw[dashed, ->, line width=0.3mm, >=latex, shorten <= 0.2cm, shorten >= 0.15cm](2_32)--(1_31);

\draw[dashed, ->, line width=0.3mm, >=latex, shorten <= 0.2cm, shorten >= 0.15cm](1_32)--(1_42);
\draw[dashed, ->, line width=0.3mm, >=latex, shorten <= 0.2cm, shorten >= 0.15cm](1_42)--(2_32);
\draw[dashed, ->, line width=0.3mm, >=latex, shorten <= 0.2cm, shorten >= 0.15cm](2_32)--(2_42);
\draw[dashed, ->, line width=0.3mm, >=latex, shorten <= 0.2cm, shorten >= 0.15cm](2_42)--(1_32);

\draw[dashed, ->, line width=0.3mm, >=latex, shorten <= 0.2cm, shorten >= 0.15cm](1_41)--(1_42);
\draw[dashed, ->, line width=0.3mm, >=latex, shorten <= 0.2cm, shorten >= 0.15cm](1_42)--(2_41);
\draw[dashed, ->, line width=0.3mm, >=latex, shorten <= 0.2cm, shorten >= 0.15cm](2_41)--(2_42);
\draw[dashed, ->, line width=0.3mm, >=latex, shorten <= 0.2cm, shorten >= 0.15cm](2_42)--(1_41);
\end{tikzpicture}
{
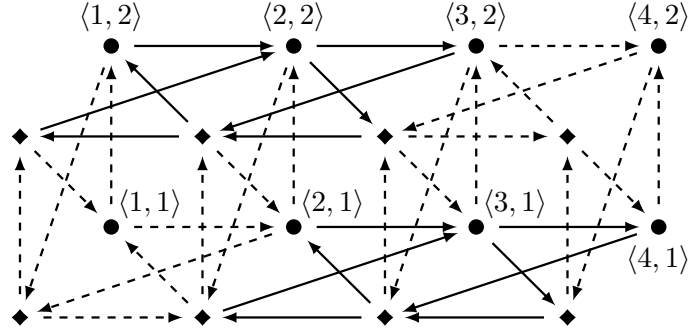
\captionof{figure}{Orientation $D_{(4,2)}$ of $(P_4\times P_2)^{(2)}$, where $d(D_{(4,2)})=4$.}\label{figA8.3.14}}
\end{center}

\indent\par We claim that $d(D_{(\lambda,2)})=d(G)$. Let $u,v\in V(G)$, where $d_G(u,v)\le d(G)-2$. By the definition of $D_{(\lambda,2)}$, we have $d_{D_{(\lambda,2)}}((p,u),(q,v))\le d_G(u,v)+2\le d(G)$ for $p,q=1,2$. Hence, it suffices to consider vertices $u,v\in V(G)$, where $d_G(u,v)= d(G)-1$ or $d_G(u,v)= d(G)$. We illustrate this for $u$ being the `top left' and $v$ being the `bottom right' vertices in Figure \ref{figA8.3.14} and the other cases can be proved similarly. That is, for $(u,v)=(\langle 1,2\rangle,\langle \lambda-1,1\rangle), (\langle 1,2\rangle,\langle \lambda,2\rangle), (\langle 1,2\rangle,\langle \lambda,1\rangle), (\langle 2,2\rangle,\langle \lambda,1\rangle)$, the claim follows by invoking Lemma \ref{lemA8.3.2} on their respective shortest paths:
\begin{align*}
&P^1=\langle 1,2\rangle \langle 2,2\rangle \ldots \langle \lambda-1,2\rangle \langle \lambda-1,1\rangle.\\
&P^2=\langle 1,2\rangle \langle 2,2\rangle \ldots \langle \lambda-1,2\rangle \langle \lambda,2\rangle.\\
&P^3=P^2 \text{ with } \langle \lambda,1\rangle.\\
&P^4=\langle 2,2\rangle \langle 2,1\rangle \langle 3,1\rangle\ldots \langle \lambda,1\rangle.
\end{align*}

\noindent Case 3. $\lambda\ge \mu \ge 3$.
\indent\par Define an orientation $D_{(\lambda,\mu)}$ for $G^{(2)}$ as follows. (See Figure \ref{figA8.3.15} when $\lambda=\mu=3$.)
\begin{align*}
&\Big \langle \Big\lceil\frac{\lambda}{2}\Big\rceil-1,\Big\lceil\frac{\mu}{2}\Big\rceil\Big \rangle\overset{1}\twoheadrightarrow \Big \langle \Big\lceil\frac{\lambda}{2}\Big\rceil,\Big\lceil\frac{\mu}{2}\Big\rceil\Big \rangle \overset{2}\twoheadleftarrow \Big \langle \Big\lceil\frac{\lambda}{2}\Big\rceil+1,\Big\lceil\frac{\mu}{2}\Big\rceil\Big \rangle\text{ and }  \\
&\Big \langle \Big\lceil\frac{\lambda}{2}\Big\rceil,\Big\lceil\frac{\mu}{2}\Big\rceil-1\Big \rangle\overset{1}\twoheadrightarrow \Big \langle \Big\lceil\frac{\lambda}{2}\Big\rceil,\Big\lceil\frac{\mu}{2}\Big\rceil\Big \rangle\overset{2}\twoheadleftarrow \Big \langle \Big\lceil\frac{\lambda}{2}\Big\rceil,\Big\lceil\frac{\mu}{2}\Big\rceil+1\Big \rangle.
\end{align*}
Except for the edges defined above,
\begin{align*}
&\langle i,j\rangle \rightsquigarrow \langle i+1,j\rangle,  \text{ and }\langle i,j\rangle\rightsquigarrow \langle i,j+1\rangle \text{ for all }1\le i\le \lambda-1 \text{ and } 1\le j\le \mu-1.
\end{align*}
\begin{center}
\tikzstyle{every node}=[circle, draw, fill=black!100,
                       inner sep=0pt, minimum width=5pt]
\begin{tikzpicture}[thick,scale=0.6]%
\draw(-10,0)node[star,star points=4, label={[yshift=0.4cm] 270:{}}](2_11){};
\draw(-6,0)node[star,star points=4, label={[yshift=0.4cm] 270:{}}](2_21){};
\draw(-2,0)node[star,star points=4, label={[yshift=0.4cm]270:{}}](2_31){};
\draw(-8,2)node[label={[yshift=-0.1cm, xshift=0.1cm]45:{\small $\langle 1,1\rangle$}}](1_11){};
\draw(-4,2)node[label={[yshift=-0.1cm, xshift=0.1cm]45:{\small $\langle 2,1\rangle$}}](1_21){};
\draw(0,2)node[label={[yshift=0.15cm, xshift=0cm]270:{\small $\langle 3,1\rangle$}}](1_31){};

\draw(-10,4)node[star,star points=4, label={[yshift=0cm, xshift=0cm] 180:{}}](2_12){};
\draw(-6,4)node[star,star points=4, label={[yshift=0.2cm, xshift=-0.18cm] -135:{}}](2_22){};
\draw(-2,4)node[star,star points=4, label={[yshift=0.2cm, xshift=-0.18cm] -135:{}}](2_32){};
\draw(-8,6)node[label={[yshift=-0.1cm, xshift=0.1cm]45:{\small $\langle 1,2\rangle$}}](1_12){};
\draw(-4,6)node[label={[yshift=-0.1cm, xshift=0.1cm]45:{\small $\langle 2,2\rangle$}}](1_22){};
\draw(0,6)node[label={[yshift=-0.1cm, xshift=0.1cm]45:{\small $\langle 3,2\rangle$}}](1_32){};

\draw(-10,8)node[star,star points=4, label={[yshift=-0.3cm, xshift=0.3cm] 90:{}}](2_13){};
\draw(-6,8)node[star,star points=4, label={[yshift=-0.3cm, xshift=0.3cm] 90:{}}](2_23){};
\draw(-2,8)node[star,star points=4, label={[yshift=-0.3cm, xshift=0.3cm] 90:{}}](2_33){};
\draw(-8,10)node[label={[yshift=-0.2cm]90:{\small $\langle 1,3\rangle$}}](1_13){};
\draw(-4,10)node[label={[yshift=-0.2cm]90:{\small $\langle 2,3\rangle$}}](1_23){};
\draw(0,10)node[label={[yshift=-0.2cm]90:{\small $\langle 3,3\rangle$}}](1_33){};

\draw[dashed, ->, line width=0.3mm, >=latex, shorten <= 0.2cm, shorten >= 0.15cm](1_11)--(1_21);
\draw[dashed, ->, line width=0.3mm, >=latex, shorten <= 0.2cm, shorten >= 0.15cm](2_11)--(2_21);
\draw[dashed, ->, line width=0.3mm, >=latex, shorten <= 0.2cm, shorten >= 0.15cm](2_21)--(1_11);
\draw[dashed, ->, line width=0.3mm, >=latex, shorten <= 0.2cm, shorten >= 0.15cm](1_21)--(2_11);

\draw[dashed, ->, line width=0.3mm, >=latex, shorten <= 0.2cm, shorten >= 0.15cm](2_21)--(2_31);
\draw[dashed, ->, line width=0.3mm, >=latex, shorten <= 0.2cm, shorten >= 0.15cm](1_21)--(1_31);
\draw[dashed, ->, line width=0.3mm, >=latex, shorten <= 0.2cm, shorten >= 0.15cm](2_31)--(1_21);
\draw[dashed, ->, line width=0.3mm, >=latex, shorten <= 0.2cm, shorten >= 0.15cm](1_31)--(2_21);

\draw[dashed, ->, line width=0.3mm, >=latex, shorten <= 0.2cm, shorten >= 0.15cm](2_12)--(2_13);
\draw[dashed, ->, line width=0.3mm, >=latex, shorten <= 0.2cm, shorten >= 0.15cm](2_13)--(1_12);
\draw[dashed, ->, line width=0.3mm, >=latex, shorten <= 0.2cm, shorten >= 0.15cm](1_12)--(1_13);
\draw[dashed, ->, line width=0.3mm, >=latex, shorten <= 0.2cm, shorten >= 0.15cm](1_13)--(2_12);

\draw[->, line width=0.3mm, >=latex, shorten <= 0.2cm, shorten >= 0.15cm](1_12)--(1_22);
\draw[->, line width=0.3mm, >=latex, shorten <= 0.2cm, shorten >= 0.15cm](2_12)--(1_22);
\draw[->, line width=0.3mm, >=latex, shorten <= 0.2cm, shorten >= 0.15cm](1_22)--(1_32);
\draw[->, line width=0.3mm, >=latex, shorten <= 0.2cm, shorten >= 0.15cm](1_22)--(2_32);
\draw[->, line width=0.3mm, >=latex, shorten <= 0.2cm, shorten >= 0.15cm](2_22)--(1_12);
\draw[->, line width=0.3mm, >=latex, shorten <= 0.2cm, shorten >= 0.15cm](2_22)--(2_12);
\draw[->, line width=0.3mm, >=latex, shorten <= 0.2cm, shorten >= 0.15cm](2_32)--(2_22);
\draw[->, line width=0.3mm, >=latex, shorten <= 0.2cm, shorten >= 0.15cm](1_32)--(2_22);

\draw[->, line width=0.3mm, >=latex, shorten <= 0.2cm, shorten >= 0.15cm](1_22)--(1_21);
\draw[->, line width=0.3mm, >=latex, shorten <= 0.2cm, shorten >= 0.15cm](1_22)--(2_21);
\draw[->, line width=0.3mm, >=latex, shorten <= 0.2cm, shorten >= 0.15cm](2_21)--(2_22);
\draw[->, line width=0.3mm, >=latex, shorten <= 0.2cm, shorten >= 0.15cm](1_21)--(2_22);

\draw[dashed, ->, line width=0.3mm, >=latex, shorten <= 0.2cm, shorten >= 0.15cm](1_11)--(1_12);
\draw[dashed, ->, line width=0.3mm, >=latex, shorten <= 0.2cm, shorten >= 0.15cm](2_11)--(2_12);
\draw[dashed, ->, line width=0.3mm, >=latex, shorten <= 0.2cm, shorten >= 0.15cm](1_12)--(2_11);
\draw[dashed, ->, line width=0.3mm, >=latex, shorten <= 0.2cm, shorten >= 0.15cm](2_12)--(1_11);

\draw[dashed, ->, line width=0.3mm, >=latex, shorten <= 0.2cm, shorten >= 0.15cm](1_31)--(1_32);
\draw[dashed, ->, line width=0.3mm, >=latex, shorten <= 0.2cm, shorten >= 0.15cm](2_31)--(2_32);
\draw[dashed, ->, line width=0.3mm, >=latex, shorten <= 0.2cm, shorten >= 0.15cm](1_32)--(2_31);
\draw[dashed, ->, line width=0.3mm, >=latex, shorten <= 0.2cm, shorten >= 0.15cm](2_32)--(1_31);

\draw[dashed, ->, line width=0.3mm, >=latex, shorten <= 0.2cm, shorten >= 0.15cm](1_32)--(1_33);
\draw[dashed, ->, line width=0.3mm, >=latex, shorten <= 0.2cm, shorten >= 0.15cm](1_33)--(2_32);
\draw[dashed, ->, line width=0.3mm, >=latex, shorten <= 0.2cm, shorten >= 0.15cm](2_32)--(2_33);
\draw[dashed, ->, line width=0.3mm, >=latex, shorten <= 0.2cm, shorten >= 0.15cm](2_33)--(1_32);

\draw[dashed, ->, line width=0.3mm, >=latex, shorten <= 0.2cm, shorten >= 0.15cm](1_13)--(1_23);
\draw[dashed, ->, line width=0.3mm, >=latex, shorten <= 0.2cm, shorten >= 0.15cm](2_13)--(2_23);
\draw[dashed, ->, line width=0.3mm, >=latex, shorten <= 0.2cm, shorten >= 0.15cm](2_23)--(1_13);
\draw[dashed, ->, line width=0.3mm, >=latex, shorten <= 0.2cm, shorten >= 0.15cm](1_23)--(2_13);

\draw[->, line width=0.3mm, >=latex, shorten <= 0.2cm, shorten >= 0.15cm](1_23)--(1_22);
\draw[->, line width=0.3mm, >=latex, shorten <= 0.2cm, shorten >= 0.15cm](2_23)--(1_22);
\draw[->, line width=0.3mm, >=latex, shorten <= 0.2cm, shorten >= 0.15cm](2_22)--(2_23);
\draw[->, line width=0.3mm, >=latex, shorten <= 0.2cm, shorten >= 0.15cm](2_22)--(1_23);

\draw[dashed, ->, line width=0.3mm, >=latex, shorten <= 0.2cm, shorten >= 0.15cm](1_23)--(1_33);
\draw[dashed, ->, line width=0.3mm, >=latex, shorten <= 0.2cm, shorten >= 0.15cm](1_33)--(2_23);
\draw[dashed, ->, line width=0.3mm, >=latex, shorten <= 0.2cm, shorten >= 0.15cm](2_23)--(2_33);
\draw[dashed, ->, line width=0.3mm, >=latex, shorten <= 0.2cm, shorten >= 0.15cm](2_33)--(1_23);
\end{tikzpicture}
{
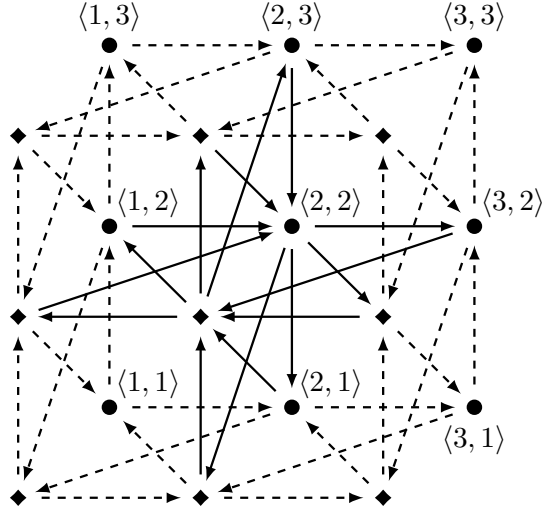
\captionof{figure}{Orientation $D_{(3,3)}$ of $(P_3\times P_3)^{(2)}$, where $d(D_{(3,3)})=4$.}\label{figA8.3.15}}
\end{center}
\indent\par We claim that $d(D_{(\lambda,\mu)})=d(G)$. Similar to Case 2, it suffices to consider $u,v\in V(G)$, where $d_G(u,v)= d(G)-1$ or $d_G(u,v)= d(G)$. We illustrate this for $u$ being the `top left' and $v$ being the `bottom right' vertices in Figure \ref{figA8.3.15} and the other cases can be proved similarly. That is, for $(u,v)=(\langle 1,\mu\rangle,\langle \lambda-1,1\rangle), (\langle 1,\mu\rangle,\langle \lambda,2\rangle), (\langle 1,\mu\rangle,\langle \lambda,1\rangle), (\langle 2,\mu\rangle,\langle \lambda,1\rangle),$ $(\langle 1,\mu-1\rangle,\langle \lambda,1\rangle)$, the claim follows by invoking Lemma \ref{lemA8.3.2} on their respective shortest paths:

\begin{align*}
P^1=&\Big \langle 1,\mu\Big \rangle \Big \langle 2,\mu\Big \rangle\ldots\Big \langle \Big\lceil\frac{\lambda}{2}\Big\rceil,\mu\Big \rangle \Big \langle \Big\lceil\frac{\lambda}{2}\Big\rceil,\mu-1\Big \rangle\ldots \Big \langle \Big\lceil\frac{\lambda}{2}\Big\rceil,1\Big \rangle \Big \langle \Big\lceil\frac{\lambda}{2}\Big\rceil+1,1\Big \rangle\ldots\Big \langle \lambda-1,1\Big \rangle.\\
P^2=&\Big \langle 1,\mu\Big \rangle \Big \langle 1,\mu-1\Big \rangle\ldots \Big \langle 1,\Big\lceil\frac{\mu}{2}\Big\rceil\Big \rangle \Big \langle 2,\Big\lceil\frac{\mu}{2}\Big\rceil\Big \rangle\ldots\Big \langle \lambda,\Big\lceil\frac{\mu}{2}\Big\rceil\Big \rangle \Big \langle \lambda,\Big\lceil\frac{\mu}{2}\Big\rceil-1\Big \rangle\ldots\Big \langle \lambda,2\Big \rangle.\\
P^3=&P^2 \text{ with }\Big \langle \lambda,1\Big \rangle.\\
P^4=&\Big \langle 2,\mu\Big \rangle \Big \langle 3,\mu\Big \rangle\ldots\Big \langle \Big\lceil\frac{\lambda}{2}\Big\rceil,\mu\Big \rangle \Big \langle \Big\lceil\frac{\lambda}{2}\Big\rceil,\mu-1\Big \rangle\ldots \Big \langle \Big\lceil\frac{\lambda}{2}\Big\rceil,1\Big \rangle \Big \langle \Big\lceil\frac{\lambda}{2}\Big\rceil+1,1\Big \rangle\ldots\Big \langle \lambda,1\Big \rangle.\\
P^5=&\Big \langle 1,\mu-1\Big \rangle \Big \langle 1,\mu-2\Big\rangle\ldots \Big \langle 1,\Big\lceil\frac{\mu}{2}\Big\rceil\Big \rangle \Big \langle 2,\Big\lceil\frac{\mu}{2}\Big\rceil\Big \rangle\ldots\Big \langle \lambda,\Big\lceil\frac{\mu}{2}\Big\rceil\Big \rangle \Big \langle \lambda,\Big\lceil\frac{\mu}{2}\Big\rceil-1\Big \rangle\ldots\Big \langle \lambda,1\Big \rangle.
\end{align*}

Hence, $G^{(2)}\in\mathscr{C}_0$ for Cases 2 and 3. To complete (b), observe that every vertex lies in a directed $C_4$ in each orientation $D_{(\lambda,\mu)}$ of all three cases and invoke Lemma \ref{lemA8.1.2}.
\qed

\indent\par If $G^{(2)}$ belongs to $\mathscr{C}_1$, a possible direction for further research is if we can increase the vertex-multiplication, $s_i$, of a particular vertex so that the resulting graph $G(s_1,s_2,\ldots,s_n)$ belongs to $\mathscr{C}_0$. This was shown possible in \cite{WHW TEG 6A} when the parent graph $G$ is a tree of diameter $4$. Propositions \ref{ppnA8.3.3} and \ref{ppnA8.3.4} show $P_3\times P_2$ as a contrasting example, i.e. however large the increase in vertex-multiplication of any particular vertex in $(P_3\times P_2)^{(2)}$, the resulting graph still lies in $\mathscr{C}_1$.
\begin{ppn}\label{ppnA8.3.3}
If $j\in\{\langle 2,1\rangle, \langle 2,2\rangle\}$ and $s_i=2$ for all $i\neq j$, then $(P_3\times P_2)(s_1,s_2,\ldots,s_n)$ $\in\mathscr{C}_1$.
\end{ppn}
\noindent\textit{Proof}: By Theorem \ref{thmA8.1.8}(b), it suffices to prove $\bar{d}((P_3\times P_2)(s_1,s_2,\ldots,s_n))>3$. Suppose there exists an orientation $D$ of $(P_3\times P_2)(s_1,s_2,\ldots, s_n)$ such that $d(D)=3$. By Theorem \ref{thmA8.1.8} and symmetry, it suffices to consider $j=\langle 2,1\rangle$, i.e. $s_{\langle 2,1\rangle}\ge 3$ (see Figure \ref{figA8.3.16} when $s_{\langle 2,1\rangle}=4$). Since $d_D((p,\langle 1,2\rangle),(q,\langle 3,2\rangle))=d_D((q,\langle 3,2\rangle),(p,\langle 1,2\rangle))=2$ for all $p,q=1,2$, it follows from Lemma \ref{lemA8.3.1} that $\langle 1,2\rangle\overset{1}\twoheadrightarrow \langle 2,2\rangle \overset{2}\twoheadleftarrow \langle 3,2\rangle$.
\indent\par For $i=1,3$, since $d_D((p,\langle i,2\rangle),(3-p,\langle i,2\rangle))\le 3$ for $p=1,2$, it follows WLOG that $\langle i,1\rangle\rightsquigarrow \langle i,2\rangle$. By $d_D((1,\langle 1,1\rangle),(1,\langle 3,1\rangle))\le 3$, we may assume WLOG that $(1,\langle 1,1\rangle)\rightarrow (1,\langle 2,1\rangle)\rightarrow (1,\langle 3,1\rangle)$. Now, it is necessary from $d_D((1,\langle 2,1\rangle),(1,\langle 1,2\rangle))\le 3$ that $(1,\langle 2,1\rangle)$ $\rightarrow (2,\langle 2,2\rangle)$. However, $d_D((2,\langle 3,2\rangle),(1,\langle 2,1\rangle))>3$, a contradiction.
\qed
\begin{center}
\tikzstyle{every node}=[circle, draw, fill=black!100,
                       inner sep=0pt, minimum width=5pt]
\begin{tikzpicture}[thick,scale=0.7]%
\draw(-10,0)node[label={[yshift=0.4cm] 270:{\small $(2,\langle 1,1\rangle)$}}](2_11){};
\draw(-4.67,1.33)node[label={[xshift=0.2cm] 0:{\small $(2,\langle 2,1\rangle)$}}](2_21){};
\draw(-5.33,0.67)node[label={[xshift=0.2cm] 0:{\small $(3,\langle 2,1\rangle)$}}](3_21){};
\draw(-6,0)node[label={[yshift=0.4cm] 270:{\small $(4,\langle 2,1\rangle)$}}](4_21){};

\draw(-2,0)node[label={[yshift=0.4cm]270:{\small $(2,\langle 3,1\rangle)$}}](2_31){};
\draw(-8,2)node[label={[yshift=-0.6cm, xshift=0.9cm]:{\small $(1,\langle 1,1\rangle)$}}](1_11){};
\draw(-4,2)node[label={[yshift=-0.6cm, xshift=0.6cm]:{\small $(1,\langle 2,1\rangle)$}}](1_21){};
\draw(0,2)node[label={[yshift=0.5cm, xshift=0cm]270:{\small $(1,\langle 3,1\rangle)$}}](1_31){};

\draw(-10,4)node[label={[yshift=-0.35cm, xshift=0.3cm] 90:{\small $(2,\langle 1,2\rangle)$}}](2_12){};
\draw(-6,4)node[label={[yshift=-0.35cm, xshift=0.3cm] 90:{\small $(2,\langle 2,2\rangle)$}}](2_22){};
\draw(-2,4)node[label={[yshift=-0.35cm, xshift=0.3cm] 90:{\small $(2,\langle 3,2\rangle)$}}](2_32){};
\draw(-8,6)node[label={[yshift=-0.5cm]90:{\small $(1,\langle 1,2\rangle)$}}](1_12){};
\draw(-4,6)node[label={[yshift=-0.5cm]90:{\small $(1,\langle 2,2\rangle)$}}](1_22){};
\draw(0,6)node[label={[yshift=-0.5cm]90:{\small $(1,\langle 3,2\rangle)$}}](1_32){};

\draw[->, line width=0.3mm, >=latex, shorten <= 0.2cm, shorten >= 0.15cm](1_11)--(1_21);
\draw[->, line width=0.3mm, >=latex, shorten <= 0.2cm, shorten >= 0.15cm](1_21)--(1_31);

\draw[->, line width=0.3mm, >=latex, shorten <= 0.2cm, shorten >= 0.15cm](1_12)--(1_22);
\draw[->, line width=0.3mm, >=latex, shorten <= 0.2cm, shorten >= 0.15cm](2_12)--(1_22);
\draw[->, line width=0.3mm, >=latex, shorten <= 0.2cm, shorten >= 0.15cm](1_22)--(1_32);
\draw[->, line width=0.3mm, >=latex, shorten <= 0.2cm, shorten >= 0.15cm](1_22)--(2_32);
\draw[->, line width=0.3mm, >=latex, shorten <= 0.2cm, shorten >= 0.15cm](2_22)--(1_12);
\draw[->, line width=0.3mm, >=latex, shorten <= 0.2cm, shorten >= 0.15cm](2_22)--(2_12);
\draw[->, line width=0.3mm, >=latex, shorten <= 0.2cm, shorten >= 0.15cm](2_32)--(2_22);
\draw[->, line width=0.3mm, >=latex, shorten <= 0.2cm, shorten >= 0.15cm](1_32)--(2_22);

\draw[dashed, ->, line width=0.3mm, >=latex, shorten <= 0.2cm, shorten >= 0.15cm](1_11)--(1_12);
\draw[dashed, ->, line width=0.3mm, >=latex, shorten <= 0.2cm, shorten >= 0.15cm](2_11)--(2_12);

\draw[dashed, ->, line width=0.3mm, >=latex, shorten <= 0.2cm, shorten >= 0.15cm](1_31)--(1_32);
\draw[dashed, ->, line width=0.3mm, >=latex, shorten <= 0.2cm, shorten >= 0.15cm](2_31)--(2_32);

\draw[dashed, ->, line width=0.3mm, >=latex, shorten <= 0.2cm, shorten >= 0.15cm](1_12)--(2_11);
\draw[dashed, ->, line width=0.3mm, >=latex, shorten <= 0.2cm, shorten >= 0.15cm](2_12)--(1_11);
\draw[dashed, ->, line width=0.3mm, >=latex, shorten <= 0.2cm, shorten >= 0.15cm](1_32)--(2_31);
\draw[dashed, ->, line width=0.3mm, >=latex, shorten <= 0.2cm, shorten >= 0.15cm](2_32)--(1_31);
\draw[->, line width=0.3mm, >=latex, shorten <= 0.2cm, shorten >= 0.15cm](1_21)--(2_22);
\end{tikzpicture}
{\captionof{figure}{Partial orientation $D$, where $d(D)=4$.}\label{figA8.3.16}}
\end{center}

\begin{ppn}\label{ppnA8.3.4}
If $s_i=2$ for all $i=\langle 2,1\rangle, \langle 2,2\rangle$, then $(P_3\times P_2)(s_1,s_2,\ldots,s_n)\in\mathscr{C}_1$.
\end{ppn}
\noindent\textit{Proof}: By Theorem \ref{thmA8.1.8}(b), it suffices to prove $\bar{d}((P_3\times P_2)(s_1,s_2,\ldots,s_n))>3$. Suppose there exists an orientation $D$ of $(P_3\times P_2)(s_1,s_2,\ldots,s_n)$ such that $d(D)=3$. By Lemma \ref{lemA8.3.1}, $\langle 1,1\rangle\overset{1}\twoheadrightarrow \langle 2,1\rangle \overset{2}\twoheadleftarrow \langle 3,1\rangle$ (see Figure \ref{figA8.3.17} when $s_{\langle 1,1\rangle}=s_{\langle 1,2\rangle}=s_{\langle 3,1\rangle}=s_{\langle 3,2\rangle}=4$). For $i=3,4,\ldots, s_{\langle 1,1\rangle}$ and $j=3,4,\ldots, s_{\langle 3,1\rangle}$, replace $(2,\langle 1,1\rangle)$ ($(2,\langle 3,1\rangle)$ resp.) by $(i,\langle 1,1\rangle)$ ($(j,\langle 3,1\rangle)$ resp.) and apply Lemma \ref{lemA8.3.1} to conclude $(p,\langle 1,1\rangle)\rightarrow (1,\langle 2,1\rangle)\rightarrow (q,\langle 3,1\rangle)\rightarrow (2,\langle 2,1\rangle)\rightarrow (p,\langle 1,1\rangle)$ for all $p=1,2,\ldots, s_{\langle 1,1\rangle}$ and $q=1,2,\ldots, s_{\langle 3,1\rangle}$. Similarly, $(p,\langle 1,2\rangle)\rightarrow (1,\langle 2,2\rangle)\rightarrow (q,\langle 3,2\rangle)\rightarrow (2,\langle 2,2\rangle)\rightarrow (p,\langle 1,2\rangle)$ for all $p=1,2,\ldots, s_{\langle 1,2\rangle}$ and $q=1,2,\ldots, s_{\langle 3,2\rangle}$. 
\indent\par By $d_D((1,\langle 2,1\rangle),(1,\langle 1,2\rangle))\le 3$, we must have $(1,\langle 2,1\rangle)\rightarrow (2,\langle 2,2\rangle)$. However, $d_D((1,\langle 3,2\rangle),(1,\langle 2,1\rangle))> 3$, a contradiction.
\qed
\begin{center}
\tikzstyle{every node}=[circle, draw, fill=black!100,
                       inner sep=0pt, minimum width=5pt]
\begin{tikzpicture}[thick,scale=0.7]%
\draw(-8,2)node[label={[xshift=-0.2cm]180:{\small $(1,\langle 1,1\rangle)$}}](1_11){};
\draw(-8.67,1.33)node[label={[xshift=-0.2cm] 180:{\small $(2,\langle 1,1\rangle)$}}](2_11){};
\draw(-9.33,0.67)node[label={[xshift=-0.2cm] 180:{\small $(3,\langle 1,1\rangle)$}}](3_11){};
\draw(-10,0)node[label={[xshift=-0.2cm] 180:{\small $(4,\langle 1,1\rangle)$}}](4_11){};

\draw(-8,6)node[label={[xshift=-0.2cm]180:{\small $(1,\langle 1,2\rangle)$}}](1_12){};
\draw(-8.67,5.33)node[label={[xshift=-0.2cm] 180:{\small $(2,\langle 1,2\rangle)$}}](2_12){};
\draw(-9.33,4.67)node[label={[xshift=-0.2cm] 180:{\small $(3,\langle 1,2\rangle)$}}](3_12){};
\draw(-10,4)node[label={[xshift=-0.2cm] 180:{\small $(4,\langle 1,2\rangle)$}}](4_12){};

\draw(-4,2)node[label={[yshift=-0.6cm, xshift=0.55cm]:{\small $(1,\langle 2,1\rangle)$}}](1_21){};
\draw(-6,0)node[label={[yshift=0.4cm] 270:{\small $(2,\langle 2,1\rangle)$}}](2_21){};

\draw(-4,6)node[label={[yshift=-0.6cm]90:{\small $(1,\langle 2,2\rangle)$}}](1_22){};
\draw(-6,4)node[label={[yshift=-0.3cm, xshift=0.3cm] 90:{\small $(2,\langle 2,2\rangle)$}}](2_22){};

\draw(0,2)node[label={[xshift=0.2cm]0:{\small $(1,\langle 3,1\rangle)$}}](1_31){};
\draw(-0.67,1.33)node[label={[xshift=0.2cm] 0:{\small $(2,\langle 3,1\rangle)$}}](2_31){};
\draw(-1.33,0.67)node[label={[xshift=0.2cm] 0:{\small $(3,\langle 3,1\rangle)$}}](3_31){};
\draw(-2,0)node[label={[xshift=0.2cm] 0:{\small $(4,\langle 3,1\rangle)$}}](4_31){};

\draw(0,6)node[label={[xshift=0.2cm]0:{\small $(1,\langle 3,2\rangle)$}}](1_32){};
\draw(-0.67,5.33)node[label={[xshift=0.2cm] 0:{\small $(2,\langle 3,2\rangle)$}}](2_32){};
\draw(-1.33,4.67)node[label={[xshift=0.2cm] 0:{\small $(3,\langle 3,2\rangle)$}}](3_32){};
\draw(-2,4)node[label={[xshift=0.2cm] 0:{\small $(4,\langle 3,2\rangle)$}}](4_32){};

\draw[->, line width=0.3mm, >=latex, shorten <= 0.2cm, shorten >= 0.15cm](1_11)--(1_21);
\draw[->, line width=0.3mm, >=latex, shorten <= 0.2cm, shorten >= 0.15cm](2_11)--(1_21);
\draw[->, line width=0.3mm, >=latex, shorten <= 0.2cm, shorten >= 0.15cm](3_11)--(1_21);
\draw[->, line width=0.3mm, >=latex, shorten <= 0.2cm, shorten >= 0.15cm](4_11)--(1_21);

\draw[->, line width=0.3mm, >=latex, shorten <= 0.2cm, shorten >= 0.15cm](1_31)--(2_21);
\draw[->, line width=0.3mm, >=latex, shorten <= 0.2cm, shorten >= 0.15cm](2_31)--(2_21);
\draw[->, line width=0.3mm, >=latex, shorten <= 0.2cm, shorten >= 0.15cm](3_31)--(2_21);
\draw[->, line width=0.3mm, >=latex, shorten <= 0.2cm, shorten >= 0.15cm](4_31)--(2_21);

\draw[->, line width=0.3mm, >=latex, shorten <= 0.2cm, shorten >= 0.15cm](2_21)--(1_11);
\draw[->, line width=0.3mm, >=latex, shorten <= 0.2cm, shorten >= 0.15cm](2_21)--(2_11);
\draw[->, line width=0.3mm, >=latex, shorten <= 0.2cm, shorten >= 0.15cm](2_21)--(3_11);
\draw[->, line width=0.3mm, >=latex, shorten <= 0.2cm, shorten >= 0.15cm](2_21)--(4_11);

\draw[->, line width=0.3mm, >=latex, shorten <= 0.2cm, shorten >= 0.15cm](1_21)--(1_31);
\draw[->, line width=0.3mm, >=latex, shorten <= 0.2cm, shorten >= 0.15cm](1_21)--(2_31);
\draw[->, line width=0.3mm, >=latex, shorten <= 0.2cm, shorten >= 0.15cm](1_21)--(3_31);
\draw[->, line width=0.3mm, >=latex, shorten <= 0.2cm, shorten >= 0.15cm](1_21)--(4_31);

\draw[->, line width=0.3mm, >=latex, shorten <= 0.2cm, shorten >= 0.15cm](1_12)--(1_22);
\draw[->, line width=0.3mm, >=latex, shorten <= 0.2cm, shorten >= 0.15cm](2_12)--(1_22);
\draw[->, line width=0.3mm, >=latex, shorten <= 0.2cm, shorten >= 0.15cm](3_12)--(1_22);
\draw[->, line width=0.3mm, >=latex, shorten <= 0.2cm, shorten >= 0.15cm](4_12)--(1_22);

\draw[->, line width=0.3mm, >=latex, shorten <= 0.2cm, shorten >= 0.15cm](1_22)--(1_32);
\draw[->, line width=0.3mm, >=latex, shorten <= 0.2cm, shorten >= 0.15cm](1_22)--(2_32);
\draw[->, line width=0.3mm, >=latex, shorten <= 0.2cm, shorten >= 0.15cm](1_22)--(3_32);
\draw[->, line width=0.3mm, >=latex, shorten <= 0.2cm, shorten >= 0.15cm](1_22)--(4_32);

\draw[->, line width=0.3mm, >=latex, shorten <= 0.2cm, shorten >= 0.15cm](2_22)--(1_12);
\draw[->, line width=0.3mm, >=latex, shorten <= 0.2cm, shorten >= 0.15cm](2_22)--(2_12);
\draw[->, line width=0.3mm, >=latex, shorten <= 0.2cm, shorten >= 0.15cm](2_22)--(3_12);
\draw[->, line width=0.3mm, >=latex, shorten <= 0.2cm, shorten >= 0.15cm](2_22)--(4_12);

\draw[->, line width=0.3mm, >=latex, shorten <= 0.2cm, shorten >= 0.15cm](1_32)--(2_22);
\draw[->, line width=0.3mm, >=latex, shorten <= 0.2cm, shorten >= 0.15cm](2_32)--(2_22);
\draw[->, line width=0.3mm, >=latex, shorten <= 0.2cm, shorten >= 0.15cm](3_32)--(2_22);
\draw[->, line width=0.3mm, >=latex, shorten <= 0.2cm, shorten >= 0.15cm](4_32)--(2_22);

\draw[->, line width=0.3mm, >=latex, shorten <= 0.2cm, shorten >= 0.15cm](1_21)--(2_22);
\end{tikzpicture}
{\captionof{figure}{Partial orientation $D$, where $d(D)=4$.}\label{figA8.3.17}}
\end{center}

\begin{eg}\label{egA8.3.5}
If $s_i=2$ for all $i\not\in \{\langle 2,1\rangle,\langle 2,2\rangle\}$ and $s_j=4$ for all $j\in\{\langle 2,1\rangle,\langle 2,2\rangle\}$, then $(P_3\times P_2)(4,4,\overbrace{2,2,\ldots,2}^8)\in\mathscr{C}_0$.
\end{eg}
\noindent\textit{Proof}: Define an orientation $D$ of $(P_3\times P_2)(4,4,2,2,\ldots,2)$. (See Figure \ref{figA8.3.18}.)
\begin{align*}
&\langle 1,2\rangle\rightrightarrows \langle 1,1\rangle, \langle 3,2\rangle\rightrightarrows \langle 3,1\rangle,\\
&(p,\langle 2,1\rangle)\rightarrow (q,\langle 2,2\rangle) \text{ for }p,q=1,2,3,4.\\
&\{(2,\langle 2,i\rangle),(4,\langle 2,i\rangle)\}\rightarrow (1,\langle 1,i\rangle) \rightarrow \{(1,\langle 2,i\rangle),(3,\langle 2,i\rangle)\}\rightarrow (2,\langle 1,i\rangle) \rightarrow \{(2,\langle 2,i\rangle),(4,\langle 2,i\rangle)\},\\
&\{(1,\langle 2,i\rangle),(4,\langle 2,i\rangle)\}\rightarrow (1,\langle 3,i\rangle) \rightarrow \{(2,\langle 2,i\rangle),(3,\langle 2,i\rangle)\}\rightarrow (2,\langle 3,i\rangle) \rightarrow \{(1,\langle 2,i\rangle),(4,\langle 2,i\rangle)\}
\end{align*}
for $i=1,2$. It is easy to verify that $d(D)=3$.
\qed
\begin{center}
\tikzstyle{every node}=[circle, draw, fill=black!100,
                       inner sep=0pt, minimum width=5pt]
\begin{tikzpicture}[thick,scale=0.7]%
\draw(-8,2)node[label={[yshift=-0.6cm, xshift=0.65cm]90:{\small $(1,\langle 1,1\rangle)$}}](1_11){};
\draw(-10,0)node[label={[yshift=0.5cm] 270:{\small $(2,\langle 1,1\rangle)$}}](2_11){};

\draw(-4,2)node[label={[yshift=-0.6cm, xshift=0.8cm]:{\small $(1,\langle 2,1\rangle)$}}](1_21){};
\draw(-4.67,1.33)node[label={[xshift=0.2cm] 0:{$$}}](2_21){};
\draw(-5.33,0.67)node[label={[xshift=0.2cm] 0:{$$}}](3_21){};
\draw(-6,0)node[label={[yshift=0.5cm] 270:{\small $(4,\langle 2,1\rangle)$}}](4_21){};

\draw(-8,6)node[label={[yshift=-0.5cm]90:{\small $(1,\langle 1,2\rangle)$}}](1_12){};
\draw(-10,4)node[label={[yshift=-0.35cm, xshift=0.3cm] 90:{\small $(2,\langle 1,2\rangle)$}}](2_12){};

\draw(-4,6)node[label={[yshift=-0.5cm]90:{\small $(1,\langle 2,2\rangle)$}}](1_22){};
\draw(-4.67,5.33)node[label={[xshift=0.2cm] 0:{$$}}](2_22){};
\draw(-5.33,4.67)node[label={[xshift=0.2cm] 0:{$$}}](3_22){};
\draw(-6,4)node[label={[yshift=-0.35cm, xshift=0.1cm] 180:{\small $(4,\langle 2,2\rangle)$}}](4_22){};

\draw(0,2)node[label={[yshift=0.4cm, xshift=0cm]270:{\small $(1,\langle 3,1\rangle)$}}](1_31){};
\draw(-2,0)node[label={[yshift=0.5cm]270:{\small $(2,\langle 3,1\rangle)$}}](2_31){};
\draw(0,6)node[label={[yshift=-0.5cm]90:{\small $(1,\langle 3,2\rangle)$}}](1_32){};
\draw(-2,4)node[label={[yshift=-0.35cm, xshift=0.3cm] 90:{\small $(2,\langle 3,2\rangle)$}}](2_32){};

\draw[->, line width=0.3mm, >=latex, shorten <= 0.2cm, shorten >= 0.15cm](1_12)--(1_22);
\draw[->, line width=0.3mm, >=latex, shorten <= 0.2cm, shorten >= 0.15cm](1_12)--(3_22);
\draw[dashed, ->, line width=0.3mm, >=latex, shorten <= 0.2cm, shorten >= 0.15cm](2_22)--(1_12);
\draw[dashed, ->, line width=0.3mm, >=latex, shorten <= 0.2cm, shorten >= 0.15cm](4_22)--(1_12);

\draw[->, line width=0.3mm, >=latex, shorten <= 0.2cm, shorten >= 0.15cm](2_12)--(2_22);
\draw[->, line width=0.3mm, >=latex, shorten <= 0.2cm, shorten >= 0.15cm](2_12)--(4_22);
\draw[dashed, ->, line width=0.3mm, >=latex, shorten <= 0.2cm, shorten >= 0.15cm](1_22)--(2_12);
\draw[dashed, ->, line width=0.3mm, >=latex, shorten <= 0.2cm, shorten >= 0.15cm](3_22)--(2_12);

\draw[->, line width=0.3mm, >=latex, shorten <= 0.2cm, shorten >= 0.15cm](1_32)--(2_22);
\draw[->, line width=0.3mm, >=latex, shorten <= 0.2cm, shorten >= 0.15cm](1_32)--(3_22);
\draw[dashed, ->, line width=0.3mm, >=latex, shorten <= 0.2cm, shorten >= 0.15cm](1_22)--(1_32);
\draw[dashed, ->, line width=0.3mm, >=latex, shorten <= 0.2cm, shorten >= 0.15cm](4_22)--(1_32);

\draw[->, line width=0.3mm, >=latex, shorten <= 0.2cm, shorten >= 0.15cm](2_32)--(1_22);
\draw[->, line width=0.3mm, >=latex, shorten <= 0.2cm, shorten >= 0.15cm](2_32)--(4_22);
\draw[dashed, ->, line width=0.3mm, >=latex, shorten <= 0.2cm, shorten >= 0.15cm](2_22)--(2_32);
\draw[dashed, ->, line width=0.3mm, >=latex, shorten <= 0.2cm, shorten >= 0.15cm](3_22)--(2_32);

\draw[->, line width=0.3mm, >=latex, shorten <= 0.2cm, shorten >= 0.15cm](1_11)--(1_21);
\draw[->, line width=0.3mm, >=latex, shorten <= 0.2cm, shorten >= 0.15cm](1_11)--(3_21);
\draw[dashed, ->, line width=0.3mm, >=latex, shorten <= 0.2cm, shorten >= 0.15cm](2_21)--(1_11);
\draw[dashed, ->, line width=0.3mm, >=latex, shorten <= 0.2cm, shorten >= 0.15cm](4_21)--(1_11);

\draw[->, line width=0.3mm, >=latex, shorten <= 0.2cm, shorten >= 0.15cm](2_11)--(2_21);
\draw[->, line width=0.3mm, >=latex, shorten <= 0.2cm, shorten >= 0.15cm](2_11)--(4_21);
\draw[dashed, ->, line width=0.3mm, >=latex, shorten <= 0.2cm, shorten >= 0.15cm](1_21)--(2_11);
\draw[dashed, ->, line width=0.3mm, >=latex, shorten <= 0.2cm, shorten >= 0.15cm](3_21)--(2_11);

\draw[->, line width=0.3mm, >=latex, shorten <= 0.2cm, shorten >= 0.15cm](1_31)--(2_21);
\draw[->, line width=0.3mm, >=latex, shorten <= 0.2cm, shorten >= 0.15cm](1_31)--(3_21);
\draw[dashed, ->, line width=0.3mm, >=latex, shorten <= 0.2cm, shorten >= 0.15cm](1_21)--(1_31);
\draw[dashed, ->, line width=0.3mm, >=latex, shorten <= 0.2cm, shorten >= 0.15cm](4_21)--(1_31);

\draw[->, line width=0.3mm, >=latex, shorten <= 0.2cm, shorten >= 0.15cm](2_31)--(1_21);
\draw[->, line width=0.3mm, >=latex, shorten <= 0.2cm, shorten >= 0.15cm](2_31)--(4_21);
\draw[dashed, ->, line width=0.3mm, >=latex, shorten <= 0.2cm, shorten >= 0.15cm](2_21)--(2_31);
\draw[dashed, ->, line width=0.3mm, >=latex, shorten <= 0.2cm, shorten >= 0.15cm](3_21)--(2_31);
\draw[densely dotted, ->, line width=0.3mm, >=latex, shorten <= 0.2cm, shorten >= 0.15cm](1_12) to [out=260, in=110](1_11);
\draw[densely dotted, ->, line width=0.3mm, >=latex, shorten <= 0.2cm, shorten >= 0.15cm](1_12)--(2_11);
\draw[densely dotted, ->, line width=0.3mm, >=latex, shorten <= 0.2cm, shorten >= 0.15cm](2_12)--(1_11);
\draw[densely dotted, ->, line width=0.3mm, >=latex, shorten <= 0.2cm, shorten >= 0.15cm](2_12)--(2_11);
\draw[densely dotted, ->, line width=0.3mm, >=latex, shorten <= 0.2cm, shorten >= 0.15cm](1_32)--(1_31);
\draw[densely dotted, ->, line width=0.3mm, >=latex, shorten <= 0.2cm, shorten >= 0.15cm](1_32)--(2_31);
\draw[densely dotted, ->, line width=0.3mm, >=latex, shorten <= 0.2cm, shorten >= 0.15cm](2_32)--(1_31);
\draw[densely dotted, ->, line width=0.3mm, >=latex, shorten <= 0.2cm, shorten >= 0.15cm](2_32) to [out=290, in=80](2_31);
\draw[densely dotted, ->, line width=0.3mm, >=latex, shorten <= 0.2cm, shorten >= 0.15cm](1_21)--(1_22);
\draw[densely dotted, ->, line width=0.3mm, >=latex, shorten <= 0.2cm, shorten >= 0.15cm](1_21)--(2_22);
\draw[densely dotted, ->, line width=0.3mm, >=latex, shorten <= 0.2cm, shorten >= 0.15cm](1_21)--(3_22);
\draw[densely dotted, ->, line width=0.3mm, >=latex, shorten <= 0.2cm, shorten >= 0.15cm](1_21)--(4_22);

\draw[densely dotted, ->, line width=0.3mm, >=latex, shorten <= 0.2cm, shorten >= 0.15cm](2_21)--(1_22);
\draw[densely dotted, ->, line width=0.3mm, >=latex, shorten <= 0.2cm, shorten >= 0.15cm](2_21)--(2_22);
\draw[densely dotted, ->, line width=0.3mm, >=latex, shorten <= 0.2cm, shorten >= 0.15cm](2_21)--(3_22);
\draw[densely dotted, ->, line width=0.3mm, >=latex, shorten <= 0.2cm, shorten >= 0.15cm](2_21)--(4_22);

\draw[densely dotted, ->, line width=0.3mm, >=latex, shorten <= 0.2cm, shorten >= 0.15cm](3_21)--(1_22);
\draw[densely dotted, ->, line width=0.3mm, >=latex, shorten <= 0.2cm, shorten >= 0.15cm](3_21)--(2_22);
\draw[densely dotted, ->, line width=0.3mm, >=latex, shorten <= 0.2cm, shorten >= 0.15cm](3_21)--(3_22);
\draw[densely dotted, ->, line width=0.3mm, >=latex, shorten <= 0.2cm, shorten >= 0.15cm](3_21)--(4_22);

\draw[densely dotted, ->, line width=0.3mm, >=latex, shorten <= 0.2cm, shorten >= 0.15cm](4_21)--(1_22);
\draw[densely dotted, ->, line width=0.3mm, >=latex, shorten <= 0.2cm, shorten >= 0.15cm](4_21)--(2_22);
\draw[densely dotted, ->, line width=0.3mm, >=latex, shorten <= 0.2cm, shorten >= 0.15cm](4_21)--(3_22);
\draw[densely dotted, ->, line width=0.3mm, >=latex, shorten <= 0.2cm, shorten >= 0.15cm](4_21)--(4_22);
\end{tikzpicture}
{
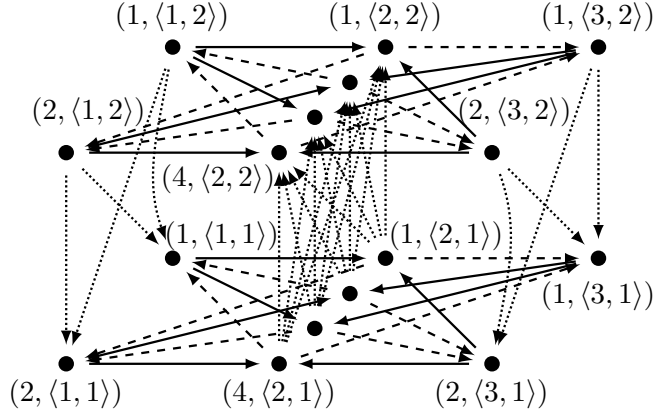
\captionof{figure}{Orientation $D$, where $d(D)=3$.}\label{figA8.3.18}}
\end{center}

\indent\par It may seem favourable to generalise Example \ref{egA8.3.5} by replacing the condition $``s_j=4"$ with $``s_j\ge 4"$ for all $j\in\{\langle 2,1\rangle,\langle 2,2\rangle\}$. Proposition \ref{ppnA8.3.7} shows that this does not hold. To this end, recall the classical Sperner Theorem.
\begin{thm}(Sperner \cite{SE})\label{Sperner} Let $n\in\mathbb{Z}^+$ and $\mathscr{A}$ be an antichain of $\mathbb{N}_n=\{1,2,\ldots,n\}$ (i.e. $A\not\subseteq B$ for all $A,B\in\mathscr{A}$). Then, $|\mathscr{A}|\le {{n}\choose{\lfloor{n/2}\rfloor}}$ with equality holding if and only if all members in $\mathscr{A}$ have the same size, ${\lfloor\frac{n}{2}\rfloor}$ or ${\lceil\frac{n}{2}\rceil}$. (The two sizes coincide if $n$ is even.)
\end{thm}
\begin{ppn}\label{ppnA8.3.7}
If $s_i=2$ for all $i\not\in \{\langle 2,1\rangle,\langle 2,2\rangle\}$ and $s_{\langle 2,2\rangle}>{{s_{\langle 2,1\rangle}+4}\choose {\lfloor (s_{\langle 2,1\rangle}+4)/2\rfloor}}$, then $H=(P_3\times P_2)(s_{\langle 2,2\rangle},s_{\langle 2,1\rangle},\overbrace{2,2,\ldots,2}^8)\in\mathscr{C}_1$.
\end{ppn}
\noindent\textit{Proof}: Let $V_1=\{(p,\langle 2,2\rangle)\mid p=1,2,\ldots,s_{\langle 2,2\rangle}\}$ and $V_2=\{(p,\langle 2,1\rangle)\mid p=1,2,\ldots,s_{\langle 2,1\rangle}\}$ $\cup \{(1,\langle 1,2\rangle),(2,\langle 1,2\rangle)\}\cup \{(1,\langle 3,2\rangle),(2,\langle 3,2\rangle)\}$. Note the subgraph induced by $V_1\cup V_2$ is a complete bipartite graph with partite sets $V_1$ and $V_2$. Since $|V_1|=s_{\langle 2,2\rangle}>{{s_{\langle 2,1\rangle}+4}\choose {\lfloor (s_{\langle 2,1\rangle}+4)/2\rfloor}}={{|V_2|}\choose{\lfloor|V_2|/2\rfloor}}$, there exists some $1\le p,q\le s_{\langle 2,2\rangle}$, $p\neq q$ such that $O((p,\langle 2,2\rangle))\subseteq O((q,\langle 2,2\rangle))$ by Sperner Theorem. This implies that $d_D((p,\langle 2,2\rangle),(q,\langle 2,2\rangle))\ge 4$. Hence, by Theorem \ref{thmA8.1.8}(b), $H\in\mathscr{C}_1$.
\qed

\indent\par We end off the section with a result on the hypercube graph.
\\
\\\textit{Proof of Proposition \ref{ppnA8.1.9}}: We shall prove $\bar{d}(Q_3^{(2)})=3=d(Q_3)$. Denote the vertices of the two disjoint copies of $C_4$ in $Q_3$ by $1,2,3,4,$ and $5,6,7,8$. Define an orientation $D$ of $Q_3^{(2)}$ as follows. (See Figure \ref{figA8.3.10}.)
\begin{align*}
& i\rightsquigarrow i+1 \rightsquigarrow i+2\rightsquigarrow i+3\text{ and }i\rightsquigarrow i+3 \text{ for } i=1,5,\\
& 4 \rightrightarrows 8, 2 \rightrightarrows 6, 5 \rightrightarrows 1,\text{ and } 7\rightrightarrows 3.
\end{align*}
It is easy to verify that $d(D)=3$. Hence, $Q_3^{(2)}\in \mathscr{C}_0$. Now, by Theorem \ref{thmA8.1.5}, $Q_\lambda^{(2)}\in \mathscr{C}_0$ for $\lambda\ge 3$. Since every vertex lies in a directed $C_4$, it follows from Lemma \ref{lemA8.1.2} that $Q_3(s_1,s_2,\ldots,s_n)\in\mathscr{C}_0\cup\mathscr{C}_1$ and $Q_\lambda(s_1,s_2,\ldots,s_n)\in\mathscr{C}_0$ for $\lambda\ge 4$.
\qed
\begin{center}
\tikzstyle{every node}=[circle, draw, fill=black!100,
inner sep=0pt, minimum width=5pt]
\begin{tikzpicture}[thick,scale=0.7]%

\draw(-8,1)node[label={[yshift=-0.1cm, xshift=0cm]270:{\small $2$}}](1_2){};
\draw(-10,-1)node[star,star points=4, label={[yshift=0.2cm] 270:{}}](2_2){};
\draw(-8,5)node[label={[yshift=-0.1cm, xshift=0.1cm]315:{\small $1$}}](1_1){};
\draw(-10,3)node[star,star points=4, label={[yshift=-0.2cm, xshift=0.3cm] 90:{}}](2_1){};
\draw(-12,5)node[label={[yshift=0.1cm]90:{\small $4$}}](1_4){};
\draw(-14,3)node[star,star points=4, label={[yshift=-0.2cm, xshift=0cm] 90:{}}](2_4){};
\draw(-12,1)node[label={[yshift=-0.1cm, xshift=0cm]270:{\small $3$}}](1_3){};
\draw(-14,-1)node[star,star points=4, label={[yshift=0.2cm] 270:{}}](2_3){};

\draw(0,4)node[label={[yshift=-0.1cm, xshift=0cm]270:{\small $6$}}](1_6){};
\draw(-2,2)node[star,star points=4, label={[yshift=0.1cm]270:{}}](2_6){};
\draw(0,8)node[label={[yshift=0.1cm]90:{\small $5$}}](1_5){};
\draw(-2,6)node[star,star points=4, label={[yshift=-0.2cm, xshift=0.3cm] 90:{}}](2_5){};
\draw(-4,8)node[label={[yshift=0.1cm]90:{\small $8$}}](1_8){};
\draw(-6,6)node[star,star points=4, label={[yshift=-0.2cm, xshift=0cm] 90:{}}](2_8){};
\draw(-4,4)node[label={[yshift=-0.1cm, xshift=0cm]270:{\small $7$}}](1_7){};
\draw(-6,2)node[star,star points=4, label={[yshift=0.7cm, xshift=0.2cm] 0:{}}](2_7){};

\draw[dashed, ->, line width=0.3mm, >=latex, shorten <= 0.2cm, shorten >= 0.15cm](1_2)--(1_3);
\draw[dashed, ->, line width=0.3mm, >=latex, shorten <= 0.2cm, shorten >= 0.15cm](2_2)--(2_3);
\draw[dashed, ->, line width=0.3mm, >=latex, shorten <= 0.2cm, shorten >= 0.15cm](1_3)--(2_2);
\draw[dashed, ->, line width=0.3mm, >=latex, shorten <= 0.2cm, shorten >= 0.15cm](2_3)--(1_2);

\draw[dashed, ->, line width=0.3mm, >=latex, shorten <= 0.2cm, shorten >= 0.15cm](1_3)--(1_4);
\draw[dashed, ->, line width=0.3mm, >=latex, shorten <= 0.2cm, shorten >= 0.15cm](2_3)--(2_4);
\draw[dashed, ->, line width=0.3mm, >=latex, shorten <= 0.2cm, shorten >= 0.15cm](1_4)--(2_3);
\draw[dashed, ->, line width=0.3mm, >=latex, shorten <= 0.2cm, shorten >= 0.15cm](2_4)--(1_3);

\draw[dashed, ->, line width=0.3mm, >=latex, shorten <= 0.2cm, shorten >= 0.15cm](1_1)--(1_2);
\draw[dashed, ->, line width=0.3mm, >=latex, shorten <= 0.2cm, shorten >= 0.15cm](2_1)--(2_2);
\draw[dashed, ->, line width=0.3mm, >=latex, shorten <= 0.2cm, shorten >= 0.15cm](1_2)--(2_1);
\draw[dashed, ->, line width=0.3mm, >=latex, shorten <= 0.2cm, shorten >= 0.15cm](2_2)--(1_1);

\draw[dashed, ->, line width=0.3mm, >=latex, shorten <= 0.2cm, shorten >= 0.15cm](1_1)--(1_4);
\draw[dashed, ->, line width=0.3mm, >=latex, shorten <= 0.2cm, shorten >= 0.15cm](2_1)--(2_4);
\draw[dashed, ->, line width=0.3mm, >=latex, shorten <= 0.2cm, shorten >= 0.15cm](1_4)--(2_1);
\draw[dashed, ->, line width=0.3mm, >=latex, shorten <= 0.2cm, shorten >= 0.15cm](2_4)--(1_1);

\draw[dashed, ->, line width=0.3mm, >=latex, shorten <= 0.2cm, shorten >= 0.15cm](1_6)--(1_7);
\draw[dashed, ->, line width=0.3mm, >=latex, shorten <= 0.2cm, shorten >= 0.15cm](2_6)--(2_7);
\draw[dashed, ->, line width=0.3mm, >=latex, shorten <= 0.2cm, shorten >= 0.15cm](1_7)--(2_6);
\draw[dashed, ->, line width=0.3mm, >=latex, shorten <= 0.2cm, shorten >= 0.15cm](2_7)--(1_6);

\draw[dashed, ->, line width=0.3mm, >=latex, shorten <= 0.2cm, shorten >= 0.15cm](1_7)--(1_8);
\draw[dashed, ->, line width=0.3mm, >=latex, shorten <= 0.2cm, shorten >= 0.15cm](2_7)--(2_8);
\draw[dashed, ->, line width=0.3mm, >=latex, shorten <= 0.2cm, shorten >= 0.15cm](1_8)--(2_7);
\draw[dashed, ->, line width=0.3mm, >=latex, shorten <= 0.2cm, shorten >= 0.15cm](2_8)--(1_7);

\draw[dashed, ->, line width=0.3mm, >=latex, shorten <= 0.2cm, shorten >= 0.15cm](1_5)--(1_6);
\draw[dashed, ->, line width=0.3mm, >=latex, shorten <= 0.2cm, shorten >= 0.15cm](2_5)--(2_6);
\draw[dashed, ->, line width=0.3mm, >=latex, shorten <= 0.2cm, shorten >= 0.15cm](1_6)--(2_5);
\draw[dashed, ->, line width=0.3mm, >=latex, shorten <= 0.2cm, shorten >= 0.15cm](2_6)--(1_5);

\draw[dashed, ->, line width=0.3mm, >=latex, shorten <= 0.2cm, shorten >= 0.15cm](1_5)--(1_8);
\draw[dashed, ->, line width=0.3mm, >=latex, shorten <= 0.2cm, shorten >= 0.15cm](2_5)--(2_8);
\draw[dashed, ->, line width=0.3mm, >=latex, shorten <= 0.2cm, shorten >= 0.15cm](1_8)--(2_5);
\draw[dashed, ->, line width=0.3mm, >=latex, shorten <= 0.2cm, shorten >= 0.15cm](2_8)--(1_5);

\draw[densely dotted, ->, line width=0.3mm, >=latex, shorten <= 0.2cm, shorten >= 0.15cm](1_4)--(1_8);
\draw[densely dotted, ->, line width=0.3mm, >=latex, shorten <= 0.2cm, shorten >= 0.15cm](2_4)--(2_8);
\draw[densely dotted, ->, line width=0.3mm, >=latex, shorten <= 0.2cm, shorten >= 0.15cm](1_4)--(2_8);
\draw[densely dotted, ->, line width=0.3mm, >=latex, shorten <= 0.2cm, shorten >= 0.15cm](2_4)--(1_8);

\draw[densely dotted, ->, line width=0.3mm, >=latex, shorten <= 0.2cm, shorten >= 0.15cm](1_2)--(1_6);
\draw[densely dotted, ->, line width=0.3mm, >=latex, shorten <= 0.2cm, shorten >= 0.15cm](2_2)--(2_6);
\draw[densely dotted, ->, line width=0.3mm, >=latex, shorten <= 0.2cm, shorten >= 0.15cm](1_2)--(2_6);
\draw[densely dotted, ->, line width=0.3mm, >=latex, shorten <= 0.2cm, shorten >= 0.15cm](2_2)--(1_6);

\draw[densely dotted, ->, line width=0.3mm, >=latex, shorten <= 0.2cm, shorten >= 0.15cm](1_5)--(1_1);
\draw[densely dotted, ->, line width=0.3mm, >=latex, shorten <= 0.2cm, shorten >= 0.15cm](2_5)--(2_1);
\draw[densely dotted, ->, line width=0.3mm, >=latex, shorten <= 0.2cm, shorten >= 0.15cm](1_5)--(2_1);
\draw[densely dotted, ->, line width=0.3mm, >=latex, shorten <= 0.2cm, shorten >= 0.15cm](2_5)--(1_1);

\draw[densely dotted, ->, line width=0.3mm, >=latex, shorten <= 0.2cm, shorten >= 0.15cm](1_7)--(1_3);
\draw[densely dotted, ->, line width=0.3mm, >=latex, shorten <= 0.2cm, shorten >= 0.15cm](2_7)--(2_3);
\draw[densely dotted, ->, line width=0.3mm, >=latex, shorten <= 0.2cm, shorten >= 0.15cm](1_7)--(2_3);
\draw[densely dotted, ->, line width=0.3mm, >=latex, shorten <= 0.2cm, shorten >= 0.15cm](2_7)--(1_3);
\end{tikzpicture}
{
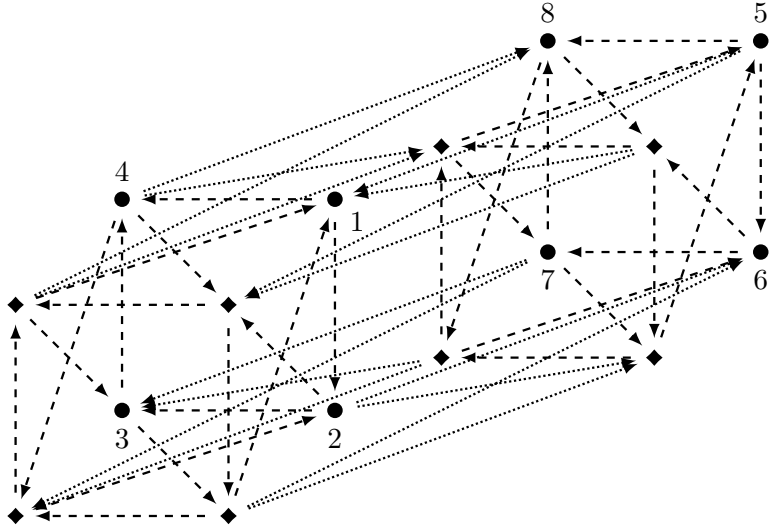
\captionof{figure}{Orientation $D$ of $Q_3^{(2)}$, where $d(D)=3$.}\label{figA8.3.19}}
\end{center}

\section{Cartesian product of trees with cycles $T_\lambda\times C_\mu$}
In this section, we consider Cartesian product of trees with cycles.
\\
\\\textit{Proof of Theorem \ref{thmA8.1.10}}:
\\Case 1. $\lambda\ge 2$ and $\mu\ge 4$.
\indent\par Let $(V_1, V_2)$ be a bipartition of $V(T_\lambda)$, i.e. $V_1$ and $V_2$ are independent sets. Let $F$ be a strong orientation of $C_\mu$, say $1\rightarrow 2\rightarrow \cdots \rightarrow \mu\rightarrow 1$, and define an orientation $D$ for $(T_\lambda\times C_\mu)^{(2)}$ as follows.
\begin{align*}
\langle u,x\rangle\rightrightarrows \langle u,y\rangle\iff x\rightarrow y \text{ in }F
\end{align*}
for any $u\in V_1$ and any $x,y\in V(C_\mu)$, i.e. the copy $C_\mu^{(2)}$ is oriented similarly to $F$.
\begin{align*}
\langle u,x\rangle\rightrightarrows \langle u,y\rangle\iff y\rightarrow x \text{ in }F
\end{align*}
for any $u\in V_2$ and any $x,y\in V(C_\mu)$, i.e. the copy $C_\mu^{(2)}$ is oriented similarly to $\tilde{F}$.
\begin{align*}
\langle u,x\rangle\rightsquigarrow\langle v,x\rangle
\end{align*}
for any $u,v\in V(T_\lambda)$ with $uv\in E(T_\lambda)$ and any $x\in V(C_\mu)$.
\indent\par We claim that $d_D((p,\langle u,x\rangle),(q,\langle v,y\rangle))\le \lambda+\lfloor\frac{\mu}{2}\rfloor=d(T_\lambda\times C_\mu)$ for any $\langle u,x\rangle,\langle v,y\rangle\in V(T_\lambda\times C_\mu)$, and $p,q=1,2$. Suppose $u=v\in V_1$. Note that either $d_F(x,y)\le \lfloor\frac{\mu}{2}\rfloor$ or $d_{\tilde{F}}(x,y)\le \lfloor\frac{\mu}{2}\rfloor$. So, there exist paths $P$ and $P'$ in $D$, each of length at most $\lfloor\frac{\mu}{2}\rfloor$, from $\{(1,\langle u,x\rangle),(2,\langle u,x\rangle)\}$ to $\{(1,\langle u,y\rangle),(2,\langle u,y\rangle)\}$ and from $\{(1,\langle w,x\rangle),(2,\langle w,x\rangle)\}$ to $\{(1,\langle w,y\rangle),$ $(2,\langle w,y\rangle)\}$, where $w\in V_2$ is some vertex adjacent to $u$ in $T_\lambda$ respectively. In the former case, $P$ suffices and we are done. In the latter case, we shall further assume $\langle u,x\rangle\rightsquigarrow\langle w,x\rangle$ for simplicity; the proof is similar otherwise. Then, $(p,\langle u,x\rangle)(p,\langle w,x\rangle)$, $P'$ and $(3-q,\langle w,y\rangle)(q,\langle u,y\rangle)$ form a $(p,\langle u,x\rangle)-(q,\langle v,y\rangle)$ path of length at most $2+\lfloor\frac{\mu}{2}\rfloor\le \lambda+\lfloor\frac{\mu}{2}\rfloor$. A similar proof follows if $u=v\in V_2$.
\indent\par Suppose $u\neq v$. Let $u w_1 w_2 \ldots w_l v$ be the unique shortest $u-v$ path in $T_\lambda$. For simplicity, we shall assume $\langle u,x\rangle\rightsquigarrow\langle w_1,x\rangle\rightsquigarrow\cdots\rightsquigarrow\langle v,x\rangle$; the proof is similar otherwise. If $x=y$, then $(p,\langle u,x\rangle)(p,\langle w_1,x\rangle)\ldots (p,\langle v,x\rangle)(3-p,\langle v,w_l\rangle)(3-p,\langle v,x\rangle)$ guarantees a $(p,\langle u,x\rangle)-(q,\langle v,y\rangle)$ path of length at most $\lambda+2\le \lambda+\lfloor\frac{\mu}{2}\rfloor$.
\indent\par Next, suppose $x\neq y$. Futhermore, we shall assume $v\in V_1$ (and hence $w_l\in V_2$); the proof is similar if $v\in V_2$. Again, consider the cases $d_F(x,y)\le \lfloor\frac{\mu}{2}\rfloor$ or $d_{\tilde{F}}(x,y)\le \lfloor\frac{\mu}{2}\rfloor$. So, there exist paths $Q$ and $Q'$ in $D$, each of length at most $\lfloor\frac{\mu}{2}\rfloor$, from $\{(1,\langle v,x\rangle),(2,\langle v,x\rangle)\}$ to $\{(1,\langle v,y\rangle),(2,\langle v,y\rangle)\}$ and from $\{(1,\langle w_l,x\rangle),(2,\langle w_l,x\rangle)\}$ to $\{(1,\langle w_l,y\rangle),(2,\langle w_l,y\rangle)\}$ respectively. In the former case, $(p,\langle u,x\rangle)(p,\langle w_1,x\rangle)\ldots (p,\langle v,x\rangle)$ and $Q$ form a $(p,\langle u,x\rangle)-(q,\langle v,y\rangle)$ path of length at most $\lambda+\lfloor\frac{\mu}{2}\rfloor$. In the latter case, $(p,\langle u,x\rangle)(p,\langle w_1,x\rangle)\ldots$ $(p,\langle w_l,x\rangle)$ with $Q'$ and $(q,\langle w_l,y\rangle) (q,\langle v,y\rangle)$ form a $(p,\langle u,x\rangle)-(q,\langle v,y\rangle)$ of length at most $\lambda+\lfloor\frac{\mu}{2}\rfloor$. Hence, $(T_\lambda\times C_\mu)^{(2)}\in \mathscr{C}_0$. 
\\
\\Case 2. $\lambda=\mu=3$.
\indent\par Define an orientation $D$ for $(T_3\times C_3)^{(2)}$ as follows. (See Figure \ref{figA8.4.20}.) For all $[i]_1\in N_T(\mathtt{c}_1)-\{\mathtt{c}_2\}$ and all $[j]_2 \in N_T(\mathtt{c}_2)-\{\mathtt{c}_1\}$,
\begin{align*}
&\langle\mathtt{c}_1,1\rangle\rightrightarrows \langle\mathtt{c}_1,2\rangle\rightrightarrows \langle\mathtt{c}_1,3\rangle\rightrightarrows \langle\mathtt{c}_1,1\rangle,\ \langle\mathtt{c}_2,3\rangle\rightrightarrows \langle\mathtt{c}_2,2\rangle\rightrightarrows \langle\mathtt{c}_2,1\rangle\rightrightarrows \langle\mathtt{c}_2,3\rangle,\\
&\langle [i]_1,y\rangle\overset{1}\twoheadrightarrow \langle\mathtt{c}_1,y\rangle \overset{2}\twoheadleftarrow \langle\mathtt{c}_2,y\rangle, \text{ and } \langle\mathtt{c}_2,y\rangle\rightsquigarrow \langle [j]_2,y\rangle \text{ for all } y=1,2,3,\\
&\langle [i]_1,1\rangle\rightsquigarrow \langle [i]_1,2\rangle\rightsquigarrow \langle [i]_1,3\rangle\rightsquigarrow \langle [i]_1,1\rangle, \text{ and }\langle [j]_2,1\rangle\rightsquigarrow \langle [j]_2,2\rangle\rightsquigarrow \langle [j]_2,3\rangle\rightsquigarrow \langle [j]_2,1\rangle.
\end{align*}
It is straightforward to verify that $d(D)=4$. In view of the symmetry of $D$, it suffices to check $D$ for $(T_3\times C_3)^{(2)}$ where $\mathtt{c}_i$ has two end-vertex neighbours $[1]_i,[2]_i$ for $i=1,2$ in $T_3$. We remark that the checking includes the distance from any vertex in the $[1]_1$-copy ($[1]_2$-copy resp.) of $C_3^{(2)}$ to any vertex in the $[2]_1$-copy ($[2]_2$-copy resp.) of $C_3^{(2)}$, though only one $[i]_1$-copy ($[j]_2$-copy resp.) is shown in Figure \ref{figA8.4.20} for brevity. Hence, $(T_3\times C_3)^{(2)}\in \mathscr{C}_0$.
\noindent\par Since every vertex lies in a directed $C_4$ in $D$ of both cases, it follows from Lemma \ref{lemA8.1.2} that $(T_\lambda\times C_\mu)(s_1,s_2,\ldots,s_n)\in\mathscr{C}_0$.
\qed

\noindent\makebox[\textwidth]{%
\tikzstyle{every node}=[circle, draw, fill=black!100,
                       inner sep=0pt, minimum width=5pt]
\begin{tikzpicture}[thick,scale=0.6]%
\draw(-11,1)node[star,star points=4, label={[yshift=0.5cm, xshift=0.2cm] 270:{}}](2_12){};
\draw(-5,1)node[star,star points=4, label={[yshift=0.5cm] 270:{}}](2_22){};
\draw(1,1)node[star,star points=4, label={[yshift=0.5cm]270:{}}](2_32){};
\draw(7,1)node[star,star points=4, label={[yshift=0.6cm]270:{}}](2_42){};

\draw(-9,2)node[label={[yshift=-0.2cm, xshift=0cm]135:{\small $\langle [i]_1,2\rangle$}}](1_12){};
\draw(-3,2)node[label={[yshift=-0.15cm, xshift=0cm]135:{\small $\langle\mathtt{c}_1,2\rangle$}}](1_22){};
\draw(3,2)node[label={[yshift=-0.15cm, xshift=0cm]135:{\small $\langle\mathtt{c}_2,2\rangle$}}](1_32){};
\draw(9,2)node[label={[yshift=-0.2cm, xshift=0cm]135:{\small $\langle [j]_2,2\rangle$}}](1_42){};

\draw(-8,5)node[star,star points=4, label={[yshift=0cm, xshift=-0.3cm] 180:{}}](2_13){};
\draw(-2,5)node[star,star points=4, label={[yshift=-0.4cm, xshift=-0.3cm] 180:{}}](2_23){};
\draw(4,5)node[star,star points=4, label={[yshift=-0.4cm, xshift=-0.3cm] 180:{}}](2_33){};
\draw(10,5)node[star,star points=4, label={[yshift=-0.4cm, xshift=-0.3cm] 180:{}}](2_43){};

\draw(-6,6)node[label={[yshift=-0.4cm]90:{\small $\langle [i]_1,3\rangle$}}](1_13){};
\draw(0,6)node[label={[yshift=-0.3cm]90:{\small $\langle\mathtt{c}_1,3\rangle$}}](1_23){};
\draw(6,6)node[label={[yshift=-0.3cm]90:{\small $\langle\mathtt{c}_2,3\rangle$}}](1_33){};
\draw(12,6)node[label={[yshift=-0.4cm]90:{\small $\langle [j]_2,3\rangle$}}](1_43){};

\draw(-8,-2)node[star,star points=4, label={[yshift=0.7cm] 270:{}}](2_11){};
\draw(-2,-2)node[star,star points=4, label={[yshift=0.6cm] 270:{}}](2_21){};
\draw(4,-2)node[star,star points=4, label={[yshift=0.6cm]270:{}}](2_31){};
\draw(10,-2)node[star,star points=4, label={[yshift=0.7cm]270:{}}](2_41){};

\draw(-6,-1)node[label={[yshift=-0.2cm, xshift=0.2cm]45:{\small $\langle [i]_1,1\rangle$}}](1_11){};
\draw(0,-1)node[label={[yshift=-0.15cm, xshift=0.2cm]45:{\small $\langle\mathtt{c}_1,1\rangle$}}](1_21){};
\draw(6,-1)node[label={[yshift=-0.15cm, xshift=0.2cm]45:{\small $\langle\mathtt{c}_2,1\rangle$}}](1_31){};
\draw(12,-1)node[label={[yshift=0.3cm, xshift=0cm]270:{\small $\langle [j]_2,1\rangle$}}](1_41){};

\draw[->, line width=0.3mm, >=latex, shorten <= 0.2cm, shorten >= 0.15cm](1_12)--(1_22);
\draw[->, line width=0.3mm, >=latex, shorten <= 0.2cm, shorten >= 0.15cm](2_12)--(1_22);
\draw[->, line width=0.3mm, >=latex, shorten <= 0.2cm, shorten >= 0.15cm](1_22)--(1_32);
\draw[->, line width=0.3mm, >=latex, shorten <= 0.2cm, shorten >= 0.15cm](1_22)--(2_32);
\draw[->, line width=0.3mm, >=latex, shorten <= 0.2cm, shorten >= 0.15cm](2_22)--(1_12);
\draw[->, line width=0.3mm, >=latex, shorten <= 0.2cm, shorten >= 0.15cm](2_22)--(2_12);
\draw[->, line width=0.3mm, >=latex, shorten <= 0.2cm, shorten >= 0.15cm](2_32)--(2_22);
\draw[->, line width=0.3mm, >=latex, shorten <= 0.2cm, shorten >= 0.15cm](1_32)--(2_22);

\draw[->, line width=0.3mm, >=latex, shorten <= 0.2cm, shorten >= 0.15cm](1_13)--(1_23);
\draw[->, line width=0.3mm, >=latex, shorten <= 0.2cm, shorten >= 0.15cm](2_13)--(1_23);
\draw[->, line width=0.3mm, >=latex, shorten <= 0.2cm, shorten >= 0.15cm](1_23)--(1_33);
\draw[->, line width=0.3mm, >=latex, shorten <= 0.2cm, shorten >= 0.15cm](1_23)--(2_33);
\draw[->, line width=0.3mm, >=latex, shorten <= 0.2cm, shorten >= 0.15cm](2_23)--(1_13);
\draw[->, line width=0.3mm, >=latex, shorten <= 0.2cm, shorten >= 0.15cm](2_23)--(2_13);
\draw[->, line width=0.3mm, >=latex, shorten <= 0.2cm, shorten >= 0.15cm](2_33)--(2_23);
\draw[->, line width=0.3mm, >=latex, shorten <= 0.2cm, shorten >= 0.15cm](1_33)--(2_23);

\draw[->, line width=0.3mm, >=latex, shorten <= 0.2cm, shorten >= 0.15cm](1_11)--(1_21);
\draw[->, line width=0.3mm, >=latex, shorten <= 0.2cm, shorten >= 0.15cm](2_11)--(1_21);
\draw[->, line width=0.3mm, >=latex, shorten <= 0.2cm, shorten >= 0.15cm](1_21)--(1_31);
\draw[->, line width=0.3mm, >=latex, shorten <= 0.2cm, shorten >= 0.15cm](1_21)--(2_31);
\draw[->, line width=0.3mm, >=latex, shorten <= 0.2cm, shorten >= 0.15cm](2_21)--(1_11);
\draw[->, line width=0.3mm, >=latex, shorten <= 0.2cm, shorten >= 0.15cm](2_21)--(2_11);
\draw[->, line width=0.3mm, >=latex, shorten <= 0.2cm, shorten >= 0.15cm](2_31)--(2_21);
\draw[->, line width=0.3mm, >=latex, shorten <= 0.2cm, shorten >= 0.15cm](1_31)--(2_21);

\draw[dashed, ->, line width=0.3mm, >=latex, shorten <= 0.2cm, shorten >= 0.15cm](1_12)--(1_13);
\draw[dashed, ->, line width=0.3mm, >=latex, shorten <= 0.2cm, shorten >= 0.15cm](2_12) to [out=100, in=190] (2_13);
\draw[dashed, ->, line width=0.3mm, >=latex, shorten <= 0.2cm, shorten >= 0.15cm](1_13) to [out=175, in=120] (2_12);
\draw[dashed, ->, line width=0.3mm, >=latex, shorten <= 0.2cm, shorten >= 0.15cm](2_13)--(1_12);

\draw[dashed, ->, line width=0.3mm, >=latex, shorten <= 0.2cm, shorten >= 0.15cm](2_11)--(2_12);
\draw[dashed, ->, line width=0.3mm, >=latex, shorten <= 0.2cm, shorten >= 0.15cm](1_11)--(1_12);
\draw[dashed, ->, line width=0.3mm, >=latex, shorten <= 0.2cm, shorten >= 0.15cm](1_12)--(2_11);
\draw[dashed, ->, line width=0.3mm, >=latex, shorten <= 0.2cm, shorten >= 0.15cm](2_12)--(1_11);

\draw[dashed, ->, line width=0.3mm, >=latex, shorten <= 0.2cm, shorten >= 0.15cm](1_11)--(2_13);
\draw[dashed, ->, line width=0.3mm, >=latex, shorten <= 0.2cm, shorten >= 0.15cm](2_11)--(1_13);
\draw[dashed, ->, line width=0.3mm, >=latex, shorten <= 0.2cm, shorten >= 0.15cm](1_13)--(1_11);
\draw[dashed, ->, line width=0.3mm, >=latex, shorten <= 0.2cm, shorten >= 0.15cm](2_13)--(2_11);

\draw[densely dotted, ->, line width=0.3mm, >=latex, shorten <= 0.2cm, shorten >= 0.15cm](1_22)--(1_23);
\draw[densely dotted, ->, line width=0.3mm, >=latex, shorten <= 0.2cm, shorten >= 0.15cm](2_22) to [out=95, in=205] (2_23);
\draw[densely dotted, ->, line width=0.3mm, >=latex, shorten <= 0.2cm, shorten >= 0.15cm](2_22) to [out=120, in=195] (1_23);
\draw[densely dotted, ->, line width=0.3mm, >=latex, shorten <= 0.2cm, shorten >= 0.15cm](1_22)--(2_23);

\draw[densely dotted, ->, line width=0.3mm, >=latex, shorten <= 0.2cm, shorten >= 0.15cm](2_21)--(2_22);
\draw[densely dotted, ->, line width=0.3mm, >=latex, shorten <= 0.2cm, shorten >= 0.15cm](1_21)--(1_22);
\draw[densely dotted, ->, line width=0.3mm, >=latex, shorten <= 0.2cm, shorten >= 0.15cm](2_21)--(1_22);
\draw[densely dotted, ->, line width=0.3mm, >=latex, shorten <= 0.2cm, shorten >= 0.15cm](1_21)--(2_22);

\draw[densely dotted, ->, line width=0.3mm, >=latex, shorten <= 0.2cm, shorten >= 0.15cm](1_23)--(1_21);
\draw[densely dotted, ->, line width=0.3mm, >=latex, shorten <= 0.2cm, shorten >= 0.15cm](2_23)--(2_21);
\draw[densely dotted, ->, line width=0.3mm, >=latex, shorten <= 0.2cm, shorten >= 0.15cm](1_23)--(2_21);
\draw[densely dotted, ->, line width=0.3mm, >=latex, shorten <= 0.2cm, shorten >= 0.15cm](2_23)--(1_21);

\draw[densely dotted, ->, line width=0.3mm, >=latex, shorten <= 0.2cm, shorten >= 0.15cm](2_32)--(2_31);
\draw[densely dotted, ->, line width=0.3mm, >=latex, shorten <= 0.2cm, shorten >= 0.15cm](1_32)--(1_31);
\draw[densely dotted, ->, line width=0.3mm, >=latex, shorten <= 0.2cm, shorten >= 0.15cm](1_32)--(2_31);
\draw[densely dotted, ->, line width=0.3mm, >=latex, shorten <= 0.2cm, shorten >= 0.15cm](2_32)--(1_31);

\draw[densely dotted, ->, line width=0.3mm, >=latex, shorten <= 0.2cm, shorten >= 0.15cm](1_33)--(1_32);
\draw[densely dotted, ->, line width=0.3mm, >=latex, shorten <= 0.2cm, shorten >= 0.15cm](2_33) to [out=205, in=95](2_32);
\draw[densely dotted, ->, line width=0.3mm, >=latex, shorten <= 0.2cm, shorten >= 0.15cm](1_33) to [out=195, in=120](2_32);
\draw[densely dotted, ->, line width=0.3mm, >=latex, shorten <= 0.2cm, shorten >= 0.15cm](2_33)--(1_32);

\draw[densely dotted, ->, line width=0.3mm, >=latex, shorten <= 0.2cm, shorten >= 0.15cm](1_31)--(2_33);
\draw[densely dotted, ->, line width=0.3mm, >=latex, shorten <= 0.2cm, shorten >= 0.15cm](2_31)--(1_33);
\draw[densely dotted, ->, line width=0.3mm, >=latex, shorten <= 0.2cm, shorten >= 0.15cm](1_31)--(1_33);
\draw[densely dotted, ->, line width=0.3mm, >=latex, shorten <= 0.2cm, shorten >= 0.15cm](2_31)--(2_33);

\draw[dashed, ->, line width=0.3mm, >=latex, shorten <= 0.2cm, shorten >= 0.15cm](1_42)--(1_43);
\draw[dashed, ->, line width=0.3mm, >=latex, shorten <= 0.2cm, shorten >= 0.15cm](2_42) to [out=100, in=200] (2_43);
\draw[dashed, ->, line width=0.3mm, >=latex, shorten <= 0.2cm, shorten >= 0.15cm](1_43) to [out=190, in=120] (2_42);
\draw[dashed, ->, line width=0.3mm, >=latex, shorten <= 0.2cm, shorten >= 0.15cm](2_43)--(1_42);

\draw[dashed, ->, line width=0.3mm, >=latex, shorten <= 0.2cm, shorten >= 0.15cm](2_41)--(2_42);
\draw[dashed, ->, line width=0.3mm, >=latex, shorten <= 0.2cm, shorten >= 0.15cm](1_41)--(1_42);
\draw[dashed, ->, line width=0.3mm, >=latex, shorten <= 0.2cm, shorten >= 0.15cm](1_42)--(2_41);
\draw[dashed, ->, line width=0.3mm, >=latex, shorten <= 0.2cm, shorten >= 0.15cm](2_42)--(1_41);

\draw[dashed, ->, line width=0.3mm, >=latex, shorten <= 0.2cm, shorten >= 0.15cm](1_43)--(1_41);
\draw[dashed, ->, line width=0.3mm, >=latex, shorten <= 0.2cm, shorten >= 0.15cm](2_43)--(2_41);
\draw[dashed, ->, line width=0.3mm, >=latex, shorten <= 0.2cm, shorten >= 0.15cm](1_41)--(2_43);
\draw[dashed, ->, line width=0.3mm, >=latex, shorten <= 0.2cm, shorten >= 0.15cm](2_41)--(1_43);

\draw[dashed, ->, line width=0.3mm, >=latex, shorten <= 0.2cm, shorten >= 0.15cm](1_33)--(1_43);
\draw[dashed, ->, line width=0.3mm, >=latex, shorten <= 0.2cm, shorten >= 0.15cm](2_33)--(2_43);
\draw[dashed, ->, line width=0.3mm, >=latex, shorten <= 0.2cm, shorten >= 0.15cm](2_43)--(1_33);
\draw[dashed, ->, line width=0.3mm, >=latex, shorten <= 0.2cm, shorten >= 0.15cm](1_43)--(2_33);

\draw[dashed, ->, line width=0.3mm, >=latex, shorten <= 0.2cm, shorten >= 0.15cm](1_32)--(1_42);
\draw[dashed, ->, line width=0.3mm, >=latex, shorten <= 0.2cm, shorten >= 0.15cm](2_32)--(2_42);
\draw[dashed, ->, line width=0.3mm, >=latex, shorten <= 0.2cm, shorten >= 0.15cm](2_42)--(1_32);
\draw[dashed, ->, line width=0.3mm, >=latex, shorten <= 0.2cm, shorten >= 0.15cm](1_42)--(2_32);

\draw[dashed, ->, line width=0.3mm, >=latex, shorten <= 0.2cm, shorten >= 0.15cm](1_31)--(1_41);
\draw[dashed, ->, line width=0.3mm, >=latex, shorten <= 0.2cm, shorten >= 0.15cm](2_31)--(2_41);
\draw[dashed, ->, line width=0.3mm, >=latex, shorten <= 0.2cm, shorten >= 0.15cm](2_41)--(1_31);
\draw[dashed, ->, line width=0.3mm, >=latex, shorten <= 0.2cm, shorten >= 0.15cm](1_41)--(2_31);
\end{tikzpicture}
}\captionsetup{justification=centering}
{
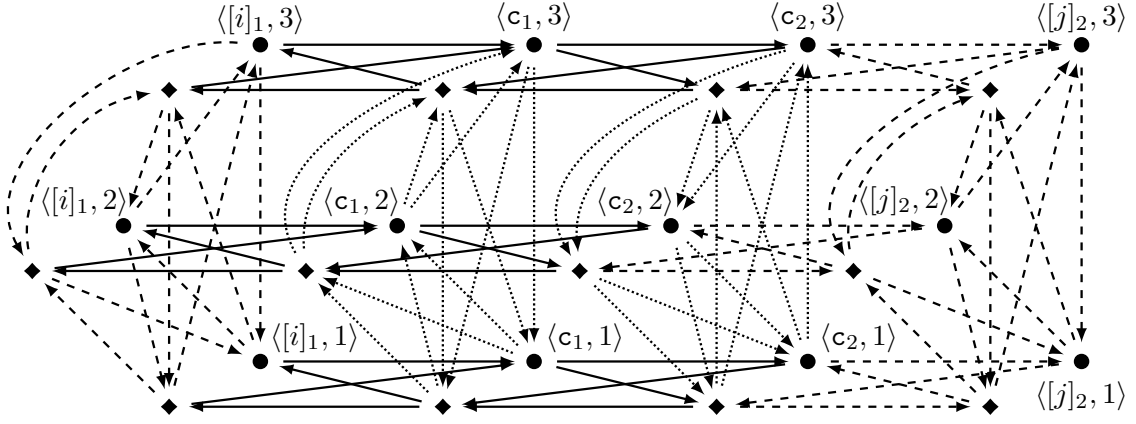
\captionof{figure}{Partial orientation $D$ of $(T_3\times C_3)^{(2)}$,\\where $[i]_1\in N_T(\mathtt{c}_1)-\{\mathtt{c}_2\}$ and $[j]_2\in N_T(\mathtt{c}_2)-\{\mathtt{c}_1\}$ and $d(D)=4$.}\label{figA8.4.20}}

\indent\par Next, we want to consider $T_2\times C_3$ and $P_2\times C_3$. Instead, we shall prove more general results involving $K_\mu$, $\mu\ge 3$, in place of $C_3$. For $T_2\times K_\mu$, we split into cases of $deg_{T_2}(\mathtt{c})=2$ (i.e. $T_2=P_3$) and $deg_{T_2}(\mathtt{c})>2$.

\begin{ppn}\label{ppnA8.4.1}
For $\mu\ge 3$,
\\(a) if $deg_{T_2}(\mathtt{c})=2$, then $(T_2\times K_\mu)(s_1,s_2,\ldots,s_n)\in\mathscr{C}_0$.
\\(b) if $deg_{T_2}(\mathtt{c})>2$, then $(T_2\times K_\mu)^{(2)}\in\mathscr{C}_1$ and $(T_2\times K_\mu)(s_1,s_2,\ldots,s_n)\in\mathscr{C}_0\cup \mathscr{C}_1$.
\end{ppn}
\noindent\textit{Proof}: Note that $d(T_2\times K_\mu)=3$. Define an orientation $D$ for $(T_2\times K_\mu)^{(2)}$ as follows. (See Figure \ref{figA8.4.21} when $deg_{T_2}(\mathtt{c})=2$ and $\mu=3$.)
\begin{align}
&\langle [1],j\rangle\overset{1}\twoheadrightarrow \langle \mathtt{c},j\rangle\overset{2}\twoheadleftarrow \langle [i],j\rangle \label{eqA8.4.1}
\end{align}
for all $[i]\in N_{T_2}(\mathtt{c})-\{[1]\}$ and $j=1,2,\ldots, \mu$.
\begin{equation}
\left.\begin{array}{@{}ll@{}}
&\langle v,j_1\rangle\rightsquigarrow \langle v,j_2\rangle \text{ whenever } 2\le j_1<j_2\le\mu,\text{ and}\\
&\langle v,j\rangle\rightsquigarrow \langle v,1\rangle\rightsquigarrow \langle v,2\rangle \text{ for }j=3,4,\ldots,\mu,
  \end{array}\right\}
\label{eqA8.4.2}
\end{equation}
for all $v\in V(T_2)$.
\indent\par (a) Suppose $deg_{T_2}(\mathtt{c})=2$. We give a brief verification of $d(D)=3$. It is easy to check that the orientation in Figure \ref{figA8.4.21} has diameter 3. Next, note for all $v\in T_2$, and all $j=4,5,\ldots,\mu$, that $\langle v,j\rangle$ really plays the same role as $\langle v,3\rangle$ in view of (\ref{eqA8.4.2}). Hence, it remains to check that the distance of any two vertices in each copy of $K_\mu^{(2)}$ is at most 3. This follows since $u\rightsquigarrow v$ or $v\rightsquigarrow u$ for all $u,v$ in each copy of $K_\mu$. Hence, $(T_2\times  K_\mu)^{(2)}\in \mathscr{C}_0$.  Since every vertex lies in a directed $C_3$, it follows from Lemma \ref{lemA8.1.2} that $(T_2\times K_\mu)(s_1,s_2,\ldots,s_n)\in\mathscr{C}_0$.

\begin{center}
\tikzstyle{every node}=[circle, draw, fill=black!100,
                       inner sep=0pt, minimum width=5pt]
\begin{tikzpicture}[thick,scale=0.7]%
\draw(-15,1)node[star,star points=4, label={[yshift=0.2cm, xshift=-0.3cm] 270:{}}](2_12){};
\draw(-9,1)node[star,star points=4, label={[yshift=0.4cm] 270:{}}](2_22){};
\draw(-3,1)node[star,star points=4, label={[yshift=0.5cm]270:{}}](2_32){};

\draw(-13,2)node[label={[yshift=-0.15cm, xshift=0cm]135:{\small $\langle [1],2\rangle$}}](1_12){};
\draw(-7,2)node[label={[yshift=-0.1cm, xshift=0cm]135:{\small $\langle \mathtt{c},2\rangle$}}](1_22){};
\draw(-1,2)node[label={[yshift=-0.15cm, xshift=0cm]135:{\small $\langle [2],2\rangle$}}](1_32){};

\draw(-12,5)node[star,star points=4, label={[yshift=0cm, xshift=-0.2cm] 180:{}}](2_13){};
\draw(-6,5)node[star,star points=4, label={[yshift=-0.4cm, xshift=-0.3cm] 180:{}}](2_23){};
\draw(-0,5)node[star,star points=4, label={[yshift=-0.4cm, xshift=-0.3cm] 180:{}}](2_33){};

\draw(-10,6)node[label={[yshift=-0.3cm]90:{\small $\langle [1],3\rangle$}}](1_13){};
\draw(-4,6)node[label={[yshift=-0.2cm]90:{\small $\langle \mathtt{c},3\rangle$}}](1_23){};
\draw(2,6)node[label={[yshift=-0.3cm]90:{\small $\langle [2],3\rangle$}}](1_33){};

\draw(-12,-2)node[star,star points=4, label={[yshift=0.5cm] 270:{}}](2_11){};
\draw(-6,-2)node[star,star points=4, label={[yshift=0.4cm] 270:{}}](2_21){};
\draw(0,-2)node[star,star points=4, label={[yshift=0.5cm]270:{}}](2_31){};

\draw(-10,-1)node[label={[yshift=-0.15cm, xshift=0.2cm]45:{\small $\langle [1],1\rangle$}}](1_11){};
\draw(-4,-1)node[label={[yshift=-0.1cm, xshift=0.2cm]45:{\small $\langle \mathtt{c},1\rangle$}}](1_21){};
\draw(2,-1)node[label={[yshift=0.3cm, xshift=0cm]270:{\small $\langle [2],1\rangle$}}](1_31){};

\draw[->, line width=0.3mm, >=latex, shorten <= 0.2cm, shorten >= 0.15cm](1_12)--(1_22);
\draw[->, line width=0.3mm, >=latex, shorten <= 0.2cm, shorten >= 0.15cm](2_12)--(1_22);
\draw[->, line width=0.3mm, >=latex, shorten <= 0.2cm, shorten >= 0.15cm](1_22)--(1_32);
\draw[->, line width=0.3mm, >=latex, shorten <= 0.2cm, shorten >= 0.15cm](1_22)--(2_32);
\draw[->, line width=0.3mm, >=latex, shorten <= 0.2cm, shorten >= 0.15cm](2_22)--(1_12);
\draw[->, line width=0.3mm, >=latex, shorten <= 0.2cm, shorten >= 0.15cm](2_22)--(2_12);
\draw[->, line width=0.3mm, >=latex, shorten <= 0.2cm, shorten >= 0.15cm](2_32)--(2_22);
\draw[->, line width=0.3mm, >=latex, shorten <= 0.2cm, shorten >= 0.15cm](1_32)--(2_22);

\draw[->, line width=0.3mm, >=latex, shorten <= 0.2cm, shorten >= 0.15cm](1_13)--(1_23);
\draw[->, line width=0.3mm, >=latex, shorten <= 0.2cm, shorten >= 0.15cm](2_13)--(1_23);
\draw[->, line width=0.3mm, >=latex, shorten <= 0.2cm, shorten >= 0.15cm](1_23)--(1_33);
\draw[->, line width=0.3mm, >=latex, shorten <= 0.2cm, shorten >= 0.15cm](1_23)--(2_33);
\draw[->, line width=0.3mm, >=latex, shorten <= 0.2cm, shorten >= 0.15cm](2_23)--(1_13);
\draw[->, line width=0.3mm, >=latex, shorten <= 0.2cm, shorten >= 0.15cm](2_23)--(2_13);
\draw[->, line width=0.3mm, >=latex, shorten <= 0.2cm, shorten >= 0.15cm](2_33)--(2_23);
\draw[->, line width=0.3mm, >=latex, shorten <= 0.2cm, shorten >= 0.15cm](1_33)--(2_23);

\draw[->, line width=0.3mm, >=latex, shorten <= 0.2cm, shorten >= 0.15cm](1_11)--(1_21);
\draw[->, line width=0.3mm, >=latex, shorten <= 0.2cm, shorten >= 0.15cm](2_11)--(1_21);
\draw[->, line width=0.3mm, >=latex, shorten <= 0.2cm, shorten >= 0.15cm](1_21)--(1_31);
\draw[->, line width=0.3mm, >=latex, shorten <= 0.2cm, shorten >= 0.15cm](1_21)--(2_31);
\draw[->, line width=0.3mm, >=latex, shorten <= 0.2cm, shorten >= 0.15cm](2_21)--(1_11);
\draw[->, line width=0.3mm, >=latex, shorten <= 0.2cm, shorten >= 0.15cm](2_21)--(2_11);
\draw[->, line width=0.3mm, >=latex, shorten <= 0.2cm, shorten >= 0.15cm](2_31)--(2_21);
\draw[->, line width=0.3mm, >=latex, shorten <= 0.2cm, shorten >= 0.15cm](1_31)--(2_21);

\draw[dashed, ->, line width=0.3mm, >=latex, shorten <= 0.2cm, shorten >= 0.15cm](1_12)--(1_13);
\draw[dashed, ->, line width=0.3mm, >=latex, shorten <= 0.2cm, shorten >= 0.15cm](2_12) to [out=75, in=210] (2_13);
\draw[dashed, ->, line width=0.3mm, >=latex, shorten <= 0.2cm, shorten >= 0.15cm](1_13) to [out=190, in=95] (2_12);
\draw[dashed, ->, line width=0.3mm, >=latex, shorten <= 0.2cm, shorten >= 0.15cm](2_13)--(1_12);

\draw[dashed, ->, line width=0.3mm, >=latex, shorten <= 0.2cm, shorten >= 0.15cm](1_22)--(1_23);
\draw[dashed, ->, line width=0.3mm, >=latex, shorten <= 0.2cm, shorten >= 0.15cm](2_22) to [out=70, in=215] (2_23);
\draw[dashed, ->, line width=0.3mm, >=latex, shorten <= 0.2cm, shorten >= 0.15cm](1_33) to [out=190, in=95] (2_32);
\draw[dashed, ->, line width=0.3mm, >=latex, shorten <= 0.2cm, shorten >= 0.15cm](2_33)--(1_32);

\draw[dashed, ->, line width=0.3mm, >=latex, shorten <= 0.2cm, shorten >= 0.15cm](1_32)--(1_33);
\draw[dashed, ->, line width=0.3mm, >=latex, shorten <= 0.2cm, shorten >= 0.15cm](2_32) to [out=75, in=210] (2_33);
\draw[dashed, ->, line width=0.3mm, >=latex, shorten <= 0.2cm, shorten >= 0.15cm](1_23) to [out=190, in=90] (2_22);
\draw[dashed, ->, line width=0.3mm, >=latex, shorten <= 0.2cm, shorten >= 0.15cm](2_23)--(1_22);

\draw[dashed, ->, line width=0.3mm, >=latex, shorten <= 0.2cm, shorten >= 0.15cm](2_11)--(2_12);
\draw[dashed, ->, line width=0.3mm, >=latex, shorten <= 0.2cm, shorten >= 0.15cm](1_11)--(1_12);
\draw[dashed, ->, line width=0.3mm, >=latex, shorten <= 0.2cm, shorten >= 0.15cm](1_12)--(2_11);
\draw[dashed, ->, line width=0.3mm, >=latex, shorten <= 0.2cm, shorten >= 0.15cm](2_12)--(1_11);

\draw[dashed, ->, line width=0.3mm, >=latex, shorten <= 0.2cm, shorten >= 0.15cm](1_11)--(2_13);
\draw[dashed, ->, line width=0.3mm, >=latex, shorten <= 0.2cm, shorten >= 0.15cm](2_11)--(1_13);
\draw[dashed, ->, line width=0.3mm, >=latex, shorten <= 0.2cm, shorten >= 0.15cm](1_13)--(1_11);
\draw[dashed, ->, line width=0.3mm, >=latex, shorten <= 0.2cm, shorten >= 0.15cm](2_13)--(2_11);

\draw[dashed, ->, line width=0.3mm, >=latex, shorten <= 0.2cm, shorten >= 0.15cm](2_21)--(2_22);
\draw[dashed, ->, line width=0.3mm, >=latex, shorten <= 0.2cm, shorten >= 0.15cm](1_21)--(1_22);
\draw[dashed, ->, line width=0.3mm, >=latex, shorten <= 0.2cm, shorten >= 0.15cm](1_22)--(2_21);
\draw[dashed, ->, line width=0.3mm, >=latex, shorten <= 0.2cm, shorten >= 0.15cm](2_22)--(1_21);

\draw[dashed, ->, line width=0.3mm, >=latex, shorten <= 0.2cm, shorten >= 0.15cm](1_23)--(1_21);
\draw[dashed, ->, line width=0.3mm, >=latex, shorten <= 0.2cm, shorten >= 0.15cm](2_23)--(2_21);
\draw[dashed, ->, line width=0.3mm, >=latex, shorten <= 0.2cm, shorten >= 0.15cm](1_21)--(2_23);
\draw[dashed, ->, line width=0.3mm, >=latex, shorten <= 0.2cm, shorten >= 0.15cm](2_21)--(1_23);

\draw[dashed, ->, line width=0.3mm, >=latex, shorten <= 0.2cm, shorten >= 0.15cm](2_31)--(2_32);
\draw[dashed, ->, line width=0.3mm, >=latex, shorten <= 0.2cm, shorten >= 0.15cm](1_31)--(1_32);
\draw[dashed, ->, line width=0.3mm, >=latex, shorten <= 0.2cm, shorten >= 0.15cm](1_32)--(2_31);
\draw[dashed, ->, line width=0.3mm, >=latex, shorten <= 0.2cm, shorten >= 0.15cm](2_32)--(1_31);

\draw[dashed, ->, line width=0.3mm, >=latex, shorten <= 0.2cm, shorten >= 0.15cm](1_31)--(2_33);
\draw[dashed, ->, line width=0.3mm, >=latex, shorten <= 0.2cm, shorten >= 0.15cm](2_31)--(1_33);
\draw[dashed, ->, line width=0.3mm, >=latex, shorten <= 0.2cm, shorten >= 0.15cm](1_33)--(1_31);
\draw[dashed, ->, line width=0.3mm, >=latex, shorten <= 0.2cm, shorten >= 0.15cm](2_33)--(2_31);
\end{tikzpicture}
{
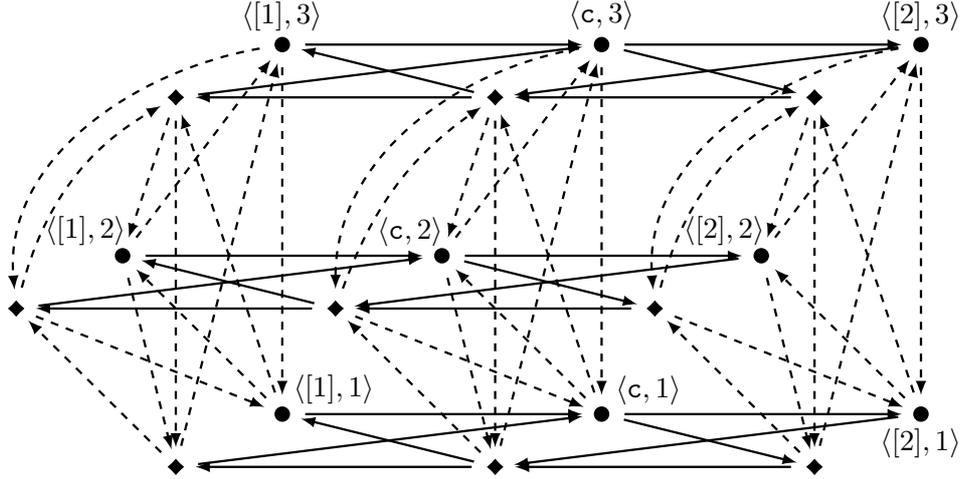
\captionof{figure}{Orientation $D$ of $(T_2\times K_3)^{(2)}$ when $deg_{T_2}(\mathtt{c})=2$, where $d(D)=3$.}\label{figA8.4.21}}
\end{center}

(b) Now, consider the case $deg_{T_2}(\mathtt{c})>2$. Suppose there exists an orientation $F$ of $(T_2\times K_\mu)^{(2)}$ with $d(F)=3$. By Lemma \ref{lemA8.3.1}, $\langle [1],1\rangle\overset{1}\twoheadrightarrow \langle \mathtt{c},1\rangle \overset{2}\twoheadleftarrow \langle [2],1\rangle$ and $\langle [1],1\rangle\overset{1}\twoheadrightarrow \langle \mathtt{c},1\rangle \overset{2}\twoheadleftarrow \langle [3],1\rangle$. However, this contradicts $\langle [3],1\rangle\overset{1}\twoheadrightarrow \langle \mathtt{c},1\rangle\overset{2}\twoheadleftarrow \langle [2],1\rangle$. Thus, $(T_2\times K_\mu)^{(2)}\in\mathscr{C}_1$.
\indent\par To show $(T_2\times K_\mu)(s_1,s_2,\ldots,s_n)\in\mathscr{C}_0\cup \mathscr{C}_1$, we need to verify $d(D)=4$. In view of (a) and its symmetry among the vertices $\langle [i],j\rangle$ for $[i]\in N_{T_2}(\mathtt{c})-\{[1]\}$ by (\ref{eqA8.4.1}), it suffices to check $d_D((p,\langle [2],j\rangle),(q,\langle [3],j\rangle))\le 4$ for $j=1,2,\ldots, \mu$, and $p,q=1,2$. That is, the partial orientation in Figure \ref{figA8.4.22} has diameter $4$, which is easy to check. Since every vertex lies in a directed $C_3$, it follows from Lemma \ref{lemA8.1.2} that $(T_2\times K_\mu)(s_1,s_2,\ldots,s_n)\in\mathscr{C}_0\cup \mathscr{C}_1$.
\qed
\begin{center}
\tikzstyle{every node}=[circle, draw, fill=black!100,
                       inner sep=0pt, minimum width=5pt]
\begin{tikzpicture}[thick,scale=0.7]%
\draw(-15,1)node[star,star points=4, label={[yshift=0.2cm, xshift=-0.3cm] 270:{}}](2_22){};
\draw(-9,1)node[star,star points=4, label={[yshift=0.4cm] 270:{}}](2_12){};
\draw(-3,1)node[star,star points=4, label={[yshift=0.5cm]270:{}}](2_32){};

\draw(-13,2)node[label={[yshift=-0.15cm, xshift=0cm]135:{\small $\langle [2],2\rangle$}}](1_22){};
\draw(-7,2)node[label={[yshift=-0.1cm, xshift=0cm]135:{\small $\langle\mathtt{c},2\rangle$}}](1_12){};
\draw(-1,2)node[label={[yshift=-0.15cm, xshift=0cm]135:{\small $\langle [3],2\rangle$}}](1_32){};

\draw(-12,5)node[star,star points=4, label={[yshift=0cm, xshift=-0.2cm] 180:{}}](2_23){};
\draw(-6,5)node[star,star points=4, label={[yshift=-0.4cm, xshift=-0.3cm] 180:{}}](2_13){};
\draw(-0,5)node[star,star points=4, label={[yshift=-0.4cm, xshift=-0.3cm] 180:{}}](2_33){};

\draw(-10,6)node[label={[yshift=-0.3cm]90:{\small $\langle [2],3\rangle$}}](1_23){};
\draw(-4,6)node[label={[yshift=-0.2cm]90:{\small $\langle\mathtt{c},3\rangle$}}](1_13){};
\draw(2,6)node[label={[yshift=-0.3cm]90:{\small $\langle [3],3\rangle$}}](1_33){};

\draw(-12,-2)node[star,star points=4, label={[yshift=0.5cm] 270:{}}](2_21){};
\draw(-6,-2)node[star,star points=4, label={[yshift=0.4cm] 270:{}}](2_11){};
\draw(0,-2)node[star,star points=4, label={[yshift=0.5cm]270:{}}](2_31){};

\draw(-10,-1)node[label={[yshift=-0.15cm, xshift=0.2cm]45:{\small $\langle [2],1\rangle$}}](1_21){};
\draw(-4,-1)node[label={[yshift=-0.1cm, xshift=0.2cm]45:{\small $\langle\mathtt{c},1\rangle$}}](1_11){};
\draw(2,-1)node[label={[yshift=0.2cm, xshift=0cm]270:{\small $\langle [3],1\rangle$}}](1_31){};

\draw[->, line width=0.3mm, >=latex, shorten <= 0.2cm, shorten >= 0.15cm](1_12)--(1_22);
\draw[->, line width=0.3mm, >=latex, shorten <= 0.2cm, shorten >= 0.15cm](1_12)--(2_22);
\draw[->, line width=0.3mm, >=latex, shorten <= 0.2cm, shorten >= 0.15cm](1_12)--(1_32);
\draw[->, line width=0.3mm, >=latex, shorten <= 0.2cm, shorten >= 0.15cm](1_12)--(2_32);
\draw[->, line width=0.3mm, >=latex, shorten <= 0.2cm, shorten >= 0.15cm](1_22)--(2_12);
\draw[->, line width=0.3mm, >=latex, shorten <= 0.2cm, shorten >= 0.15cm](2_22)--(2_12);
\draw[->, line width=0.3mm, >=latex, shorten <= 0.2cm, shorten >= 0.15cm](2_32)--(2_12);
\draw[->, line width=0.3mm, >=latex, shorten <= 0.2cm, shorten >= 0.15cm](1_32)--(2_12);

\draw[->, line width=0.3mm, >=latex, shorten <= 0.2cm, shorten >= 0.15cm](1_13)--(1_23);
\draw[->, line width=0.3mm, >=latex, shorten <= 0.2cm, shorten >= 0.15cm](1_13)--(2_23);
\draw[->, line width=0.3mm, >=latex, shorten <= 0.2cm, shorten >= 0.15cm](1_13)--(1_33);
\draw[->, line width=0.3mm, >=latex, shorten <= 0.2cm, shorten >= 0.15cm](1_13)--(2_33);
\draw[->, line width=0.3mm, >=latex, shorten <= 0.2cm, shorten >= 0.15cm](1_23)--(2_13);
\draw[->, line width=0.3mm, >=latex, shorten <= 0.2cm, shorten >= 0.15cm](2_23)--(2_13);
\draw[->, line width=0.3mm, >=latex, shorten <= 0.2cm, shorten >= 0.15cm](2_33)--(2_13);
\draw[->, line width=0.3mm, >=latex, shorten <= 0.2cm, shorten >= 0.15cm](1_33)--(2_13);

\draw[->, line width=0.3mm, >=latex, shorten <= 0.2cm, shorten >= 0.15cm](1_11)--(1_21);
\draw[->, line width=0.3mm, >=latex, shorten <= 0.2cm, shorten >= 0.15cm](1_11)--(2_21);
\draw[->, line width=0.3mm, >=latex, shorten <= 0.2cm, shorten >= 0.15cm](1_11)--(1_31);
\draw[->, line width=0.3mm, >=latex, shorten <= 0.2cm, shorten >= 0.15cm](1_11)--(2_31);
\draw[->, line width=0.3mm, >=latex, shorten <= 0.2cm, shorten >= 0.15cm](1_21)--(2_11);
\draw[->, line width=0.3mm, >=latex, shorten <= 0.2cm, shorten >= 0.15cm](2_21)--(2_11);
\draw[->, line width=0.3mm, >=latex, shorten <= 0.2cm, shorten >= 0.15cm](2_31)--(2_11);
\draw[->, line width=0.3mm, >=latex, shorten <= 0.2cm, shorten >= 0.15cm](1_31)--(2_11);

\draw[dashed, ->, line width=0.3mm, >=latex, shorten <= 0.2cm, shorten >= 0.15cm](1_22)--(1_23);
\draw[dashed, ->, line width=0.3mm, >=latex, shorten <= 0.2cm, shorten >= 0.15cm](2_22) to [out=85, in=205] (2_23);
\draw[dashed, ->, line width=0.3mm, >=latex, shorten <= 0.2cm, shorten >= 0.15cm](1_12)--(1_13);
\draw[dashed, ->, line width=0.3mm, >=latex, shorten <= 0.2cm, shorten >= 0.15cm](2_12) to [out=75, in=210] (2_13);
\draw[dashed, ->, line width=0.3mm, >=latex, shorten <= 0.2cm, shorten >= 0.15cm](1_32)--(1_33);
\draw[dashed, ->, line width=0.3mm, >=latex, shorten <= 0.2cm, shorten >= 0.15cm](2_32) to [out=75, in=210] (2_33);

\draw[dashed, ->, line width=0.3mm, >=latex, shorten <= 0.2cm, shorten >= 0.15cm](1_23) to [out=190, in=100] (2_22);
\draw[dashed, ->, line width=0.3mm, >=latex, shorten <= 0.2cm, shorten >= 0.15cm](2_23)--(1_22);
\draw[dashed, ->, line width=0.3mm, >=latex, shorten <= 0.2cm, shorten >= 0.15cm](1_33) to [out=190, in=95] (2_32);
\draw[dashed, ->, line width=0.3mm, >=latex, shorten <= 0.2cm, shorten >= 0.15cm](2_33)--(1_32);
\draw[dashed, ->, line width=0.3mm, >=latex, shorten <= 0.2cm, shorten >= 0.15cm](1_13) to [out=190, in=95] (2_12);
\draw[dashed, ->, line width=0.3mm, >=latex, shorten <= 0.2cm, shorten >= 0.15cm](2_13)--(1_12);

\draw[dashed, ->, line width=0.3mm, >=latex, shorten <= 0.2cm, shorten >= 0.15cm](2_21)--(2_22);
\draw[dashed, ->, line width=0.3mm, >=latex, shorten <= 0.2cm, shorten >= 0.15cm](1_21)--(1_22);
\draw[dashed, ->, line width=0.3mm, >=latex, shorten <= 0.2cm, shorten >= 0.15cm](1_22)--(2_21);
\draw[dashed, ->, line width=0.3mm, >=latex, shorten <= 0.2cm, shorten >= 0.15cm](2_22)--(1_21);

\draw[dashed, ->, line width=0.3mm, >=latex, shorten <= 0.2cm, shorten >= 0.15cm](1_21)--(2_23);
\draw[dashed, ->, line width=0.3mm, >=latex, shorten <= 0.2cm, shorten >= 0.15cm](2_21)--(1_23);
\draw[dashed, ->, line width=0.3mm, >=latex, shorten <= 0.2cm, shorten >= 0.15cm](1_23)--(1_21);
\draw[dashed, ->, line width=0.3mm, >=latex, shorten <= 0.2cm, shorten >= 0.15cm](2_23)--(2_21);

\draw[dashed, ->, line width=0.3mm, >=latex, shorten <= 0.2cm, shorten >= 0.15cm](2_11)--(2_12);
\draw[dashed, ->, line width=0.3mm, >=latex, shorten <= 0.2cm, shorten >= 0.15cm](1_11)--(1_12);
\draw[dashed, ->, line width=0.3mm, >=latex, shorten <= 0.2cm, shorten >= 0.15cm](1_12)--(2_11);
\draw[dashed, ->, line width=0.3mm, >=latex, shorten <= 0.2cm, shorten >= 0.15cm](2_12)--(1_11);

\draw[dashed, ->, line width=0.3mm, >=latex, shorten <= 0.2cm, shorten >= 0.15cm](1_13)--(1_11);
\draw[dashed, ->, line width=0.3mm, >=latex, shorten <= 0.2cm, shorten >= 0.15cm](2_13)--(2_11);
\draw[dashed, ->, line width=0.3mm, >=latex, shorten <= 0.2cm, shorten >= 0.15cm](1_11)--(2_13);
\draw[dashed, ->, line width=0.3mm, >=latex, shorten <= 0.2cm, shorten >= 0.15cm](2_11)--(1_13);

\draw[dashed, ->, line width=0.3mm, >=latex, shorten <= 0.2cm, shorten >= 0.15cm](2_31)--(2_32);
\draw[dashed, ->, line width=0.3mm, >=latex, shorten <= 0.2cm, shorten >= 0.15cm](1_31)--(1_32);
\draw[dashed, ->, line width=0.3mm, >=latex, shorten <= 0.2cm, shorten >= 0.15cm](1_32)--(2_31);
\draw[dashed, ->, line width=0.3mm, >=latex, shorten <= 0.2cm, shorten >= 0.15cm](2_32)--(1_31);

\draw[dashed, ->, line width=0.3mm, >=latex, shorten <= 0.2cm, shorten >= 0.15cm](1_31)--(2_33);
\draw[dashed, ->, line width=0.3mm, >=latex, shorten <= 0.2cm, shorten >= 0.15cm](2_31)--(1_33);
\draw[dashed, ->, line width=0.3mm, >=latex, shorten <= 0.2cm, shorten >= 0.15cm](1_33)--(1_31);
\draw[dashed, ->, line width=0.3mm, >=latex, shorten <= 0.2cm, shorten >= 0.15cm](2_33)--(2_31);
\end{tikzpicture}
{\captionof{figure}{Partial orientation $D$ of $(T_2\times K_3)^{(2)}$ when $deg_{T_2}(\mathtt{c})>2$, where $d(D)=4$.}\label{figA8.4.22}}
\end{center}
\indent\par In Proposition \ref{ppnA8.4.2}, we generalise the sufficient condition in Proposition \ref{ppnA8.4.1}(b), ``$deg_{T_2}(\mathtt{c})$ $>2$", for the vertex-multiplication of $T_2\times K_\mu$ to be in $\mathscr{C}_1$.

\begin{ppn}\label{ppnA8.4.2}
Let $\mu \ge 3$ and $m=\min\{s_{\langle\mathtt{c},v\rangle}\mid v\in V(K_\mu)\}$. If $deg_{T_2}(\mathtt{c})>{{m}\choose{\lfloor{m/2}\rfloor}}$, then $(T_2\times K_\mu)(s_1,s_2,\ldots,s_n)\in\mathscr{C}_1$.
\end{ppn}
\noindent\textit{Proof}: Suppose $D$ is an orientation of $(T_2\times K_\mu)(s_1,s_2,\ldots,s_n)$ with $d(D)=3=d(T_2\times K_\mu)$. In view of parity, $d_D((p,\langle [i],v\rangle), (q,\langle [j],v\rangle))=2$ for any $p=1,2,\dots,s_{\langle [i],v\rangle}$, $q=1,2,\dots,s_{\langle [j],v\rangle}$ and all $[i],[j]\in N_{T_2}(\mathtt{c})$ with $[i]\neq [j]$. For any $(p,\langle [i],v\rangle)\in V(D)$, let $O^{\langle \mathtt{c},v\rangle}((p,\langle [i],v\rangle))=O((p,\langle [i],v\rangle))\cap \{(r,\langle \mathtt{c}, v\rangle)\mid r=1,2,\ldots, s_{\langle \mathtt{c}, v\rangle)}\}$.
\indent\par Since $deg_{T_2}(\mathtt{c})>{{m}\choose{\lfloor{m/2}\rfloor}}$, there exists some $v^*\in V(K_\mu)$ such that $deg_{T_2}(\mathtt{c})>{{s_{\langle\mathtt{c},v^*\rangle}}\choose{\lfloor{s_{\langle\mathtt{c},v^*\rangle}/2}\rfloor}}$. By Sperner Theorem, for some $p^*=1,2,\dots,s_{\langle [i^*],v^*\rangle}$, some $q^*=1,2,\dots,s_{\langle [j^*],v^*\rangle}$ and some $[i^*],[j^*]\in N_{T_2}(\mathtt{c})$ with $[i^*]\neq [j^*]$,
$O^{\langle \mathtt{c},v^*\rangle}((p^*,\langle [i^*],v^*\rangle))\subseteq O^{\langle \mathtt{c},v^*\rangle}((q^*,\langle [j^*],v^*\rangle))$. Hence, it follows that $d_D((p^*,\langle [i^*],v^*\rangle), (q^*,\langle [j^*],v^*\rangle))>2$, a contradiction. Hence,  $(T_2\times K_\mu)(s_1,s_2,\ldots,s_n)\not\in\mathscr{C}_0$. By Proposition \ref{ppnA8.4.1}(b), $(T_2\times K_\mu)(s_1,s_2,\ldots,s_n)\in\mathscr{C}_1$.
\qed

\begin{rmk}
The same proof and notation as Proposition \ref{ppnA8.4.2} shows that if $deg_{T_2}(\mathtt{c})>{{m}\choose{\lfloor{m/2}\rfloor}}$, then $(T_2\times K_2)(s_1,s_2,\ldots,s_n)\not\in\mathscr{C}_0$.
\end{rmk}
~\\
\begin{ppn}\label{ppnA8.4.4}
For $\mu\ge 3$, $(P_2\times K_\mu)(s_1,s_2,\ldots,s_n)\in\mathscr{C}_1$.
\end{ppn}
\noindent\textit{Proof}: Suppose $F$ is an orientation of $(P_2\times K_\mu)(s_1,s_2,\ldots,s_n)$ with $d(F)=2=d(P_2\times K_\mu)$. It follows from $d_F((p,\langle u,x\rangle),(q,\langle v,x\rangle))\le 2$ that $(p,\langle u,x\rangle)\rightarrow (q,\langle v,x\rangle)$ for $u,v\in V(P_2)$, $x\in V(K_\mu)$, $p=1,2,\ldots, s_{\langle u,x\rangle}$, $q=1,2,\ldots, s_{\langle v,x\rangle}$. Then, $d_F((q,\langle v,x\rangle),(p,\langle u,x\rangle))$ $>2$, a contradiction. Hence, $(P_2\times K_\mu)(s_1,s_2,\ldots,s_n)\not\in\mathscr{C}_0$.
\indent\par Define an orientation $D$ of $(P_2\times K_\mu)^{(2)}$ as follows. (See Figure \ref{figA8.4.23} when $\mu=3$.)
\begin{align*}
&\langle 2,1\rangle \rightrightarrows \langle 1,1\rangle, \langle 1,2\rangle \rightrightarrows \langle 2,2\rangle, \langle 1,i\rangle\rightsquigarrow \langle 2,i\rangle \text{ for }i=3,4,\ldots, \mu.\\
&\langle k,j_1\rangle\rightsquigarrow \langle k,j_2\rangle \text{ whenever } 2\le j_1<j_2\le\mu,\text{ and}\\
&\langle k,j\rangle\rightsquigarrow \langle k,1\rangle \rightsquigarrow \langle k,2\rangle\text{ for }j=3,4,\ldots,\mu,
\end{align*}
for $k=1,2$.
\indent\par It can be verified easily that $d(D)=3$. Hence, $(P_2\times K_\mu)^{(2)}\in \mathscr{C}_1$. Furthermore, since every vertex lies in a directed $C_3$, it follows from Lemma \ref{lemA8.1.2} that $(P_2\times K_\mu)(s_1,s_2,\ldots,s_n)\in\mathscr{C}_1$.
\qed
\begin{center}
\tikzstyle{every node}=[circle, draw, fill=black!100,
                       inner sep=0pt, minimum width=5pt]
\begin{tikzpicture}[thick,scale=0.7]%
\draw(-15,1)node[star,star points=4, label={[yshift=0.4cm, xshift=-0.5cm] 270:{}}](2_12){};
\draw(-7,1)node[star,star points=4, label={[yshift=0.4cm, xshift=-0.5cm] 270:{}}](2_22){};

\draw(-13,2)node[label={[yshift=-0.15cm, xshift=0cm]135:{\small $\langle 1,2\rangle$}}](1_12){};
\draw(-5,2)node[label={[yshift=-0.15cm, xshift=0cm]135:{\small $\langle 2,2\rangle$}}](1_22){};

\draw(-12,5)node[star,star points=4, label={[yshift=-0.6cm, xshift=0.3cm] 90:{}}](2_13){};
\draw(-4,5)node[star,star points=4, label={[yshift=-0.6cm, xshift=0.3cm] 90:{}}](2_23){};

\draw(-10,6)node[label={[yshift=-0.2cm]90:{\small $\langle 1,3\rangle$}}](1_13){};
\draw(-2,6)node[label={[yshift=-0.2cm]90:{\small $\langle 2,3\rangle$}}](1_23){};

\draw(-12,-2)node[star,star points=4, label={[yshift=0.4cm] 270:{}}](2_11){};
\draw(-4,-2)node[star,star points=4, label={[yshift=0.4cm] 270:{}}](2_21){};

\draw(-10,-1)node[label={[yshift=-0.1cm, xshift=0.1cm]45:{\small $\langle 1,1\rangle$}}](1_11){};
\draw(-2,-1)node[label={[yshift=0.15cm, xshift=0cm]270:{\small $\langle 2,1\rangle$}}](1_21){};

\draw[dashed, ->, line width=0.3mm, >=latex, shorten <= 0.2cm, shorten >= 0.15cm](1_13)--(1_23);
\draw[dashed, ->, line width=0.3mm, >=latex, shorten <= 0.2cm, shorten >= 0.15cm](2_13)--(2_23);
\draw[dashed, ->, line width=0.3mm, >=latex, shorten <= 0.2cm, shorten >= 0.15cm](2_23)--(1_13);
\draw[dashed, ->, line width=0.3mm, >=latex, shorten <= 0.2cm, shorten >= 0.15cm](1_23)--(2_13);

\draw[dashed, ->, line width=0.3mm, >=latex, shorten <= 0.2cm, shorten >= 0.15cm](1_22)--(1_23);
\draw[dashed, ->, line width=0.3mm, >=latex, shorten <= 0.2cm, shorten >= 0.15cm](2_22) to [out=85, in=205] (2_23);
\draw[dashed, ->, line width=0.3mm, >=latex, shorten <= 0.2cm, shorten >= 0.15cm](1_23) to [out=190, in=100] (2_22);
\draw[dashed, ->, line width=0.3mm, >=latex, shorten <= 0.2cm, shorten >= 0.15cm](2_23)--(1_22);

\draw[dashed, ->, line width=0.3mm, >=latex, shorten <= 0.2cm, shorten >= 0.15cm](1_12)--(1_13);
\draw[dashed, ->, line width=0.3mm, >=latex, shorten <= 0.2cm, shorten >= 0.15cm](2_12) to [out=75, in=210] (2_13);
\draw[dashed, ->, line width=0.3mm, >=latex, shorten <= 0.2cm, shorten >= 0.15cm](1_13) to [out=190, in=95] (2_12);
\draw[dashed, ->, line width=0.3mm, >=latex, shorten <= 0.2cm, shorten >= 0.15cm](2_13)--(1_12);

\draw[dashed, ->, line width=0.3mm, >=latex, shorten <= 0.2cm, shorten >= 0.15cm](1_23)--(1_21);
\draw[dashed, ->, line width=0.3mm, >=latex, shorten <= 0.2cm, shorten >= 0.15cm](2_23)--(2_21);
\draw[dashed, ->, line width=0.3mm, >=latex, shorten <= 0.2cm, shorten >= 0.15cm](1_21)--(2_23);
\draw[dashed, ->, line width=0.3mm, >=latex, shorten <= 0.2cm, shorten >= 0.15cm](2_21)--(1_23);

\draw[dashed, ->, line width=0.3mm, >=latex, shorten <= 0.2cm, shorten >= 0.15cm](1_11)--(1_12);
\draw[dashed, ->, line width=0.3mm, >=latex, shorten <= 0.2cm, shorten >= 0.15cm](1_12)--(2_11);
\draw[dashed, ->, line width=0.3mm, >=latex, shorten <= 0.2cm, shorten >= 0.15cm](2_12)--(1_11);
\draw[dashed, ->, line width=0.3mm, >=latex, shorten <= 0.2cm, shorten >= 0.15cm](2_11)--(2_12);

\draw[dashed, ->, line width=0.3mm, >=latex, shorten <= 0.2cm, shorten >= 0.15cm](1_11)--(2_13);
\draw[dashed, ->, line width=0.3mm, >=latex, shorten <= 0.2cm, shorten >= 0.15cm](2_11)--(1_13);
\draw[dashed, ->, line width=0.3mm, >=latex, shorten <= 0.2cm, shorten >= 0.15cm](1_13)--(1_11);
\draw[dashed, ->, line width=0.3mm, >=latex, shorten <= 0.2cm, shorten >= 0.15cm](2_13)--(2_11);

\draw[dashed, ->, line width=0.3mm, >=latex, shorten <= 0.2cm, shorten >= 0.15cm](2_21)--(2_22);
\draw[dashed, ->, line width=0.3mm, >=latex, shorten <= 0.2cm, shorten >= 0.15cm](1_21)--(1_22);
\draw[dashed, ->, line width=0.3mm, >=latex, shorten <= 0.2cm, shorten >= 0.15cm](1_22)--(2_21);
\draw[dashed, ->, line width=0.3mm, >=latex, shorten <= 0.2cm, shorten >= 0.15cm](2_22)--(1_21);

\draw[densely dotted, ->, line width=0.3mm, >=latex, shorten <= 0.2cm, shorten >= 0.15cm](1_12)--(1_22);
\draw[densely dotted, ->, line width=0.3mm, >=latex, shorten <= 0.2cm, shorten >= 0.15cm](2_12)--(2_22);
\draw[densely dotted, ->, line width=0.3mm, >=latex, shorten <= 0.2cm, shorten >= 0.15cm](1_12)--(2_22);
\draw[densely dotted, ->, line width=0.3mm, >=latex, shorten <= 0.2cm, shorten >= 0.15cm](2_12)--(1_22);

\draw[densely dotted, ->, line width=0.3mm, >=latex, shorten <= 0.2cm, shorten >= 0.15cm](1_21)--(1_11);
\draw[densely dotted, ->, line width=0.3mm, >=latex, shorten <= 0.2cm, shorten >= 0.15cm](2_21)--(2_11);
\draw[densely dotted, ->, line width=0.3mm, >=latex, shorten <= 0.2cm, shorten >= 0.15cm](1_21)--(2_11);
\draw[densely dotted, ->, line width=0.3mm, >=latex, shorten <= 0.2cm, shorten >= 0.15cm](2_21)--(1_11);
\end{tikzpicture}
{
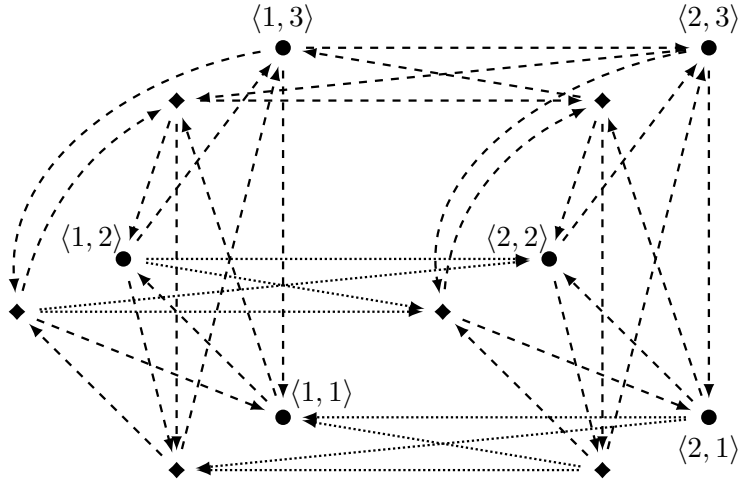
\captionof{figure}{Orientation $D$ of $(P_2\times K_3)^{(2)}$, where $d(D)=3$.}\label{figA8.4.23}}
\end{center}

\section{Cartesian product of two cycles $C_\lambda\times C_\mu$}
In this section, we prove Theorem \ref{thmA8.1.11}.
\begin{ppn}\label{ppnA8.5.1}
If $\lambda\ge\mu\ge 4$, then $(C_\lambda\times C_\mu)(s_1,s_2,\ldots,s_n)\in\mathscr{C}_0$.
\end{ppn}
\noindent\textit{Proof}: We shall use a similar strategy as in Theorem \ref{thmA8.1.10}. Partition $V(C_\lambda)$ into $V_1=\{v \mid v$ is odd$\}$ and $V_2=\{v \mid v$ is even$\}$. Let $F$ be a strong orientation of $C_\mu$, say $1\rightarrow 2\rightarrow \cdots \rightarrow \mu\rightarrow 1$, and define an orientation $D$ for $(C_\lambda\times C_\mu)^{(2)}$ as follows.
\begin{align*}
\langle u,x\rangle\rightrightarrows \langle u,y\rangle\iff x\rightarrow y \text{ in }F
\end{align*}
for any $u\in V_1$, and any $x,y\in V(C_\mu)$, i.e. the copy $C_\mu^{(2)}$ is oriented similarly to $F$.
\begin{align*}
\langle u,x\rangle\rightrightarrows \langle u,y\rangle\iff y\rightarrow x \text{ in }F
\end{align*}
for any $u\in V_2$, and any $x,y\in V(C_\mu)$, i.e. the copy $C_\mu^{(2)}$ is oriented similarly to $\tilde{F}$.
\begin{align*}
\langle u,x\rangle\rightsquigarrow \langle u+1,x\rangle\text{ (addition is taken modulo }\lambda\text{)}
\end{align*}
for any $u\in V(C_\lambda)$ and any $x\in V(C_\mu)$.
\indent\par We claim that $d_D((p,\langle u,x\rangle),(q,\langle v,y\rangle))\le \lfloor\frac{\lambda}{2}\rfloor+\lfloor\frac{\mu}{2}\rfloor=d(C_\lambda \times C_\mu)$ for any $\langle u,x\rangle,\langle v,y\rangle$ $\in V(C_\lambda\times C_\mu)$, and $p,q=1,2$. Suppose $u=v\in V_1$. Note that either $d_F(x,y)\le \lfloor\frac{\mu}{2}\rfloor$ or $d_{\tilde{F}}(x,y)\le \lfloor\frac{\mu}{2}\rfloor$. So, there exist paths $P$ and $P'$ in $D$, each of length at most $\lfloor\frac{\mu}{2}\rfloor$, from $\{(1,\langle u,x\rangle),(2,\langle u,x\rangle)\}$ to $\{(1,\langle v,y\rangle),(2,\langle v,y\rangle)\}$ and from $\{(1,\langle w,x\rangle),(2,\langle w,x\rangle)\}$ to $\{(1,\langle w,y\rangle),$ $(2,\langle w,y\rangle)\}$, where $w\in V_2$ is some vertex adjacent to $u$ in $C_\lambda$ respectively. In the former case, $P$ suffices and we are done. In the latter case, we shall further assume $w=u+1\pmod{\lambda}$ for simplicity; the proof is similar if $u=w+1\pmod{\lambda}$. Then, $(p,\langle u,x\rangle) (p,\langle w,x\rangle)$ with $P'$ and $(3-q,\langle w,y\rangle) (q,\langle u,y\rangle)$ form a $(p,\langle u,x\rangle)-(q,\langle v,y\rangle)$ path of length at most $2+\lfloor\frac{\mu}{2}\rfloor\le \lfloor\frac{\lambda}{2}\rfloor+\lfloor\frac{\mu}{2}\rfloor$. A similar proof follows if $u=v\in V_2$.
\indent\par Suppose $u\neq v$. For simplicity, we shall assume $u (u+1) \ldots (u+l) v$ to be a shortest $u-v$ path in $C_\lambda$; the proof is similar if the shortest path is $u (u-1)\ldots (u-l) v$. If $x=y$, then $(p,\langle u,x\rangle) (p,\langle u+1,x\rangle) \ldots (p,\langle v,x\rangle) (3-p,\langle u+l,x\rangle) (3-p,\langle v,x\rangle)$ guarantees a $(p,\langle u,x\rangle)-(q,\langle v,y\rangle)$ path of length at most $\lfloor\frac{\lambda}{2}\rfloor+2\le \lfloor\frac{\lambda}{2}\rfloor+\lfloor\frac{\mu}{2}\rfloor$.
\indent\par Next, suppose $x\neq y$. Furthermore, we shall assume $v\in V_1$; the proof is similar if $v\in V_2$. Again, consider the cases 
$d_F(x,y)\le \lfloor\frac{\mu}{2}\rfloor$ or $d_{\tilde{F}}(x,y)\le \lfloor\frac{\mu}{2}\rfloor$. In the former case, there is a path $Q$ of length at most $\lfloor\frac{\mu}{2}\rfloor$ from $\{(1,\langle v,x\rangle),(2,\langle v,x\rangle)\}$ to $\{(1,\langle v,y\rangle),(2,\langle v,y\rangle)\}$ in $D$. So, $(p,\langle u,x\rangle) (p,\langle u+1,x\rangle) \ldots (p,\langle u+l,x\rangle) (p,\langle v,x\rangle)$ and $Q$ form a $(p,\langle u,x\rangle)-(q,\langle v,y\rangle)$ path of length at most $\lfloor\frac{\lambda}{2}\rfloor+\lfloor\frac{\mu}{2}\rfloor$. In the latter case, unless $u=\lambda$ is odd and $v=1$, there exists some $i\in\{0,1\}$ such that $u+i\in V_2$. Moreover, there is a path $Q'$ of length at most $\lfloor\frac{\mu}{2}\rfloor$ from $\{(1,\langle u+i,x\rangle),(2,\langle u+i,x\rangle)\}$ to $\{(1,\langle u+i,y\rangle),(2,\langle u+i,y\rangle)\}$ in $D$ so that $(p,\langle u,x\rangle) (p,\langle u+1,x\rangle) \ldots (p,\langle u+i,x\rangle)$ with $Q'$ and $(q,\langle u+i,y\rangle) (q,\langle u+i+1,y\rangle) \ldots (q,\langle v,y\rangle)$ form a $(p,\langle u,x\rangle)-(q,\langle v,y\rangle)$ of length at most $\lfloor\frac{\lambda}{2}\rfloor+\lfloor\frac{\mu}{2}\rfloor$. 
\indent\par Finally, if $u=\lambda$ is odd, $v=1$, and $y-x\le \lfloor\frac{\mu}{2}\rfloor +1\pmod{\mu}$, then $(p,\langle \lambda,x\rangle)\rightarrow (p,\langle 1,x\rangle)\rightarrow \{(1,\langle 1,x+1\rangle),(2,\langle 1,x+1\rangle)\} \rightarrow \{(1,\langle 1,x+2\rangle),(2,\langle 1,x+2\rangle)\}\rightarrow \cdots \rightarrow \{(1,\langle 1,y\rangle),(2,\langle 1,y\rangle)\}$ ensures a path of length at most $1+\lfloor\frac{\mu}{2}\rfloor+1\le \lfloor\frac{\lambda}{2}\rfloor+\lfloor\frac{\mu}{2}\rfloor$.
If $u=\lambda$ is odd, $v=1$, and $y-x> \lfloor\frac{\mu}{2}\rfloor +1\pmod{\mu}$, then $(p,\langle \lambda,x\rangle)\rightarrow (3-p,\langle \lambda-1,x\rangle)\rightarrow \{(1,\langle \lambda-1,x-1\rangle),(2,\langle \lambda-1,x-1\rangle)\}\rightarrow \cdots \rightarrow \{(1,\langle \lambda-1,y\rangle),(2,\langle \lambda-1,y\rangle)\}$ and $(q,\langle \lambda-1,y\rangle) (q,\langle \lambda,y\rangle) (q,\langle 1,y\rangle)$ form a $(p,\langle \lambda,x\rangle)-(q,\langle 1,y\rangle)$ path of length at most $3+\lceil\frac{\mu}{2}\rceil-2\le \lfloor\frac{\lambda}{2}\rfloor+\lfloor\frac{\mu}{2}\rfloor$.
\indent\par Since every vertex lies in a directed $C_4$, it follows from Lemma \ref{lemA8.1.2} that $(C_\lambda\times C_\mu)(s_1,s_2,\ldots, s_n) \in \mathscr{C}_0$. 
\qed

\begin{cor}
$(C_3\times C_3)(s_1,s_2,\ldots,s_n)\in\mathscr{C}_0\cup \mathscr{C}_1$.
\end{cor}
\noindent\textit{Proof}: We claim that $d(D)=3=d(C_3\times C_3)+1$ where $D$ is as defined in Proposition \ref{ppnA8.5.1}. For any $\langle u,x\rangle,\langle v,y\rangle\in V(C_\lambda\times C_\mu)$, observe that either $\langle u,x\rangle\rightsquigarrow \langle v,x\rangle$ or $\langle v,x\rangle\rightsquigarrow \langle u,x\rangle$ and $\langle v,x\rangle\rightrightarrows \langle v,x+1\rangle\rightrightarrows \langle v,x+2\rangle$ or $\langle v,x\rangle\rightrightarrows \langle v,x-1\rangle\rightrightarrows \langle v,x-2\rangle$, where the addition is taken modulo 3, proves the claim. Hence, $(C_3\times C_3)^{(2)}\in\mathscr{C}_0\cup\mathscr{C}_1$. Since every vertex lies in a directed $C_3$, it follows from Lemma \ref{lemA8.1.2} that $(C_3\times C_3)(s_1,s_2,\ldots, s_n)\in\mathscr{C}_0\cup \mathscr{C}_1$. 
\qed
\begin{ppn}
$(C_4\times C_3)^{(2)}\in\mathscr{C}_0$ and $(C_4\times C_3)(s_1,s_2,\ldots,s_n) \in\mathscr{C}_0\cup \mathscr{C}_1$.
\end{ppn}
\noindent\textit{Proof}: Define an orientation $D$ for $(C_4\times C_3)^{(2)}$ as follows. (See Figure \ref{figA8.5.24}.)
\begin{align*}
&\langle 2,i\rangle\rightsquigarrow \langle 1,i\rangle\text{ and }\langle 3,i\rangle\rightsquigarrow \langle 4,i\rangle\text{ for }i=1,2,3.\\
&\langle 1,2\rangle\rightrightarrows \langle 1,1\rangle\rightrightarrows \langle 4,1\rangle\rightrightarrows \langle 4,2\rangle\rightrightarrows \langle 1,2\rangle\rightrightarrows \langle 1,3\rangle\rightrightarrows \langle 4,3\rangle\rightrightarrows \langle 4,2\rangle,\\
&\langle 3,2\rangle\rightrightarrows \langle 3,1\rangle\rightrightarrows \langle 2,1\rangle\rightrightarrows \langle 2,2\rangle\rightrightarrows \langle 3,2\rangle\rightrightarrows \langle 3,3\rangle\rightrightarrows \langle 2,3\rangle\rightrightarrows \langle 2,2\rangle,\\
&\langle 1,3\rangle\rightrightarrows \langle 1,1\rangle, \langle 4,1\rangle\rightrightarrows \langle 4,3\rangle, \langle 3,3\rangle\rightrightarrows \langle 3,1\rangle,\text{ and } \langle 2,1\rangle\rightrightarrows \langle 2,3\rangle.
\end{align*}
It is easy to check $d(D)=3=d(C_4\times C_3)$. Since every vertex lies in a directed $C_4$, it follows from Lemma \ref{lemA8.1.2} that $(C_4\times C_3)(s_1,s_2,\ldots,s_n)\in\mathscr{C}_0\cup \mathscr{C}_1$.
\qed
\begin{center}
\tikzstyle{every node}=[circle, draw, fill=black!100,
                       inner sep=0pt, minimum width=5pt]
\begin{tikzpicture}[thick,scale=0.7]%
\draw(-11,6)node[label={[yshift=-0.2cm]90:{\small $\langle 1,1\rangle$}}](1_11){};
\draw(-12,5)node[star,star points=4, label={[yshift=-0.4cm, xshift=-0.1cm] 180:{}}](2_11){};
\draw(-7,-1)node[label={[yshift=-0.15cm, xshift=0cm]45:{\small $\langle 2,1\rangle$}}](1_21){};
\draw(-8,-2)node[star,star points=4, label={[yshift=0.5cm, xshift=0.9cm] 270:{}}](2_21){};
\draw(-7,-6)node[label={[yshift=-0.1cm, xshift=0.1cm]45:{\small $\langle 3,1\rangle$}}](1_31){};
\draw(-8,-7)node[star,star points=4, label={[yshift=0.5cm, xshift=0cm] 270:{}}](2_31){};
\draw(-11,1)node[label={[yshift=-0.1cm, xshift=-0.2cm]135:{\small $\langle 4,1\rangle$}}](1_41){};
\draw(-12,0)node[star,star points=4, label={[yshift=0.6cm, xshift=-0.7cm] 270:{}}](2_41){};

\draw(-5,6)node[label={[yshift=-0.2cm]90:{\small $\langle 1,2\rangle$}}](1_12){};
\draw(-6,5)node[star,star points=4, label={[yshift=-0.4cm, xshift=-0.1cm] 180:{}}](2_12){};
\draw(-1,-1)node[label={[yshift=-0.15cm, xshift=0cm]45:{\small $\langle 2,2\rangle$}}](1_22){};
\draw(-2,-2)node[star,star points=4, label={[yshift=0.5cm, xshift=0.9cm] 270:{}}](2_22){};
\draw(-1,-6)node[label={[yshift=-0.1cm, xshift=0.1cm]45:{\small $\langle 3,2\rangle$}}](1_32){};
\draw(-2,-7)node[star,star points=4, label={[yshift=0.5cm, xshift=0cm] 270:{}}](2_32){};
\draw(-5,1)node[label={[yshift=-0.1cm, xshift=-0.2cm]135:{\small $\langle 4,2\rangle$}}](1_42){};
\draw(-6,0)node[star,star points=4, label={[yshift=0.6cm, xshift=0.9cm] 270:{}}](2_42){};

\draw(1,6)node[label={[yshift=-0.2cm]90:{\small $\langle 1,3\rangle$}}](1_13){};
\draw(0,5)node[star,star points=4, label={[yshift=-0.4cm, xshift=-0.1cm] 180:{}}](2_13){};
\draw(5,-1)node[label={[yshift=-0.15cm, xshift=0cm]45:{\small $\langle 2,3\rangle$}}](1_23){};
\draw(4,-2)node[star,star points=4, label={[yshift=0.5cm, xshift=0.9cm] 270:{}}](2_23){};
\draw(5,-6)node[label={[yshift=-0.1cm, xshift=0.1cm]45:{\small $\langle 3,3\rangle$}}](1_33){};
\draw(4,-7)node[star,star points=4, label={[yshift=0.5cm, xshift=0cm] 270:{}}](2_33){};
\draw(1,1)node[label={[yshift=-0.1cm, xshift=-0.2cm]135:{\small $\langle 4,3\rangle$}}](1_43){};
\draw(0,0)node[star,star points=4, label={[yshift=0.6cm, xshift=0.9cm] 270:{}}](2_43){};
\draw[densely dotted, ->, line width=0.3mm, >=latex, shorten <= 0.2cm, shorten >= 0.15cm](1_11)--(1_41);
\draw[densely dotted, ->, line width=0.3mm, >=latex, shorten <= 0.2cm, shorten >= 0.15cm](1_11) to [out=200, in=140] (2_41);
\draw[densely dotted, ->, line width=0.3mm, >=latex, shorten <= 0.2cm, shorten >= 0.15cm](2_11)--(1_41);
\draw[densely dotted, ->, line width=0.3mm, >=latex, shorten <= 0.2cm, shorten >= 0.15cm](2_11) to [out=250, in=120](2_41);
\draw[densely dotted, ->, line width=0.3mm, >=latex, shorten <= 0.2cm, shorten >= 0.15cm](1_31)--(1_21);
\draw[densely dotted, ->, line width=0.3mm, >=latex, shorten <= 0.2cm, shorten >= 0.15cm](1_31)--(2_21);
\draw[densely dotted, ->, line width=0.3mm, >=latex, shorten <= 0.2cm, shorten >= 0.15cm](2_31)--(1_21);
\draw[densely dotted, ->, line width=0.3mm, >=latex, shorten <= 0.2cm, shorten >= 0.15cm](2_31)--(2_21);
\draw[dashed, ->, line width=0.3mm, >=latex, shorten <= 0.2cm, shorten >= 0.15cm](1_11)--(2_21);
\draw[dashed, ->, line width=0.3mm, >=latex, shorten <= 0.2cm, shorten >= 0.15cm](2_11)--(1_21);
\draw[dashed, ->, line width=0.3mm, >=latex, shorten <= 0.2cm, shorten >= 0.15cm](1_21)--(1_11);
\draw[dashed, ->, line width=0.3mm, >=latex, shorten <= 0.2cm, shorten >= 0.15cm](2_21)--(2_11);
\draw[dashed, ->, line width=0.3mm, >=latex, shorten <= 0.2cm, shorten >= 0.15cm](1_41)--(2_31);
\draw[dashed, ->, line width=0.3mm, >=latex, shorten <= 0.2cm, shorten >= 0.15cm](2_41)--(1_31);
\draw[dashed, ->, line width=0.3mm, >=latex, shorten <= 0.2cm, shorten >= 0.15cm](1_31)--(1_41);
\draw[dashed, ->, line width=0.3mm, >=latex, shorten <= 0.2cm, shorten >= 0.15cm](2_31)--(2_41);
\draw[densely dotted, ->, line width=0.3mm, >=latex, shorten <= 0.2cm, shorten >= 0.15cm](1_42)--(1_12);
\draw[densely dotted, ->, line width=0.3mm, >=latex, shorten <= 0.2cm, shorten >= 0.15cm](1_42)--(2_12);
\draw[densely dotted, ->, line width=0.3mm, >=latex, shorten <= 0.2cm, shorten >= 0.15cm](2_42) to [out=140, in=200](1_12);
\draw[densely dotted, ->, line width=0.3mm, >=latex, shorten <= 0.2cm, shorten >= 0.15cm](2_42) to [out=120, in=250](2_12);
\draw[densely dotted, ->, line width=0.3mm, >=latex, shorten <= 0.2cm, shorten >= 0.15cm](1_22)--(1_32);
\draw[densely dotted, ->, line width=0.3mm, >=latex, shorten <= 0.2cm, shorten >= 0.15cm](1_22)--(2_32);
\draw[densely dotted, ->, line width=0.3mm, >=latex, shorten <= 0.2cm, shorten >= 0.15cm](2_22)--(1_32);
\draw[densely dotted, ->, line width=0.3mm, >=latex, shorten <= 0.2cm, shorten >= 0.15cm](2_22)--(2_32);
\draw[dashed, ->, line width=0.3mm, >=latex, shorten <= 0.2cm, shorten >= 0.15cm](1_12)--(2_22);
\draw[dashed, ->, line width=0.3mm, >=latex, shorten <= 0.2cm, shorten >= 0.15cm](2_12)--(1_22);
\draw[dashed, ->, line width=0.3mm, >=latex, shorten <= 0.2cm, shorten >= 0.15cm](1_22)--(1_12);
\draw[dashed, ->, line width=0.3mm, >=latex, shorten <= 0.2cm, shorten >= 0.15cm](2_22)--(2_12);
\draw[dashed, ->, line width=0.3mm, >=latex, shorten <= 0.2cm, shorten >= 0.15cm](1_42)--(2_32);
\draw[dashed, ->, line width=0.3mm, >=latex, shorten <= 0.2cm, shorten >= 0.15cm](2_42)--(1_32);
\draw[dashed, ->, line width=0.3mm, >=latex, shorten <= 0.2cm, shorten >= 0.15cm](1_32)--(1_42);
\draw[dashed, ->, line width=0.3mm, >=latex, shorten <= 0.2cm, shorten >= 0.15cm](2_32)--(2_42);
\draw[densely dotted, ->, line width=0.3mm, >=latex, shorten <= 0.2cm, shorten >= 0.15cm](1_13)--(1_43);
\draw[densely dotted, ->, line width=0.3mm, >=latex, shorten <= 0.2cm, shorten >= 0.15cm](1_13) to [out=200, in=140](2_43);
\draw[densely dotted, ->, line width=0.3mm, >=latex, shorten <= 0.2cm, shorten >= 0.15cm](2_13)--(1_43);
\draw[densely dotted, ->, line width=0.3mm, >=latex, shorten <= 0.2cm, shorten >= 0.15cm](2_13) to [out=250, in=120](2_43);
\draw[densely dotted, ->, line width=0.3mm, >=latex, shorten <= 0.2cm, shorten >= 0.15cm](1_33)--(1_23);
\draw[densely dotted, ->, line width=0.3mm, >=latex, shorten <= 0.2cm, shorten >= 0.15cm](1_33)--(2_23);
\draw[densely dotted, ->, line width=0.3mm, >=latex, shorten <= 0.2cm, shorten >= 0.15cm](2_33)--(1_23);
\draw[densely dotted, ->, line width=0.3mm, >=latex, shorten <= 0.2cm, shorten >= 0.15cm](2_33)--(2_23);
\draw[dashed, ->, line width=0.3mm, >=latex, shorten <= 0.2cm, shorten >= 0.15cm](1_13)--(2_23);
\draw[dashed, ->, line width=0.3mm, >=latex, shorten <= 0.2cm, shorten >= 0.15cm](2_13)--(1_23);
\draw[dashed, ->, line width=0.3mm, >=latex, shorten <= 0.2cm, shorten >= 0.15cm](1_23)--(1_13);
\draw[dashed, ->, line width=0.3mm, >=latex, shorten <= 0.2cm, shorten >= 0.15cm](2_23)--(2_13);
\draw[dashed, ->, line width=0.3mm, >=latex, shorten <= 0.2cm, shorten >= 0.15cm](1_43)--(2_33);
\draw[dashed, ->, line width=0.3mm, >=latex, shorten <= 0.2cm, shorten >= 0.15cm](2_43)--(1_33);
\draw[dashed, ->, line width=0.3mm, >=latex, shorten <= 0.2cm, shorten >= 0.15cm](1_33)--(1_43);
\draw[dashed, ->, line width=0.3mm, >=latex, shorten <= 0.2cm, shorten >= 0.15cm](2_33)--(2_43);
\draw[densely dotted, ->, line width=0.3mm, >=latex, shorten <= 0.2cm, shorten >= 0.15cm](1_12)--(1_13);
\draw[densely dotted, ->, line width=0.3mm, >=latex, shorten <= 0.2cm, shorten >= 0.15cm](1_12)--(2_13);
\draw[densely dotted, ->, line width=0.3mm, >=latex, shorten <= 0.2cm, shorten >= 0.15cm](2_12)--(1_13);
\draw[densely dotted, ->, line width=0.3mm, >=latex, shorten <= 0.2cm, shorten >= 0.15cm](2_12)--(2_13);

\draw[densely dotted, ->, line width=0.3mm, >=latex, shorten <= 0.2cm, shorten >= 0.15cm](1_32)--(1_33);
\draw[densely dotted, ->, line width=0.3mm, >=latex, shorten <= 0.2cm, shorten >= 0.15cm](1_32)--(2_33);
\draw[densely dotted, ->, line width=0.3mm, >=latex, shorten <= 0.2cm, shorten >= 0.15cm](2_32)--(1_33);
\draw[densely dotted, ->, line width=0.3mm, >=latex, shorten <= 0.2cm, shorten >= 0.15cm](2_32)--(2_33);

\draw[densely dotted, ->, line width=0.3mm, >=latex, shorten <= 0.2cm, shorten >= 0.15cm](1_43)--(1_42);
\draw[densely dotted, ->, line width=0.3mm, >=latex, shorten <= 0.2cm, shorten >= 0.15cm](1_43)--(2_42);
\draw[densely dotted, ->, line width=0.3mm, >=latex, shorten <= 0.2cm, shorten >= 0.15cm](2_43)--(1_42);
\draw[densely dotted, ->, line width=0.3mm, >=latex, shorten <= 0.2cm, shorten >= 0.15cm](2_43)--(2_42);

\draw[densely dotted, ->, line width=0.3mm, >=latex, shorten <= 0.2cm, shorten >= 0.15cm](1_23)--(1_22);
\draw[densely dotted, ->, line width=0.3mm, >=latex, shorten <= 0.2cm, shorten >= 0.15cm](1_23)--(2_22);
\draw[densely dotted, ->, line width=0.3mm, >=latex, shorten <= 0.2cm, shorten >= 0.15cm](2_23)--(1_22);
\draw[densely dotted, ->, line width=0.3mm, >=latex, shorten <= 0.2cm, shorten >= 0.15cm](2_23)--(2_22);
\draw[densely dotted, ->, line width=0.3mm, >=latex, shorten <= 0.2cm, shorten >= 0.15cm](1_12)--(1_11);
\draw[densely dotted, ->, line width=0.3mm, >=latex, shorten <= 0.2cm, shorten >= 0.15cm](1_12)--(2_11);
\draw[densely dotted, ->, line width=0.3mm, >=latex, shorten <= 0.2cm, shorten >= 0.15cm](2_12)--(1_11);
\draw[densely dotted, ->, line width=0.3mm, >=latex, shorten <= 0.2cm, shorten >= 0.15cm](2_12)--(2_11);

\draw[densely dotted, ->, line width=0.3mm, >=latex, shorten <= 0.2cm, shorten >= 0.15cm](1_32)--(1_31);
\draw[densely dotted, ->, line width=0.3mm, >=latex, shorten <= 0.2cm, shorten >= 0.15cm](1_32)--(2_31);
\draw[densely dotted, ->, line width=0.3mm, >=latex, shorten <= 0.2cm, shorten >= 0.15cm](2_32)--(1_31);
\draw[densely dotted, ->, line width=0.3mm, >=latex, shorten <= 0.2cm, shorten >= 0.15cm](2_32)--(2_31);

\draw[densely dotted, ->, line width=0.3mm, >=latex, shorten <= 0.2cm, shorten >= 0.15cm](1_41)--(1_42);
\draw[densely dotted, ->, line width=0.3mm, >=latex, shorten <= 0.2cm, shorten >= 0.15cm](1_41)--(2_42);
\draw[densely dotted, ->, line width=0.3mm, >=latex, shorten <= 0.2cm, shorten >= 0.15cm](2_41)--(1_42);
\draw[densely dotted, ->, line width=0.3mm, >=latex, shorten <= 0.2cm, shorten >= 0.15cm](2_41)--(2_42);

\draw[densely dotted, ->, line width=0.3mm, >=latex, shorten <= 0.2cm, shorten >= 0.15cm](1_21)--(1_22);
\draw[densely dotted, ->, line width=0.3mm, >=latex, shorten <= 0.2cm, shorten >= 0.15cm](1_21)--(2_22);
\draw[densely dotted, ->, line width=0.3mm, >=latex, shorten <= 0.2cm, shorten >= 0.15cm](2_21)--(1_22);
\draw[densely dotted, ->, line width=0.3mm, >=latex, shorten <= 0.2cm, shorten >= 0.15cm](2_21)--(2_22);
\draw[densely dotted, ->, line width=0.3mm, >=latex, shorten <= 0.2cm, shorten >= 0.15cm](1_13)--(6,6);
\draw[densely dotted, ->, line width=0.3mm, >=latex, shorten <= 0.2cm, shorten >= 0.15cm](1_13)--(6,5);
\draw[densely dotted, ->, line width=0.3mm, >=latex, shorten <= 0.2cm, shorten >= 0.15cm](2_13)--(6,6);
\draw[densely dotted, ->, line width=0.3mm, >=latex, shorten <= 0.2cm, shorten >= 0.15cm](2_13)--(6,5);

\draw[densely dotted, ->, line width=0.3mm, >=latex, shorten <= 0.2cm, shorten >= 0.15cm](-15,6)--(1_11);
\draw[densely dotted, ->, line width=0.3mm, >=latex, shorten <= 0.2cm, shorten >= 0.15cm](-15,6)--(2_11);
\draw[densely dotted, ->, line width=0.3mm, >=latex, shorten <= 0.2cm, shorten >= 0.15cm](-15,5)--(1_11);
\draw[densely dotted, ->, line width=0.3mm, >=latex, shorten <= 0.2cm, shorten >= 0.15cm](-15,5)--(2_11);

\draw[densely dotted, ->, line width=0.3mm, >=latex, shorten <= 0.2cm, shorten >= 0.15cm](1_33)--(7,-6);
\draw[densely dotted, ->, line width=0.3mm, >=latex, shorten <= 0.2cm, shorten >= 0.15cm](1_33)--(7,-7);
\draw[densely dotted, ->, line width=0.3mm, >=latex, shorten <= 0.2cm, shorten >= 0.15cm](2_33)--(7,-6);
\draw[densely dotted, ->, line width=0.3mm, >=latex, shorten <= 0.2cm, shorten >= 0.15cm](2_33)--(7,-7);

\draw[densely dotted, ->, line width=0.3mm, >=latex, shorten <= 0.2cm, shorten >= 0.15cm](-13,-6)--(1_31);
\draw[densely dotted, ->, line width=0.3mm, >=latex, shorten <= 0.2cm, shorten >= 0.15cm](-13,-6)--(2_31);
\draw[densely dotted, ->, line width=0.3mm, >=latex, shorten <= 0.2cm, shorten >= 0.15cm](-13,-7)--(1_31);
\draw[densely dotted, ->, line width=0.3mm, >=latex, shorten <= 0.2cm, shorten >= 0.15cm](-13,-7)--(2_31);

\draw[densely dotted, ->, line width=0.3mm, >=latex, shorten <= 0.2cm, shorten >= 0.15cm](1_41)--(-15,1);
\draw[densely dotted, ->, line width=0.3mm, >=latex, shorten <= 0.2cm, shorten >= 0.15cm](1_41)--(-15,0);
\draw[densely dotted, ->, line width=0.3mm, >=latex, shorten <= 0.2cm, shorten >= 0.15cm](2_41)--(-15,1);
\draw[densely dotted, ->, line width=0.3mm, >=latex, shorten <= 0.2cm, shorten >= 0.15cm](2_41)--(-15,0);

\draw[densely dotted, ->, line width=0.3mm, >=latex, shorten <= 0.2cm, shorten >= 0.15cm](6,1)--(1_43);
\draw[densely dotted, ->, line width=0.3mm, >=latex, shorten <= 0.2cm, shorten >= 0.15cm](6,1)--(2_43);
\draw[densely dotted, ->, line width=0.3mm, >=latex, shorten <= 0.2cm, shorten >= 0.15cm](6,0)--(1_43);
\draw[densely dotted, ->, line width=0.3mm, >=latex, shorten <= 0.2cm, shorten >= 0.15cm](6,0)--(2_43);

\draw[densely dotted, ->, line width=0.3mm, >=latex, shorten <= 0.2cm, shorten >= 0.15cm](1_21)--(-13,-1);
\draw[densely dotted, ->, line width=0.3mm, >=latex, shorten <= 0.2cm, shorten >= 0.15cm](1_21)--(-13,-2);
\draw[densely dotted, ->, line width=0.3mm, >=latex, shorten <= 0.2cm, shorten >= 0.15cm](2_21)--(-13,-1);
\draw[densely dotted, ->, line width=0.3mm, >=latex, shorten <= 0.2cm, shorten >= 0.15cm](2_21)--(-13,-2);

\draw[densely dotted, ->, line width=0.3mm, >=latex, shorten <= 0.2cm, shorten >= 0.15cm](7,-1)--(1_23);
\draw[densely dotted, ->, line width=0.3mm, >=latex, shorten <= 0.2cm, shorten >= 0.15cm](7,-1)--(2_23);
\draw[densely dotted, ->, line width=0.3mm, >=latex, shorten <= 0.2cm, shorten >= 0.15cm](7,-2)--(1_23);
\draw[densely dotted, ->, line width=0.3mm, >=latex, shorten <= 0.2cm, shorten >= 0.15cm](7,-2)--(2_23);
\end{tikzpicture}
{
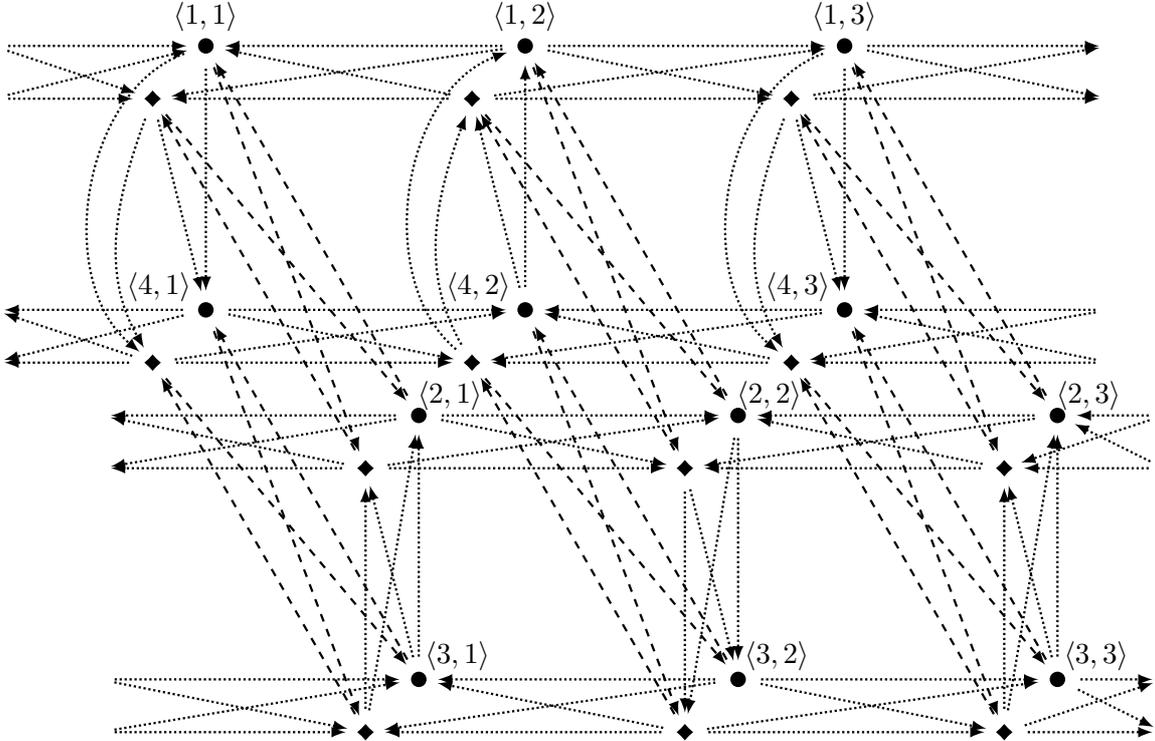
\captionof{figure}{Orientation $D$ of $(C_4\times C_3)^{(2)}$, where $d(D)=3$.}\label{figA8.5.24}}
\end{center}

\section{Concluding remarks}
In this paper, we considered primarily vertex-multiplications of Cartesian products involving trees, paths and cycles as they are some special families of graphs studied for orientations (see \cite {KKM LKT 1, KKM TEG 1, KKM TEG 2, KKM TEG 9}). We refer the interested reader to a good survey on orientations of graphs \cite{KKM TEG 10} by Koh and Tay. 
\indent\par It can be shown that $(T_2\times T_2)(s_1,s_2,\ldots, s_n)\in\mathscr{C}_0\cup\mathscr{C}_1$. We believe its characterisation likely involves notions and techniques of Extremal Set Theory such as antichains. This is akin to Proposition \ref{ppnA8.4.2} and vertex-multiplications of trees with diameter $4$ (see \cite{WHW TEG 6A}). Hence, we conclude by proposing the following problem.
\begin{prob}
Characterise the vertex-multiplications $(T_2\times T_2)(s_1,s_2,\ldots,s_n)$ that belong to $\mathscr{C}_0$.  
\end{prob}
\section*{Acknowledgements}
The first author would like to thank the National Institute of Education, Nanyang Technological University of Singapore, for the generous support of the Nanyang Technological University Research Scholarship.

\end{document}